\documentclass[10pt, a4paper, headings=standardclasses]{scrartcl}

\usepackage{packages-styles/general-math}
\usepackage{packages-styles/field-probabilitytheory}
\newcommand\Language{ngerman,english}
\usepackage{packages-styles/general-text}
\usepackage{packages-styles/general-scientific}
\usepackage{packages-styles/general-scientificpaper}
\usepackage{packages-styles/figures}
\usepackage{packages-styles/field-algorithmic}

\usepackage[frozencache, cachedir=mintedcache]{minted}
\usepackage{xcolor}
\usepackage{fancyvrb}
\usepackage{xcolor}
\usepackage{fancyvrb}

\usepackage{xcolor}
\usepackage{fancyvrb}

\makeatletter

\def\PYG@bf#1{\textbf{#1}}
\def\PYG@it#1{\textit{#1}}
\def\PYG@ul#1{\underline{#1}}
\def\PYG@tc#1{\textcolor{black}{#1}} 
\def\PYG@bc#1{#1}
\def\PYG@ff#1{#1}

\def\PYG@reset{\let\PYG@it=\relax \let\PYG@bf=\relax%
    \let\PYG@ul=\relax \let\PYG@tc=\relax%
    \let\PYG@bc=\relax \let\PYG@ff=\relax}
\def\PYG@tok#1{\csname PYG@tok@#1\endcsname}
\def\PYG@toks#1+{\ifx\relax#1\empty\else%
    \PYG@tok{#1}\expandafter\PYG@toks\fi}
\def\PYG@do#1{\PYG@bc{\PYG@tc{\PYG@ul{%
    \PYG@it{\PYG@bf{\PYG@ff{#1}}}}}}}
\def\PYG#1#2{\PYG@reset\PYG@toks#1+\relax+\PYG@do{#2}}

\expandafter\def\csname PYG@tok@w\endcsname{\def\PYG@tc##1{\textcolor[rgb]{0.73,0.73,0.73}{##1}}}
\expandafter\def\csname PYG@tok@c\endcsname{\let\PYG@it=\textit\def\PYG@tc##1{\textcolor[rgb]{0.25,0.50,0.50}{##1}}}
\expandafter\def\csname PYG@tok@cp\endcsname{\def\PYG@tc##1{\textcolor[rgb]{0.74,0.48,0.00}{##1}}}
\expandafter\def\csname PYG@tok@k\endcsname{\let\PYG@bf=\textbf\def\PYG@tc##1{\textcolor[rgb]{0.00,0.50,0.00}{##1}}}
\expandafter\def\csname PYG@tok@kp\endcsname{\def\PYG@tc##1{\textcolor[rgb]{0.00,0.50,0.00}{##1}}}
\expandafter\def\csname PYG@tok@kt\endcsname{\def\PYG@tc##1{\textcolor[rgb]{0.69,0.00,0.25}{##1}}}
\expandafter\def\csname PYG@tok@o\endcsname{\def\PYG@tc##1{\textcolor[rgb]{0.40,0.40,0.40}{##1}}}
\expandafter\def\csname PYG@tok@ow\endcsname{\let\PYG@bf=\textbf\def\PYG@tc##1{\textcolor[rgb]{0.67,0.13,1.00}{##1}}}
\expandafter\def\csname PYG@tok@nb\endcsname{\def\PYG@tc##1{\textcolor[rgb]{0.00,0.50,0.00}{##1}}}
\expandafter\def\csname PYG@tok@nf\endcsname{\def\PYG@tc##1{\textcolor[rgb]{0.00,0.00,1.00}{##1}}}
\expandafter\def\csname PYG@tok@nc\endcsname{\let\PYG@bf=\textbf\def\PYG@tc##1{\textcolor[rgb]{0.00,0.00,1.00}{##1}}}
\expandafter\def\csname PYG@tok@nn\endcsname{\let\PYG@bf=\textbf\def\PYG@tc##1{\textcolor[rgb]{0.00,0.00,1.00}{##1}}}
\expandafter\def\csname PYG@tok@s\endcsname{\def\PYG@tc##1{\textcolor[rgb]{0.73,0.13,0.13}{##1}}}
\expandafter\def\csname PYG@tok@si\endcsname{\let\PYG@bf=\textbf\def\PYG@tc##1{\textcolor[rgb]{0.73,0.40,0.53}{##1}}}
\expandafter\def\csname PYG@tok@s1\endcsname{\def\PYG@tc##1{\textcolor[rgb]{0.73,0.13,0.13}{##1}}}
\expandafter\def\csname PYG@tok@m\endcsname{\def\PYG@tc##1{\textcolor[rgb]{0.40,0.40,0.40}{##1}}}
\expandafter\def\csname PYG@tok@mi\endcsname{\def\PYG@tc##1{\textcolor[rgb]{0.40,0.40,0.40}{##1}}}
\expandafter\def\csname PYG@tok@err\endcsname{\def\PYG@bc##1{{\setlength{\fboxsep}{0pt}\fcolorbox[rgb]{1.00,0.00,0.00}{1,1,1}{\strut ##1}}}}
\expandafter\def\csname PYG@tok@kn\endcsname{\let\PYG@bf=\textbf\def\PYG@tc##1{\textcolor[rgb]{0.00,0.50,0.00}{##1}}}
\expandafter\def\csname PYG@tok@c1\endcsname{\let\PYG@it=\textit\def\PYG@tc##1{\textcolor[rgb]{0.25,0.50,0.50}{##1}}}
\makeatother

\definecolor{crimson2143940}{RGB}{214,39,40}
\definecolor{darkgray176}{RGB}{176,176,176}
\definecolor{forestgreen4416044}{RGB}{44,160,44}
\definecolor{goldenrod18818934}{RGB}{188,189,34}
\definecolor{lightgray204}{RGB}{204,204,204}
\definecolor{orchid227119194}{RGB}{227,119,194}
\definecolor{steelblue31119180}{RGB}{31,119,180}

\Crefname{property}{Property}{Properties}
\Crefname{postulate}{Postulate}{Postulates}
\setlist[description]{font=\rmfamily}

\usepackage{arxiv}

\begin{document}

  \title{Optimization of Sums of Bivariate Functions}
  \subtitle{An Introduction to Relaxation-Based Methods\\[0.06em]%
  for the Case of Finite Domains}
  \date{\today}
\author{
    {\normalfont
    \begin{tabular}[t]{c@{\extracolsep{8em}}} 
    {\bfseries Nils M{\"u}ller}  \\
    {\itshape Max Planck Institute for Informatics}  \\
    {\itshape Saarbrücken, Germany;}\\
    {\itshape Department of Mathematics} \\
     {\itshape Ruhr University Bochum}  \\
     {\itshape Bochum, Germany} \\
     \href{mailto:nmueller@mpi-inf.mpg.de} {\texttt{nmueller@mpi-inf.mpg.de}}
    \end{tabular}}}
\maketitle
\thispagestyle{empty}

\begin{abstract}
We study the optimization of functions with $n>2$ arguments that have a representation as a sum of several functions that have only $2$ of the $n$ arguments each, termed \emph{sums of bivariates}, on finite domains.

\noindent The complexity of optimizing sums of bivariates is shown to be \emph{\text{NP}-equivalent} and it is shown that there exists \emph{free lunch} in the optimization of sums of bivariates.\\

\noindent Based on measure-valued extensions of the objective function, so-called \emph{relaxations}, $\ell^2$-approximation, and entropy-regularization, we derive several tractable problem formulations solvable with linear programming, coordinate ascent as well as with closed-form solutions.

\noindent The limits of applying tractable versions of such relaxations to sums of bivariates are investigated using general results for reconstructing measures from their bivariate marginals.\\

\noindent Experiments in which the derived algorithms are applied to random functions, vertex coloring, and signal reconstruction problems provide insights into qualitatively different function classes that can be modeled as sums of bivariates.
\end{abstract}
\vspace*{5mm}

\newpage
\thispagestyle{empty}
\tableofcontents

\newpage

\section{Introduction}
In this work, we study the optimization of functions with $n>2$ arguments that have a representation as a sum of several functions that have only $2$ of the $n$ arguments each. \cref{def:biv} introduces the rigorous definition that will be used throughout the work.

\begin{definition}[Sum of Bivariates]
\label{def:biv}
    Let the function $F: \Omega \to \R$ be defined on a product space $\Omega := \Omega_1 \times \dots \times \Omega_n$ of factors $\Omega_i \subset \R$ with a finite number of elements, where $\abs{\Omega_i} = K_i$ for all $i \in \mathcal{V} \coloneqq \N_{\leq n}$ and for some $n,K_i \in \N$.

    We call the function $F$ a \textbf{sum of bivariates} if there exists an index set\\ $\mathcal{E} \subseteq \{(i,j) \in \mathcal{V} \times \mathcal{V} \mid i < j\}$
    and bivariate functions $f_{i,j}: \Omega_i \times \Omega_j \to \R, (i,j) \in \mathcal{E}$
    with
    \[
    F(x_1, \dots, x_n) = \sum_{(i,j) \in \mathcal{E}} f_{i,j}(x_i, x_j) \qquad \forall (x_1, \dots, x_n) \in \Omega \,.
\]
\end{definition}

The function model can be interpreted as indexed by a graph: The vertices index arguments and the edges index bivariate summands. If the functions $f_{i,j}, (i,j) \in \mathcal{E}$ were instead univariate, due to the additive structure, the problem would be separable into $n$ individual univariate optimization problems. In this sense, the function model at hand presents the simplest non-trivial structure of its kind that makes no further assumptions on the functions $f_{i,j}, (i,j) \in \mathcal{E}$. Yet, it will turn out that both the structure of the bivariates $f_{i,j}, (i,j) \in \mathcal{E}$ as well as the edge set $\mathcal{E}$ determine the hardness of optimizing $F$.

We restrict ourselves to finite candidate spaces $\Omega$ for the benefit of a simpler and tractable problem. Generalizations to uncountable candidate spaces are formally possible and practically useful.\\

The main focus of this work lies on so-called relaxation-based optimization methods for sums of bivariates. A \emph{relaxation} is a transformation of a function $F: \Omega \to \R$ to a function
    \[
         \mu \in \mathcal{M}_1^+(\Omega) \mapsto \int_\Omega F \, \mathrm{d}\mu \,,
    \]
    where (probability) measures $\mu \in \mathcal{M}_1^+(\Omega)$ replace the domain $\Omega$ and evaluate to the integral of $F$ with respect to $\mu$. Next to being linear functionals on a closed convex space and due to the sparse structure of sums of bivariates, relaxations of sums of bivariates will be shown to have useful properties for optimization. This perspective allows us to derive a wide range of optimization methods for sums of bivariates employing linear programming, coordinate ascent, dynamic programming, and even closed-form solution approaches. To this end, we also study the $\ell^2$-approximation and complexity theory of sums of bivariates.

    Experiments on random functions, vertex coloring, and signal reconstruction problems provide insights into qualitatively different function classes that can be modeled as sums of bivariates.

\paragraph{Motivation.} From a theoretical perspective, the motivation in studying the optimization of sums of bivariates lies in the pursuit of tractable models for global optimization. Making no assumptions on the structure of individual bivariates $f_{i,j}, (i,j) \in \mathcal{E}$ is an intriguing alternative to tractable models in the literature that often feature some sort of linear, convex, Lipschitz, or smooth structure.\\

From a practical perspective, sums of bivariates model a wide range of applications, where tractable instances often arise as inverse problems in signal reconstruction. In such problems, we often encounter a structure
    \[
        F(x_1, \dots, x_n) =  \sum_{i\in \N_{\leq n}} H_i(x_i) + \sum_{i,j \in \N_{\leq n}} G_{i,j}(x_i, x_j) \qquad \forall (x_1, \dots, x_n) \in \Omega \,,
    \]
    where $\Omega := \Omega_1 \times \dots \times \Omega_n$ models the possible reconstructions and $F: \Omega \to \R$ models their quality based on argument-wise errors $H_i: \Omega_i \to \R, i \in \N_{\leq n}$ that are based on measured data, as well as pairwise regularization terms $G_{i,j}: \Omega_i \times \Omega_j \to \R, (i,j) \in \N_{\leq n}^2$. Clearly, such problems are sums of bivariates.
    
    Sums of bivariates also appear in various other applications, such as, in the \emph{Markowitz model for portifolio selection} as the hybrid objective that is the linear combination of portfolio risk and return \cite{markowitz1952}, and as the Hamiltonian of so-called \emph{Sherrington--Kirkpatrick spin glasses} \cite{sherrington1975} as well as that of \emph{Hopfield networks} \cite{amari1972}.\\

    Due to its linearity and its commutativity with the projection operation for bivariate functions, relaxation is a particularly promising transformation for the function class at hand.
    
    The pursuit of an elementary and consistent derivation of relaxation-based results from principles of mathematical optimization, approximation, and complexity theory inspires the focus of this work.

\paragraph{Related Work.}
Relaxation techniques for sums of bivariates were originally developed by \cite{schlesinger1976}. A review of associated linear programming formulations of optimization problems encountered in this work is given by \cite{werner2007}.

State-of-the-art methods that have a structure similar to those that we derive from relaxations of sums of bivariates have been studied in \cite{kolmogorov2005,kolmogorov2014,kappes2015,tourani2018,tourani2020}. Subgradient methods that can serve as an alternative to the block coordinate ascent-based algorithms, which we introduce, have been developed by \cite{schlesinger2007,komodakis2007}. Conversely, block-coordinate ascent methods, similar to those developed in this work, have also been used to create scalable solvers for integer linear programming \cite{lange2021}. Methods similar to our entropy-regularized relaxations have been developed before \cite{globerson07a,liu2013}. 

Complexity results for various problem formulations encountered in this work are derived in \cite{li2016}. In particular, \cite{prusa2015} show that any linear program can be reduced in linear time to what we call the \emph{dual linear program} associated with a relaxation of sums of bivariates.

Extensive collections of results on the subject from a perspective of mathematical optimization as well as from a probabilistic one can be found in \cite{savchynskyy2019discrete} and \cite{wainwright2008}, respectively.\\
Very similar coordinate ascent-based optimization approaches for the Markov Random Field model have been described by \cite{ochs2024}. 

\paragraph{Outline.} In \cref{sec:approx}, we cover the basic $\ell^2$-approximation of functions by sums of bivariates. This leads us to derive an elementary characterization of parameterizations of zero and thereby the dual variables encountered later on.

We describe hard as well as non-trivial easy instances of sums of bivariates using complexity theory and dynamic programming in \cref{sec:fundamental}. An application of a no-free-lunch theorem will also be presented.

In \cref{subsec:relaxation}, we focus on the central results on relaxations of sums of bivariates and the reconstructions of measures from bivariate marginals. We extend this analysis with a regularized version of relaxation for sums of bivariates in \cref{subsec:entropyrelaxation}. We conclude our theory with a method that allows us to reconstruct solutions to the optimization of sums of bivariates from relaxation-associated problems in \cref{subsec:dual2primal}.\\

\cref{sec:algorithms} contains the derivation of algorithms based on the various problem formulations and results encountered throughout the previous sections. Pseudocode for all discussed algorithms is included in this section.

The experiments in \cref{sec:experiments} consider the minimization of sums of random bivariates, vertex coloring, and signal reconstruction problems, which constitute qualitatively different problem classes. We verify qualitative properties of the algorithms and contrast hypotheses about modeling choices.\\

A conclusion and a description of promising future work can be found in \cref{sec:conc}.

The appendix contains extensive lemmata that would otherwise inhibit readability, a computational proof, pseudocode for a state-of-the-art algorithm for the minimization of sums of bivariates, as well as important referenced results.

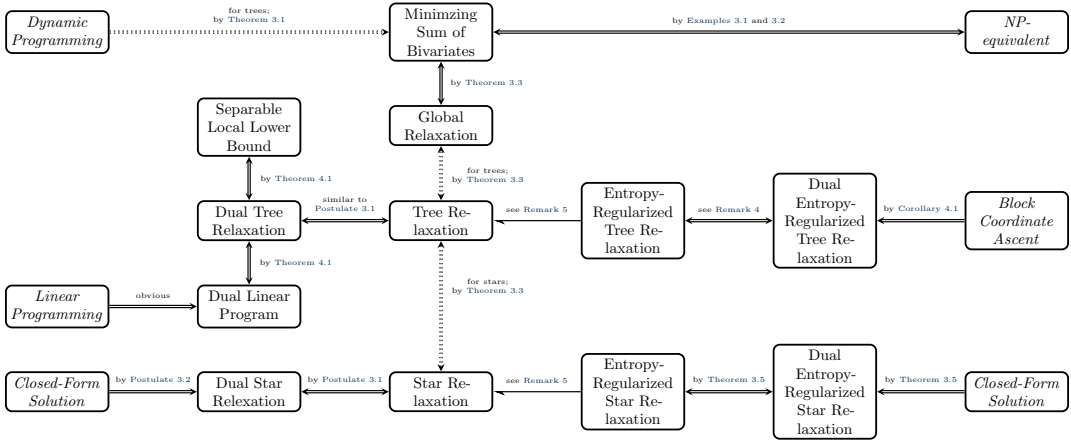
\begin{figure}[ht]
\begin{center}
\resizebox{\linewidth}{!}{\begin{tikzpicture}[auto]
     \node (center) {};
     \node[box, above=0.0cm of center] (SB) {Minimzing Sum of Bivariates};
     \node[box, below=1cm of SB] (GR) {Global Relaxation}
        edge[double, thick, <->, >=stealth] node[xshift=2cm] {\tiny{by \cref{thm:relaxation}}} (SB);
    \node[box, left=6.38cm of SB] (DP) {\emph{Dynamic Programming}}
        edge[double, thick, ->, >=stealth, dotted]  node[xshift=0.0cm, yshift=-0.2cm, label={[font=\tiny,align=center] for trees;\\by \cref{thm:dynprog}}] {} (SB);
     \node[box, below=1.25cm of GR] (TR) {Tree Relaxation}
        edge[double, thick, <->, >=stealth, dotted] node[xshift=1.2cm, yshift=-0.5cm, label={[font=\tiny,align=center] for trees;\\by \cref{thm:relaxation}}] {}
        (GR);
     \node[box, below=3cm of TR] (SR) {Star Relaxation}
        edge[double, thick, <->, >=stealth, dotted]  node[xshift=1.2cm, label={[font=\tiny,align=center] for stars;\\by \cref{thm:relaxation}}] {} (TR);
    \node[box, left=2cm of TR] (DTR) {Dual Tree Relaxation}
        edge[double, thick, <->, >=stealth]
        node[yshift=-0.175cm, label={[font=\tiny,align=center] similar to\\\cref{thm:dual}}] {} (TR);
    \node[box, left=2cm of SR] (DSR) {Dual Star Relexation}
        edge[double, thick, <->, >=stealth] node[yshift=-.2cm, label={[font=\tiny,align=center] by \cref{thm:dual}}] {} (SR);
    \node[box, above=1cm of DTR] (SLLB) {Separable Local Lower Bound}
        edge[double, thick, <->, >=stealth] node[xshift=0.95cm,yshift=-0.35cm, label={[font=\tiny,align=center] by \cref{thm:duallinprog}}] {}  (DTR);
    \node[box, below=1cm of DTR] (DLP) {Dual Linear Program}
        edge[double, thick, <->, >=stealth] node[xshift=1.2cm,yshift=-0.35cm, label={[font=\tiny,align=center] by \cref{thm:duallinprog}}] {}  (DTR);
    \node[box, left=2cm of DLP] (LP) {\emph{Linear Programming}}
        edge[double, thick, ->, >=stealth] node[yshift=-0.2cm, label={[font=\tiny,align=center] obvious}] {} (DLP);
    \node[box, left=2cm of DSR] (CF) {\emph{Closed-Form Solution}}
        edge[double, thick, ->, >=stealth] node[yshift=-0.2cm, label={[font=\tiny,align=center] by \cref{thm:nscsdv}}] {} (DSR);
    \node[box, right=2.cm of TR] (ETR) {Entropy-Regularized Tree Relaxation}
        edge[thick, ->, >={Stealth[harpoon]}]
        node[yshift=0.075cm, label={[font=\tiny,align=center] see \cref{rem:entropyapprox}}] {} (TR);
    \node[box, right=2cm of ETR] (DETR) {Dual Entropy-Regularized Tree Relaxation}
        edge[double, thick, <->, >=stealth] node[yshift=.075cm, label={[font=\tiny,align=center] see \cref{rem:reprdetr}}] {} (ETR);
    \node[box, right=2cm of DETR] (CA) {\emph{Block Coordinate Ascent}}
        edge[double, thick, ->, >=stealth] node[yshift=0.05cm, label={[font=\tiny,align=center] by \cref{cor:convergence}}] {} (DETR);
    \node[box, right=2cm of SR] (ESR) {Entropy-Regularized Star Relaxation}
        edge[thick, ->, >={Stealth[harpoon]}]
        node[yshift=0.075cm, label={[font=\tiny,align=center] see \cref{rem:entropyapprox}}] {} (SR);
    \node[box, right=2cm of ESR] (DESR) {Dual Entropy-Regularized Star Relaxation}
        edge[double, thick, <->, >=stealth] node[yshift=.05cm, label={[font=\tiny,align=center] by \cref{thm:entropydual}}] {} (ESR);
    \node[box, right=2cm of DESR] (CF2) {\emph{Closed-Form Solution}}
        edge[double, thick, ->, >=stealth] node[yshift=.05cm, label={[font=\tiny,align=center] by \cref{thm:entropydual}}] {} (DESR);
    \node[box, right=10.77cm of SB] (NPE) {\emph{\text{NP}-equivalent}}
        edge[double, thick, <->, >=stealth] node[yshift=.05cm, label={[font=\tiny,align=center] by \cref{ex:hamilton,ex:linintprog}}] {} (SB);
\end{tikzpicture}}
\end{center}%
\caption{Hierarchy of Optimization Problems: The reducibility results considered in this work. An arrow symbolizes that solving one problem yields the optimal value of another problem. Half arrows are approximations. Dotted arrows are conditional. Known generic problems are set in cursive font.}
\label{fig:hierarchy}%
\end{figure}

\cref{fig:hierarchy} describes the hierarchy of optimization problems that emerges from the reducibility results encountered in this work.

\paragraph{Acknowledgment.} This work was submitted as part of the author's Master's Thesis at Saarland University in December 2024. The author thanks Peter Ochs and Andreas Karrenbauer, who supervised this work, for many discussions and a review that was helpful in improving the manuscript.

\paragraph{Notation and Terminology.}
\begin{itemize}
    \item Given a measure space $(\Omega, \mathcal{A}, \mu)$ we write $\langle f, \mu \rangle \coloneqq \int f \, \mathrm{d}\mu$ for a $\mu$-integrable function $f: (\Omega, \mathcal{A}) \to \R$. We define the \emph{entropy} of $\mu$ as $H(\mu) \coloneqq -\E_\mu\big(\log \mu(\{\cdot\})\big)$.
    \item Given two measures $\mu, \lambda$ with $\mu = f \lambda$, we define $\langle \mu, \lambda \rangle \coloneqq \langle f, \lambda \rangle$. We define their \emph{$(\lambda,\mu)$-cross-entropy} as $H(\lambda, \mu) \coloneqq -\E_\lambda\big(\log \mu(\{\cdot\})\big)$. Their \emph{$(\lambda,\mu)$-$\mathrm{KL}$-divergence} is defined as $D_{\mathrm{KL}} \coloneqq H(\lambda, \mu) - H(\lambda)$.
    \item If $\Omega$ is a product space and $f$ is defined on a factor thereof, we often implicitly consider $f$ to be a function on $\Omega$ that is constant along the other factors.
    \item $\mathcal{M}_1^+(\Omega, \mathcal{A})$ denotes the set probability measures, $\mathcal{M}^+(\Omega, \mathcal{A})$ the set of measures, and $\mathcal{M}(\Omega, \mathcal{A})$ the set of signed measures on the measurable space $(\Omega, \mathcal{A})$.
    \item We will only consider (finite) discrete measures; therefore, we have $\abs{\Omega} \in \N$, $\mathcal{A} = \mathcal{P}(\Omega)$, $\mathcal{M}(\Omega) \cong \R^\Omega$, $\mathcal{M}^+(\Omega) \cong \R_{\geq 0}^\Omega$ and $\mathcal{M}_1^+(\Omega) \cong \Delta_{\abs{\Omega}}$, where $\Delta_{\abs{\Omega}}$ denotes the $\abs{\Omega}$-dimensional simplex. Thus, we also adopt the notation $\mu(x) \coloneqq \mu(\{x\})$, where $x \in \Omega$ and $\mu \in \mathcal{M}(\Omega)$.
    \item Given a product space $\Omega = \Omega_1 \times \dots \times \Omega_n, n\in \N$, we define the projection $\pi_i^{\Omega}: \Omega \to \Omega_i, x \mapsto x_i$. Often we will drop the superscript and write $\pi_i \coloneqq \pi_i^\Omega$, if the domain is clear by context. In the same setting, where additionally $\Omega_i = \N_{\leq m}$, we often denote $f \in \R^{\Omega_i}$ by $\begin{bmatrix}f(1) & \dots & f(m)\end{bmatrix}_i$.
    \item A function $F: \Omega \to \R$ is called \emph{separable} if there exist $f_i \in \R^{\Omega_i}, i \in \N_{\leq n}$ with $F \equiv f_1 + \dots + f_n$.
    \item We often replace universally quantified variables by $\cdot$ and do not explicitly denote the quantifier.
    \item Given a graph $\mathcal{G}=(\mathcal{V}, \mathcal{E})$ with vertices $\mathcal{V} = \N_{\leq n}$ and edges $\mathcal{E} \subseteq \{(i,j) \in \mathcal{V} \times \mathcal{V} \mid i < j \}$, we define $\mathcal{N}(i) \coloneqq \{ (i,j) \in \mathcal{E}\}$ and $\widetilde{\mathcal{N}}(i) \coloneqq \{ j \in \mathcal{V} \mid (i,j) \in \mathcal{E} \lor (j,i) \in \mathcal{E} \}, i \in \mathcal{V}$. In oriented trees, we will assume that a parent vertex is in the left factor of an edge and a child vertex to be in the right factor of an edge.
    \item For a set $\Omega$ with $\abs{\Omega} \in \N$ and $f \in \R^{\Omega}$, we define as $\mathrm{lse}_\varepsilon(f) \coloneqq \varepsilon\log \sum_{x \in \Omega}\exp\big(f(x)/\varepsilon\big)$ the \emph{$\varepsilon$-$\mathrm{LogSumExp}$} for all $\varepsilon > 0$. We will also use the notation $\mathrm{lse}^x_\varepsilon\big(f(x)\big) \coloneqq \mathrm{lse}_\varepsilon(f)$ if $f$ has further parameters.
\end{itemize}

\clearpage
\section{Approximation}
\label{sec:approx}%

Initially, we would like to better understand the parameterization that is implicit to the defining functional equation of sums of bivariates in \cref{def:biv}. This will yield not only a characterization of the defining functional equation but also a way to efficiently approximate functions by sums of bivariates, i.e.~to project onto the space of sums of bivariates.\\

\begin{figure}[ht]%
        \centering
        \includegraphics[width=1.\textwidth]{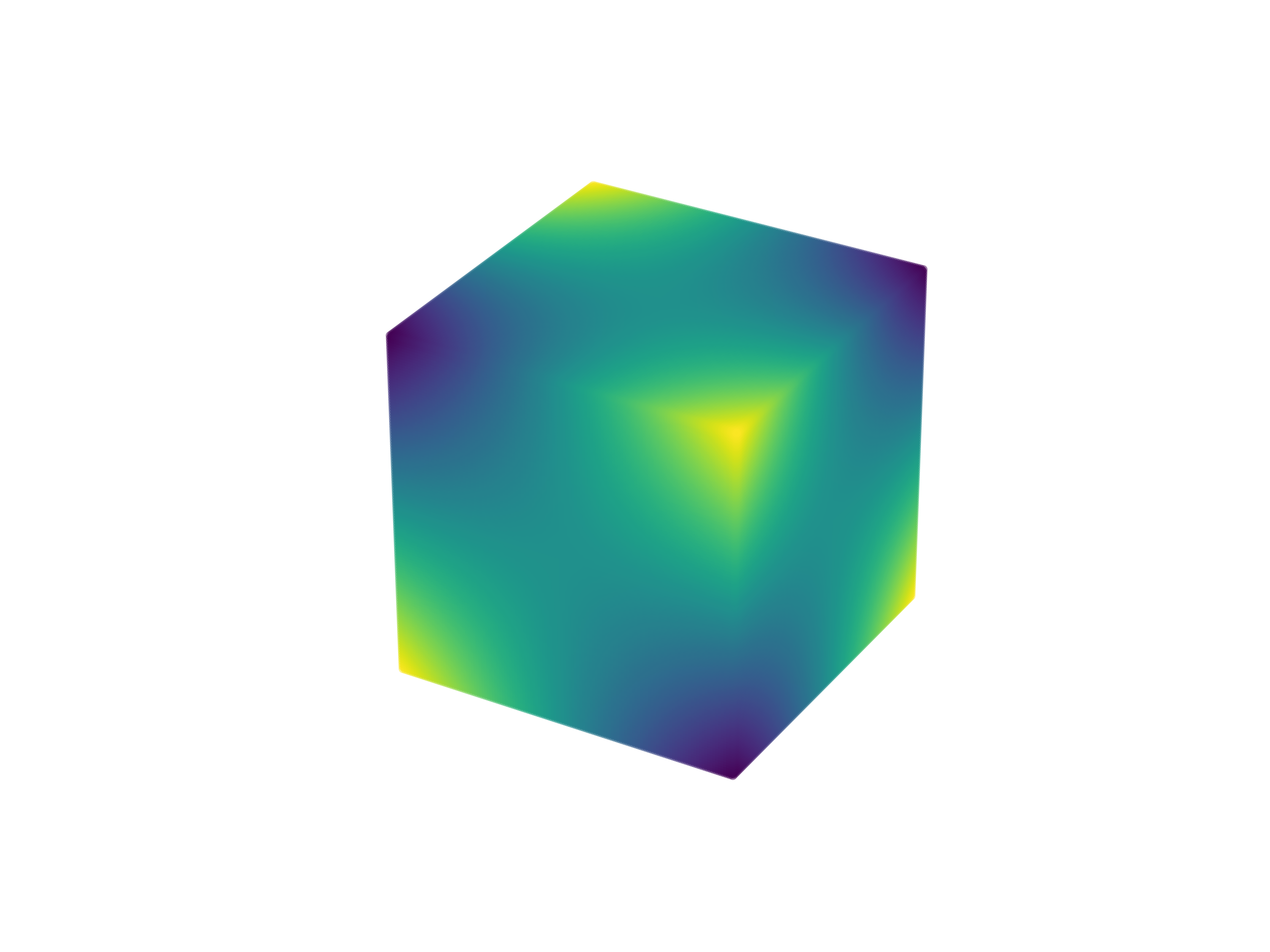}
        \vspace*{-2cm}
        \caption{Visualization of a function $L^2$-orthogonal to the sums of bivariates on $\Omega = [0,1]^3$.}
\label{fig:biv-orth}
\end{figure}%

For a given function $G: \N_{\leq m}^n \to \R$, we call the function that sums $G$ over all values that all but its $k$-th and $\ell$-th argument can take, where $k<\ell$, the \emph{$(k,\ell)$-marginal} of $G$, i.e.
\[
    a,b \in \N_{\leq m} \longmapsto \sum_{\substack{y_i \in \N_{\leq m}\\
                    i \in \N_{\leq n}\\
                    i \neq k,l}} G\big(y_1, \dots, y_{k-1}, a, y_{k+1}, \dots, y_{\ell-1}, b, y_{\ell+1}, \dots, y_{n}\big) \,.
\]
The first theorem provides a necessary and sufficient condition for the approximation of a function $G$ by sums of bivariates in terms of all $(k, \ell)$-marginals. By that, the theorem also provides a characterization of the defining functional when replacing $G$ by a function that is postulated to be a sum of bivariates. The theorem features a more compact notation of the $(k,\ell)$-marginal.

Insights into conditions and limits for bivariates to be the marginals of a possibly unknown function $G$ will be covered in \cref{thm:rfm,ex:counterrecon}.

\begin{theorem}[Approximation]
\label{thm:approx}
        Given $G: \N_{\leq m}^n \to \R$, $m \in \N$, a sum of bivariates $F: \N_{\leq m}^n \to \R$, where $F \equiv \sum_{(i,j) \in \mathcal{E}} f_{i,j}$, $\mathcal{E} = \{(i,j) \in \mathcal{V} \times \mathcal{V} \mid i < j \}$, and $\mathcal{V} = \N_{\leq n}$, is a best $\ell^2$-approximation of $G$ if and only if
    \begingroup%
    \allowdisplaybreaks%
        \begin{align*}
            \sum_{y\in\N_{\leq m}^{n-2}} G\big(z^{k,l}(a,b,y) \big) &= \sum_{y\in\N_{\leq m}^{n-2}} F\big(z^{k,l}(a,b,y) \big) \qquad \forall a,b \in \N_{\leq m} \quad \forall (k, \ell) \in \mathcal{E}
            \\
            &=
            m^{n-2} f_{k,\ell}(a,b)
         + m^{n-3} \bigg(\big(\textstyle\sum_{j<k} \overline{f_{j,k}}^1(a) + \overline{f_{j, \ell}}^1(b) \big)
         \\
         &
         \qquad\qquad\qquad\qquad\qquad\qquad
         + \big(\textstyle\sum_{k<j<\ell} \overline{f_{k,j}}^2(a) + \overline{f_{j, \ell}}^1(b)\big)
         \\
         &
         \qquad\qquad\qquad\qquad\qquad\qquad
         + \big(\textstyle\sum_{\ell<j} \overline{f_{k,j}}^2(a) + \overline{f_{\ell,j}}^2(b) \big)\bigg)
             \\
         &
         \quad\quad
         +
         m^{n-4} \!\!\!\sum_{i,j\notin \{k,\ell\}} \overline{\overline{f_{i,j}}}\,,
    \tag*{\small{(this eq.~only holds if $n\geq 4$)}}
        \end{align*}%
        \endgroup%
        where%
        \begin{gather*}%
            z^{k,\ell}(a,b,y) \coloneqq (y_1, \dots, y_{k-1}, a, y_k, \dots, y_{\ell-2}, b, y_{\ell-1}, \dots, y_{n-2}) \,,
            \\
            \overline{f_{i,j}}^1(s) \coloneqq \sum_{x \in \N_{\leq m}} f_{i,j}(x,s) \,,\,\,
            \overline{f_{i,j}}^2(s) \coloneqq \sum_{x \in \N_{\leq m}} f_{i,j}(s,x)\,,\,\, \forall s \in \{a,b\}  \,,\,\, \text{and}
            \\
            \overline{\overline{f_{i,j}}} \coloneqq \sum_{x,s \in \N_{\leq m}} f_{i,j}(x,s) \,,\,\, \forall (i,j) \in \mathcal{E} \,.
        \end{gather*}
    \end{theorem}%
\begin{proof}
    We want to find a bivariate function $F \equiv \sum_{(i,j) \in \mathcal{E}} f_{i,j}$,  where $\mathcal{E} = \{(i,j) \in \mathcal{V} \times \mathcal{V} \mid i < j \}$, and $\mathcal{V} = \N_{\leq m}$ that minimizes
    \[
        \norm{G - F}_{\ell^2} = \sqrt{\sum_{x\in\N_{\leq m}^n} \big( G(x) - \sum_{(i,j)\in\mathcal{E}} f_{i,j}(x_i, x_j) \big)^2} \,.
    \]
    Due to convexity, the first-order condition
    \begingroup
    \allowdisplaybreaks
    \begin{align*}
        &0 \overset{!}{=} \frac{\partial}{\partial f_{k,\ell}(a,b)} \sum_{x\in\N_{\leq m}^n} \big( G(x) - \sum_{(i,j)\in\mathcal{E}} f_{i,j}(x_i, x_j) \big)^2
         \qquad \forall a,b \in \N_{\leq m} \quad \forall (k,\ell) \in \mathcal{E}
    \tag*{\small{(first-order condition)}}
         \\
         &\phantom{0}=
         \frac{\partial}{\partial f_{k,\ell}(a,b)} \!\!\sum_{y_k,y_\ell\in \N_{\leq m}} \sum_{y\in\N_{\leq m}^{n-2}}\!\!\! \big( G(z^{k,\ell}(y_k,y_\ell,y)) - \!\!\!\!\sum_{(i,j)\in\mathcal{E}}\!\! f_{i,j}(z^{k,\ell}_i(y_k,y_\ell,y), z^{k,\ell}_j(y_k,y_\ell,y)) \big)^2
    \tag*{\small{(reorder summation; definition $z^{k,\ell}$)}}
         \\
         &\phantom{0}=
         -2\sum_{y\in\N_{\leq m}^{n-2}}  G(z^{k,\ell}(a,b,y)) - \sum_{(i,j)\in\mathcal{E}} f_{i,j}(z^{k,\ell}_i(a,b,y), z^{k,\ell}_j(a,b,y))
    \tag*{\small{(differentiation)}}
         \\
         \iff\!\!\!\!\!\!&
         \sum_{y\in\N_{\leq m}^{n-2}}  G\big(z^{k,\ell}(a,b,y)\big) = \!\!\!\sum_{y\in\N_{\leq m}^{n-2}} \sum_{(i,j)\in\mathcal{E}} f_{i,j}\big(z^{k,\ell}_i(a,b,y), z^{k,\ell}_j(a,b,y)\big)
    \tag*{\small{(multiply by $-1/2$; add marginalized sum of bivariates)}}
         \\
         &
         \phantom{\sum_{y\in\N_{\leq m}^{n-2}}  G(z^{k,\ell}(a,b,y))} = \!\!\!\sum_{y\in\N_{\leq m}^{n-2}} F\big(z^{k,\ell}(a,b,y)\big)
    \tag*{\small{(definition of $F$)}}
         \\
         &
         \phantom{\sum_{y\in\N_{\leq m}^{n-2}}  G(z^{k,\ell}(a,b,y))} =
         m^{n-2} f_{k,\ell}(a,b)
         + m^{n-3} \bigg(\big(\textstyle\sum_{j<k} \overline{f_{j,k}}^1(a) + \overline{f_{j, \ell}}^1(b) \big)
         \\
         &
         \phantom{\sum_{y\in\N_{\leq m}^{n-2}}  G(z^{k,\ell}(a,b,y))}\qquad\qquad\qquad\qquad\qquad\qquad
         + \big(\textstyle\sum_{k<j<\ell} \overline{f_{k,j}}^2(a) + \overline{f_{j, \ell}}^1(b)\big)
         \\
         &
         \phantom{\sum_{y\in\N_{\leq m}^{n-2}}  G(z^{k,\ell}(a,b,y))}\qquad\qquad\qquad\qquad\qquad\qquad
         + \big(\textstyle\sum_{\ell<j} \overline{f_{k,j}}^2(a) + \overline{f_{\ell,j}}^2(b) \big)\bigg)
             \\
         &
         \phantom{\sum_{y\in\N_{\leq m}^{n-2}}  G(z^{k,\ell}(a,b,y))}\quad\quad
         +
         m^{n-4} \!\!\!\sum_{i,j\notin \{k,\ell\}} \overline{\overline{f_{i,j}}}
    \tag*{\small{(reorder summation based on occurrences of $a$ and $b$; definitions $\overline{f_{\cdot,\cdot}}^\cdot$, $\overline{\overline{f_{\cdot,\cdot}}}$)}}
    \end{align*}
    \endgroup
    is necessary and sufficient.
\end{proof}
\cref{fig:biv-orth} shows the residual of (numerically) approximating $(x,y,z) \in [0,1]^3 \mapsto xyz$ by sums of bivariates. The resulting function has a $\ell^2$-best approximant by sums of bivariates that is the function $0$, therefore, it is $\ell^2$-orthogonal to the sums of bivariates. There are, however, still multiple sets of bivariates that sum to $0$ and, therefore, parameterize the best approximant.\\

Thus, \cref{cor:dv} characterizes the parameterizations of $0$ by sums of bivariates. As we will see in \cref{sec:optimization,sec:algorithms}, the result is central to solvers that use dual formulations and---sometimes without loss---transforms the search space of the minimization of sums of bivariates into one with comparatively low dimensions.

A relevant feature of the following result is the representation not only of the constant sum of bivariates $F$ as a sum of univariates, but, in particular, a representation of the individual bivariates $f_{i,j}$ as a sum of univariates.
\begin{corollary}[Dual~Variables]
    \label{cor:dv}
        In the setting of \cref{thm:approx}, we have
        \begingroup
        \allowdisplaybreaks
        \begin{align*}
            & \sum_{y\in\N_{\leq m}^{n-2}} F\big(z^{k,l}(a,b,y) \big) = 0 \qquad \forall a,b \in \{1, \dots, m\} \quad \forall (k, \ell) \in \mathcal{E}
        \tag*{\small{(marginal null constraint)}}
            \\
            &\quad\iff
            F \equiv 0
        \tag*{\small{(global null constraint)}}
            \\
            &\quad\iff
            \forall (i,j) \in \mathcal{E}: \exists \rho_{i,j}, \rho_{j, i} \in \R^{\N_{\leq m}}: \bigg(f_{i,j} \equiv \rho_{i,j}(\cdot_i) + \rho_{j, i}(\cdot_j)\bigg)
            \\
            &\qquad\qquad\qquad
            \land \bigg(\forall i \in \mathcal{V}: \exists \rho_i \in \R: \sum_{j \in \mathcal{V}} \rho_{i,j} \equiv \rho_i\bigg)
            \\
            &\qquad\qquad\qquad
            \land \bigg(\sum_{i\in \mathcal{V}} \rho_i = 0\bigg) \,.
        \tag*{\small{(dual null constraint)}}
        \end{align*}
        \endgroup
    \end{corollary}
    \begin{remark}
        Clearly, a similar result holds for constant sums of bivariates.
    \end{remark}
    \begin{proof}
        We prove the statement as a Ringschluss.
        \begin{description}
        \item[\normalfont{\underline{\enquote{marginal null constraint $\Rightarrow$ dual null constraint}}:}]
        Consider, that for all $(k,\ell) \in \mathcal{E}$ and for all $a,b \in \N_{\leq m}$, we have
        \begingroup
        \allowdisplaybreaks
        \begin{align*}
            &\qquad0
            =
            \sum_{y\in\N_{\leq m}^{n-2}} F\big(z^{k,l}(a,b,y) \big)
        \tag*{\small{(marginal null constraint)}}
            \\
            &\phantom{\qquad0}=
            \sum_{y\in\N_{\leq m}^{n-2}} \sum_{(i,j) \in \mathcal{E}} f_{i,j}\big(z^{k,\ell}_i(a,b,y), z^{k,\ell}_j(a,b,y) \big)
        \tag*{\small{(representation of sum of bivariates)}}
            \\
            &\iff
            f_{k,\ell}(a,b)
            =- \frac{1}{m^{n-2}}\sum_{y\in\N_{\leq m}^{n-2}} \sum_{(i,j) \in \mathcal{E}\setminus\{(k,\ell)\}} f_{i,j}\big(z^{k,\ell}_i(a,b,y), z^{k,\ell}_j(a,b,y) \big) \,.
        \tag*{\small{(multiplication by $-1/m^{n-2}$; addition of $f_{k,\ell}(a,b)$)}}
        \end{align*}
        \endgroup
        Since the right-hand side of the previous equation is a sum of univariate functions in either $a$ or $b$, or of constants, for all $k, \ell \in \mathcal{E}$, we get that
        \[
            \exists \rho_{k, \ell} , \rho_{\ell, k} \in \R^{\N_{\leq m}}: f_{k, \ell} \equiv \rho_{k, \ell} + \rho_{\ell, k} \,.
        \]
        Further, for all $a\in \N_{\leq n}$ and for all $k \in \mathcal{V}$, we have
        \begingroup
        \allowdisplaybreaks
        \begin{align*}
            &\qquad0
            =
            \sum_{b \in \N_{\leq m}}\sum_{y\in\N_{\leq m}^{n-2}} F\big(z^{k,l}(a,b,y) \big)
        \tag*{\small{(sum of marginal null constraint)}}
            \\
            &\phantom{\qquad0}=
            \sum_{b\in \N_{\leq m}}\sum_{y\in\N_{\leq m}^{n-2}} \sum_{(i,j) \in \mathcal{E}} f_{i,j}\big(z^{k,\ell}_i(a,b,y), z^{k,\ell}_j(a,b,y) \big)
        \tag*{\small{(representation of sum of bivariates)}}
            \\
            &\phantom{\qquad0}=
            \sum_{b\in \N_{\leq m}}\sum_{y\in\N_{\leq m}^{n-2}} \sum_{(i,j) \in \mathcal{E}} \rho_{i,j}\big(z^{k,\ell}_i(a,b,y)\big)+\rho_{j,i}\big(z^{k,\ell}_j(a,b,y) \big)
        \tag*{\small{(previous result)}}
            \\
            &\phantom{\qquad0}=
            \sum_{b\in \N_{\leq m}}\sum_{y\in\N_{\leq m}^{n-2}} \sum_{i \in \mathcal{V}} \sum_{j \in \mathcal{V}} \rho_{i,j}\big(z^{k,\ell}_i(a,b,y)\big)
        \tag*{\small{(reorder summation)}}
            \\
            &\iff
            \sum_{j\in \mathcal{V}} \rho_{k, j}(a)
            =
            -\frac{1}{m^{n-1}}\sum_{b\in \N_{\leq m}}\sum_{y\in\N_{\leq m}^{n-2}} \sum_{i \in \mathcal{V}\setminus\{k\}} \sum_{j \in \mathcal{V}} \rho_{i,j}\big(z^{k,\ell}_i(a,b,y)\big) \,.
        \tag*{\small{(multiplication by $-1/m^{n-1}$; addition of $\sum_{j \in \mathcal{V}} \rho_{k, j}(a)$)}}
        \end{align*}
        \endgroup
        Since the right-hand side of the previous equation is a constant, and by a similar argument for the variable $b$, we get for all $k \in \mathcal{V}$ that
        \[
            \exists \rho_k \in \R: \sum_{j \in \mathcal{V}} \rho_{k,j} \equiv \rho_k \,.
        \]%
        Finally, we get the last claim for an arbitrary choice of $a,b \in \N_{\leq m}$ and $k, \ell \in \mathcal{E}$ by
    \begingroup
    \allowdisplaybreaks
        \begin{align*}
            &\qquad0
            =
            \sum_{y\in\N_{\leq m}^{n-2}} F\big(z^{k,l}(a,b,y) \big)
        \tag*{\small{(marginal null constraint)}}
            \\
            &\phantom{\qquad0}=
            \sum_{y\in\N_{\leq m}^{n-2}} \sum_{(i,j) \in \mathcal{E}} \rho_{i,j}\big(z^{k,\ell}_i(a,b,y)\big)+\rho_{j,i}\big(z^{k,\ell}_j(a,b,y) \big)
        \tag*{\small{(previous result)}}
            \\
            &\phantom{\qquad0}=
            \sum_{y\in\N_{\leq m}^{n-2}} \sum_{i \in \mathcal{V}} \bigg(\sum_{j \in \mathcal{V}} \rho_{i,j}\big(z^{k,\ell}_i(a,b,y)\big)\bigg)
        \tag*{\small{(reorder summation)}}
            \\
            &\phantom{\qquad0}\equiv
            \sum_{y\in\N_{\leq m}^{n-2}} \sum_{i \in \mathcal{V}} \rho_i
        \tag*{\small{(previous result)}}
            \\
            &\phantom{\qquad0}=
            \cancel{m^{n-2}} \sum_{i \in \mathcal{V}} \rho_i \,.
        \tag*{\small{(constant summation)}}
        \end{align*}
    \endgroup
        \item[\normalfont{\underline{\enquote{dual null constraint $\Rightarrow$ global null constraint}}:}] We have for all $x_1, \dots, x_n \in \N_{\leq m}$ that
        \begin{align*}
            F(x_1, \dots, x_n)
            &=
            \sum_{(i,j) \in \mathcal{E}} f_{i,j}(x_i, x_j)
        \tag*{\small{(representation of sum of bivariates)}}
            \\
            &=
            \sum_{(i,j) \in \mathcal{E}} \rho_{i,j}(x_i) + \rho_{j, i}(x_j)
        \tag*{\small{(first property of dual null constraint)}}
            \\
            &=
            \sum_{i \in \mathcal{V}} \bigg(\sum_{j \in \mathcal{V}} \rho_{i,j}(x_i) \bigg)
        \tag*{\small{(reorder summation)}}
            \\
            &\equiv
            \sum_{i \in \mathcal{V}} \rho_i
        \tag*{\small{(second property of dual null constraint)}}
            \\
            &=
            0 \,.
        \tag*{\small{(third property of dual null constraint)}}
        \end{align*}
        \item[\normalfont{\underline{\enquote{global null constraint $\Rightarrow$ marginal null constraint}}:}] Clear.\qedhere
        \end{description}
    \end{proof}
Since the main goal of this work is the analysis of relaxation-based solvers, further approximation-based insights into sums of bivariates will be deferred.
    
\clearpage%
\section{Optimization}
\label{sec:optimization}

\subsection{Fundamental Perspective}
\label{sec:fundamental}%

We now attempt to classify as well as understand some hard- and easy-to-optimize instances of sums of bivariates. Such insight may help in modeling, managing our expectations for the optimization methods we design, and in interpreting their performance.\\

The Hamiltonian cycle problem asks to determine whether for a given graph there exists a sequence of successively adjacent vertices that are nonrepetitive. The Hamiltonian cycle problem is among the hardest decision problems for which a solution is polynomially sized and efficiently verifiable.

In \cref{ex:hamilton}, we recite the reduction of the Hamiltonian cycle problem to checking whether the minimal value of a sum of bivariates is $0$. By that we show that efficiently minimizing sums of bivariates is at least as hard as such hardest decision problems.

\begin{example}[{The Hamiltonian Cycle Problem as a Sum of Bivariates \cite[p.~13~Sec.~1.3]{savchynskyy2019discrete}}]
\label{ex:hamilton}%
        Let $(\mathcal{V}', \mathcal{E}')$ be a graph for which we want prove existence of a Hamiltonian cycle.\\
        Define $\mathcal{V} \coloneqq \N_{\leq n}, n = \abs{\mathcal{V}'}$, $\mathcal{E} = \{(i,j) \in \mathcal{V} \times \mathcal{V} \mid i< j\}$, $\Omega_i \coloneqq \mathcal{V}', i \in \N_{\leq n}$, $\Omega \coloneqq \Omega_1\times \cdots \times \Omega_n$, and for all $(i,j) \in \mathcal{E}$ for all $(x_i, x_j) \in \Omega_i \times \Omega_j$, that
        \[
            f_{i,j}(x_i, x_j) \coloneqq \begin{cases}
                0 & \text{if } x_i \neq x_j \land \big( i+1 =j \implies (x_i, x_j) \in \mathcal{E}' \big) \bigg)
                \\
                1 & \text{else.}
            \end{cases}
        \]%
        Therefore,
        \begin{align*}
            &\min_{x \in \Omega} \textstyle\sum_{(i,j) \in \mathcal{E}} f_{i,j}(x_i,x_j) = 0
            \\
            &\qquad\qquad\iff
            \exists x \in \Omega: \bigg( \forall (i,j) \in \mathcal{E}: x_i \neq x_j \land \big( i+1 =j \implies (x_i, x_j) \in \mathcal{E}' \big) \bigg)
            \\
            &\qquad\qquad\iff
            \text{$(\mathcal{V}', \mathcal{E}')$ has a Hamiltonian cycle.}
        \end{align*}%
    \end{example}
    The problem of deciding whether the graph $\mathcal{G}' = (\mathcal{V}', \mathcal{E}')$ has a Hamiltonian cycle is \text{NP}-complete. Therefore, since the Hamiltonian cycle problem can be reduced in polynomial time to checking the minimal value of a sum of bivariates, we know that minimization of the sum of bivariates is \text{NP}-hard.\\

    Next, we show in \cref{ex:linintprog} that the minimization of sums of bivariates can be efficiently reduced to \emph{integer linear programming}. This implies that minimizing rationally-valued sums of bivariates is also \text{NP}-easy. In conjunction with the result of \text{NP}-hardness by \cref{ex:hamilton}, we therefore know that the minimization of rationally-valued sums of bivariates is \text{NP}-equivalent. This means that there exists an exact and efficient algorithm for the minimization of rationally-valued sums of bivariates if and only if there is an efficient algorithm for the hardest decision problems for which solutions are polynomially sized and efficiently verifiable, i.e.~if $\text{P} = \text{NP}$ \cite[p.~368~Prop.~15.35]{vygen2006}.

    \begin{example}[Sums of Bivariates as Integer Linear Programs]%
    \label{ex:linintprog}%
        In the setting of \cref{def:biv}, let $F \equiv \sum_{(i,j) \in \mathcal{E}} f_{i,j}: \Omega_1 \times \dots \times \Omega_n \to \mathbb{Q}$ encode a sum of bivariates.
        
        Let, w.l.o.g., $\abs{\Omega_1} =  \dots = \abs{\Omega_i} = m \in \N$. Consider the binary variables of the integer linear program to be
        \[
            x_{i,j} \in \{0,1\}^{m \times m} \qquad \forall (i,j) \in \mathcal{E} \,,
        \]
        and the constraints to be
        \[
            \begin{cases}
            \sum_{s,t=1}^m x_{i,j; s,t} = 1 & \forall (i,j) \in \mathcal{E}
            \\
            \sum_{s=1}^m x_{i,j; s, t} = \sum_{s=1}^m x_{k,j; s, t} & \forall t \in \N_{\leq m} \,\, \forall (i,j), (k,j) \in \mathcal{E}
            \\
            \sum_{s=1}^m x_{i,j; t, s} = \sum_{s=1}^m x_{k,j; t, s} & \forall t \in \N_{\leq m} \,\, \forall (i,j), (i,k) \in \mathcal{E}
            \\
            \sum_{s=1}^m x_{i,j; s, t} = \sum_{s=1}^m x_{j,k; t, s} & \forall t \in \N_{\leq m} \,\, \forall (i,j), (j,k) \in \mathcal{E} \,.
            \end{cases}
        \]
        The objective to be minimized can be modeled as
        \[
            \sum_{(i,j) \in \mathcal{E}} \sum_{s,t=1}^m f_{i,j}(s,t) \cdot x_{i,j; s,t} \,,
        \]
        which leaves us with a linear integer program that has a polynomial size w.r.t.~to the size of the sum of bivariates. Clearly, we obtain an optimal solution $y^* \in \Omega$ to the minimization of sums of bivariates given an optimal solution to the integer linear program $x^* \in \big(\{0,1\}^{m \times m}\big)^{\abs{\mathcal{E}}}$ by picking
        \[
            y^*_i \in \argmax_{s \in \N_{\leq m}} \sum_{t=1}^m x^*_{i,j;s,t} \quad \text{and} \quad y^*_j \in \argmax_{j \in \N_{\leq m}} \sum_{s=1}^m x^*_{i,j;s,t} \qquad \forall (i,j) \in \mathcal{E} \,,
        \]
        which is well-defined by our constraints.
    \end{example}
In contrast to the described hard instances, we now introduce relatively easy-to-minimize instances of sums of bivariates, which are widely known. Interestingly, this class of functions is not characterized by requirements on the structure of the bivariates but on the structure that encodes the dependencies, that is, the graph $\mathcal{G}=(\mathcal{V}, \mathcal{E})$ that indexes the sum of bivariates $F \equiv \sum_{(i,j) \in \mathcal{E}} f_{i,j}$ in \cref{def:biv}. Such easy instances are, in fact, characterized by tree-structured index graphs $\mathcal{G}$.\\

The following \cref{thm:dynprog} is constructive in the sense that it describes an algorithm to efficiently solve such easy instances. The minimization of sums of bivariates with these tree-structured indices can be interpreted as a dynamic program that sequentially allows the reduction of the problem to a smaller instance.

\begin{theorem}[Dynamic~Programming~for~Tree-Indexed~Instances]
\label{thm:dynprog}
    Let $F \equiv \sum_{(i,j) \in \mathcal{E}} f_{i, j}$ be a sum of bivariates, where $\mathcal{G} = (\mathcal{V}, \mathcal{E})$ is a tree, and let $(i_t, j_t) \in \mathcal{E}, j_t\coloneqq\abs{\mathcal{V}}+1-t, t \in \N_{\leq \abs{\mathcal{V}}-1}$ be an edge counting, such that
    \begin{itemize}
        \item $\mathcal{E}_{t} \coloneqq \mathcal{E}_{t-1} \setminus \{(i_t, j_t)\}$, where $\mathcal{E}_0 \coloneqq \mathcal{E}\cup\{(*, i_{\abs{\mathcal{V}}-1})\}$ and $*$ is a placeholder symbol,
        \item $\mathcal{V}_t \coloneqq \mathcal{V}_{t-1} \setminus\{j_t\}$, where $\mathcal{V}_0 \coloneqq \mathcal{V}$,
        \item $j_t$ is a leaf of $\mathcal{G}_{t-1} \coloneqq (\mathcal{V}_{t-1}, \mathcal{E}_{t-1})$,
        \item define $f_{0;*, i_{\abs{\mathcal{V}}-1}} :\equiv 0$, $f_{0; k, \ell} \coloneqq f_{k,\ell}$ and
        \[
            f_{t;k,\ell} \coloneqq \begin{cases}
                f_{t-1;k,\ell} + 
                \min_{x_{j_t}} f_{t-1;i_t, j_t}(\cdot, x_{j_t}) & \text{if } \ell = i_t
                \\
                f_{t-1;k,\ell} & \text{else,}
            \end{cases}
        \]
        for all $(k,\ell) \in \mathcal{E}_t$ for all $t \in \N_{\leq \abs{\mathcal{V}}-1}$, and
        \item $F_t :\equiv \sum_{(i,j)\in \mathcal{E}_t} f_{t;i,j}$ for all $t\in \N_{[0,\abs{V}-1]}$.
    \end{itemize}
\newpage
    Then,
    \begin{itemize}
        \item $\min F = \min F_t$ for all $t \in \N_{\leq \abs{\mathcal{V}}-1}$, and
        \item $x^* \in \argmin_{x \in \Omega} F(x)$ if $x_t^* \in \argmin_{x_t \in \Omega_t} F_{\abs{V}-t}(x_1^*, \dots, x_{t-1}^*, x_t)$ for all $t \in \N_{\leq \abs{\mathcal{V}}}$.
    \end{itemize}
\end{theorem}
\begin{proof} The first result can be inductively concluded from the following equation. We have for all $t \in \N_{\leq \abs{\mathcal{V}}-1}$ that
    \begingroup
    \allowdisplaybreaks
    \begin{align*}
        F_t
        &=
        \sum_{(i,j)\in \mathcal{E}_t} f_{t;i,j}
    \tag*{\small{(definition of $F_t$)}}
        \\
        &= \bigg(\sum_{(i,j)\in \mathcal{E}_t}f_{t-1;i,j}\bigg)+\sum_{(i,j)\in \mathcal{E}_t} f_{t;i,j} - f_{t-1;i,j}
    \tag*{\small{(reorder summation; partition of zero)}}
        \\
        &=
        \bigg(\sum_{(i,j)\in \mathcal{E}_t}f_{t-1;i,j}\bigg) + \min_{x_{j_t}} f_{t-1;i_t, j_t}(\cdot, x_{j_t})
    \tag*{\small{(definition of $f_{t;\cdot,\cdot}$)}}
        \\
        &=
        \min_{x_{j_t}} \sum_{(i,j)\in \mathcal{E}_{t-1}}f_{t-1;i,j}
    \tag*{\small{(definition of $\mathcal{E}_t$; $j_t$ is leaf of $\mathcal{G}_{t-1}$)}}
        \\
        &=
        \min_{x_{j_t}} F_{t-1} \,.
    \tag*{\small{(definition of $F_{t-1}$)}}
    \end{align*}
    \endgroup
    The second result is implied by
    \[
        (x_1^*, \dots,x_t^*) \in \argmin_{x_1, \dots,x_t} F_{\abs{\mathcal{V}}-t}(x_1^*, \dots,x_{t-1}^*,x_t) \quad \forall t \in \N_{\leq \abs{\mathcal{V}}} \,,
    \]
    which is also true inductively by
    \begin{align*}
        F_{\abs{\mathcal{V}}-t}(x_1^*, \dots, x^*_{t})
        &=
        \min_{x_t} F_{\abs{\mathcal{V}}-t}(x_1^*, \dots,x_{t-1}^*,x_{t})
    \tag*{\small{(definition of $x^*_t$)}}
        \\
        &=
        F_{\abs{\mathcal{V}}+1-t}(x_1^*, \dots,x_{t-1}^*)
    \tag*{\small{(previous equation)}}
        \\
        &=
        \min F_{\abs{\mathcal{V}}+1-t}
    \tag*{\small{(induction hypothesis)}}
        \\
        &=
        \min F_{\abs{\mathcal{V}}-t} \,,
    \tag*{\small{(first result)}}
    \end{align*}
    where $t \in \N_{[2,\abs{\mathcal{V}}]}$ and $x_1^* \in \argmin_{x_1} F_{\abs{\mathcal{V}}-1}(x_1)$ by definition.
\end{proof}

In \cref{ex:dynprog}, we present a small but detailed example of a sum of bivariates that we minimize using the method described in \cref{thm:dynprog}.

\begin{example}[Dynamic~Programming~for~Tree-Indexed~Instances]
\label{ex:dynprog}
In the setting of \cref{def:biv}, consider
\begin{gather*}
    \mathcal{V} \coloneqq \{1,2,3,4,5\}\,,\quad 
    \mathcal{E} \coloneqq \{(1,2),(2,3),(3,4),(3,5)\} \,, 
    \\
    \Omega_1 = \dots = \Omega_5 = \{1,2,3\} \,,
\end{gather*}%
and a sum of bivariates
\[
    F \equiv \sum_{(i,j)\in \mathcal{E}} f_{i,j} \,,
\]
where
\begin{align*}%
        f_{1,2} &\coloneqq \begin{bmatrix}
        3 & 8 & 0 \\
        6 & 8 & 9 \\
        6 & 1 & 9
        \end{bmatrix}
        \in \R^{\Omega_1\times \Omega_2}\,,&
        f_{2,3} &\coloneqq 
        \begin{bmatrix}
        8 & 4 & 2 \\
        6 & 9 & 9 \\
        6 & 4 & 2
        \end{bmatrix}
        \in \R^{\Omega_2\times \Omega_3}
        \,,
        \\
        f_{3,4} &\coloneqq 
        \begin{bmatrix}
        7 & 0 & 0 \\
        2 & 2 & 1 \\
        4 & 2 & 6
        \end{bmatrix}
        \in \R^{\Omega_3\times \Omega_4}
        \,,&
        f_{3,5} &\coloneqq 
        \begin{bmatrix}
        8 & 9 & 6 \\
        0 & 2 & 7 \\
        9 & 3 & 3
        \end{bmatrix}
        \in \R^{\Omega_3\times \Omega_5}\,.
    \end{align*}
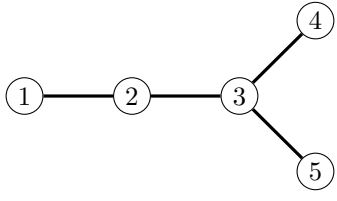
\begin{figure}[ht]%
\begin{center}%
\begin{tikzpicture}[scale=1, every node/.style={circle, draw, fill=white, inner sep=2pt}]%
    \node (1) at (-1.414, 1) {1};
    \node (2) at (0, 1) {2};
    \node (3) at (1.414, 1) {3};
    \node (4) at (2.414, 2) {4};
    \node (5) at (2.414, 0) {5};
    
    \draw[very thick] (1) -- (2);
    \draw[very thick] (2) -- (3);
    \draw[very thick] (3) -- (4);
    \draw[very thick] (3) -- (5);
\end{tikzpicture}\end{center}%
    \caption{Visualization of the variable dependencies of the sum of bivariates of \cref{ex:dynprog}.}%
\end{figure}%
Then, in the setting of \cref{thm:dynprog}, we have
\[
    j_1 = 5\,,\,\, j_2 = 4\,,\,\, j_3 = 3\,,\,\,\text{and } j_4 = 2 \,,
\]
as well as
\begingroup
\allowdisplaybreaks
\begin{align*}
    F_0 &\equiv f_{1,2} + f_{2,3} + f_{3,4} + f_{3,5}    
    \\
    F_1 &\equiv f_{1,2} + f_{2,3} + f_{3,4} + \min_{x_5 \in \Omega_5} f_{3,5}(\cdot,x_5)
    \\
    &\equiv
    f_{1,2} + f_{2,3} + f_{3,4} + \begin{bmatrix}
        6 & 0 &3
    \end{bmatrix}_3
    \\
    F_2&\equiv f_{1,2} + f_{2,3} + \big(\min_{x_4 \in \Omega_4} f_{3,4}(\cdot,x_4)\big) + \big(\min_{x_5 \in \Omega_5} f_{3,5}(\cdot,x_5)\big)
    \\
    &\equiv
    f_{1,2} + f_{2,3} + \begin{bmatrix}
        0 & 1 &2
    \end{bmatrix}_3 + \begin{bmatrix}
        6 & 0 &3
    \end{bmatrix}_3
    \\
    F_3 &\equiv f_{1,2} + \bigg(\min_{x_3 \in \Omega_3} f_{2,3}(\cdot, x_3) + \big(\min_{x_4 \in \Omega_4} f_{3,4}(x_3,x_4)\big) + \big(\min_{x_5 \in \Omega_5} f_{3,5}(x_3,x_5)\big)\bigg)
    \\
    &\equiv
    f_{1,2} + \bigg( \min_{x_3 \in \Omega_3} f_{2,3}(\cdot, x_3) + \begin{bmatrix}
        6 & 1 &5
    \end{bmatrix}_3(x_3) \bigg)
    \\
    &\equiv
    f_{1,2} + \min_{x_3 \in \Omega_3} \begin{bmatrix}
        14 & 5 & 7
        \\
        12 & 10 & 14
        \\
        12 & 5 & 7
    \end{bmatrix}_{2,3}\!\!\!\!(\cdot, x_3)
    \\
    &\equiv
    f_{1,2} +\begin{bmatrix}
        5& 10 & 5
    \end{bmatrix}_2
    \\
    F_4 &\equiv  \min_{x_2 \in \Omega_2} f_{1,2}(\cdot, x_2) + \bigg(\min_{x_3 \in \Omega_3} f_{2,3}(x_2, x_3) + \big(\min_{x_4 \in \Omega_4} f_{3,4}(x_3,x_4)\big) + \big(\min_{x_5 \in \Omega_5} f_{3,5}(x_3,x_5)\big)\bigg)
    \\
    &\equiv
    \min_{x_2 \in \Omega_2}
    \begin{bmatrix}
        8 & 18 & 5
        \\
        11 & 18 & 14
        \\
        11 & 11 & 14
    \end{bmatrix}_{1,2}\!\!\!\!(\cdot, x_2)
    \\
    &\equiv
    \begin{bmatrix}
        5& 11 & 11
    \end{bmatrix}_1 \,.
\end{align*}
\endgroup
By \cref{thm:dynprog}, we know $\min F = \min F_4 = 5$ and for
\begingroup
\allowdisplaybreaks
\begin{align*}
    &x_1^* \in \argmin_{x_1\in\Omega_1} F_4(x_1) = \{ 1\}\,,
    \\
    &x_2^* \in \argmin_{x_2\in\Omega_2} F_3(x_1^*, x_2) =  \argmin_{x_2\in\Omega_2} \begin{bmatrix}
        8& 18 & 5
    \end{bmatrix}_2 = \{3\}\,,
    \\
    &x_3^* \in \argmin_{x_3\in\Omega_3} F_2(x_1^*, x_2^*, x_3) =  \argmin_{x_3\in\Omega_3}
    \begin{bmatrix}
        12& 5 & 7
    \end{bmatrix}_3 = \{2\}\,,
    \\
    &x_4^* \in \argmin_{x_4\in\Omega_4} F_1(x_1^*, x_2^*, x_3^*, x_4) =  \argmin_{x_4\in\Omega_4}
    \begin{bmatrix}
        6& 6 & 5
    \end{bmatrix}_4 = \{3\}\,, \text{and}
    \\
    &x_5^* \in \argmin_{x_4\in\Omega_4} F(x_1^*, x_2^*, x_3^*, x_4^*, x_5) =  \argmin_{x_5\in\Omega_5}
    \begin{bmatrix}
        6& 8 & 13
    \end{bmatrix}_5 = \{1\} \,,
\end{align*}
\endgroup
we have $(x_1^*, \dots, x_5^*) \in \argmin F$.
\end{example}
We identified some hard and easy instances of sums of bivariates. However, it does not have to be our goal to solve the hardest problems optimally. We may be happy to solve other instances overproportionally efficient with a particular optimization algorithm.

Related to this perspective are results that constrain the overperformance of some algorithms over others when applied to a class of objective functions, which are called \emph{no-free-lunch results}. Such a result could, in principle, give us a formal reason to abandon the goal of minimizing general sums of bivariates---or at the very least a reason to constrain to particular subclass of sums of bivariates. The following \cref{thm:freelunch}, presented rather informally, gives us a characterization of classes of objectives for which all arguably reasonable algorithms have the same performance distribution with respect to all arguably reasonable performance measures.

\begin{theorem}[{(No-)Free-Lunch~Theorem~\cite[p.~3~Thm.~2]{igel2004} \& \cite{schumacher2001}}]
\label{thm:freelunch}
    Let $\Omega$, $Y\subset \R$ be finite sets.
    If and only if $\mathcal{F} \subseteq Y^{\Omega}$ is closed under composition with permutations of $\Omega$, then for
        \begin{itemize}
            \item any two algorithms $a,b$, mapping (candidate, objective value)-sequences to non-repeating successive candidate-objective value points,
            \item any budget $k \in \N_{\leq \abs{\Omega}}$, and
            \item any performance measure, mapping objective value-sequences to real numbers,
        \end{itemize}
        the distribution of performances of length-$k$ sequences', generated by $a$ and $b$, when applied to $\mathcal{F}$, is equivalent.
    \end{theorem}

We conclude this section with the construction of a counterexample. That is, we show in \cref{ex:freelunch} that sums of bivariates are not closed under permutation of their domain. This means that we find an instance of the sums of bivariates $F \equiv \sum_{(i,j) \in \mathcal{E}} f_{i,j}: \Omega \to \R$, such that when we compose it with a permutation of $\Omega$, we obtain a function that does not have a representation as a sum of bivariates.

Admittedly, the negation of \cref{thm:freelunch} is rather weak: It just implies that a performance measure exists, such that there exist two algorithms with non-equal performance measures when applied to sums of bivariates. However, this fact may still serve as a sanity check, as we would hope this is true if the endeavor of finding good optimizers for arbitrary sums of bivariates is meant to be worthwhile.

\begin{example}[Sums of Bivariate Functions Not Closed Under Permutation]
\label{ex:freelunch}
        Consider
        \[
            \begin{cases}
                F: \{0,1\}^6 \longrightarrow \Z_{[-64, 64]}
                \\
                x_1x_2\dots x_6 \longmapsto \sum_{\substack{i,j\in \N_{\leq 6}\\i < j}} \chi_{x_i = x_j = 1}(x_i, x_j)
            \end{cases}
        \]%
        and the permutation of $\N_{[0, 63]}$ in cycle notation given by
        \begin{align*}
            p = &(0, 8, 18, 11, 54, 41, 28, 26, 55, 59, 48, 40, 60, 24, 47, 12, 33, 63, 13, 22, 25, 16, 23, 32, 7, 36, 21,\\
            &\phantom{(} 6, 1, 52, 44, 50, 42, 17, 10, 53, 37, 14, 39, 9, 58, 46, 38, 51, 5, 27, 56, 31, 15, 49, 35, 61, 45, 3,\\
            &\phantom{(} 30, 19, 57, 34, 4, 43, 2, 62, 20, 29) \,.
        \end{align*}%
        Then, let $\iota: \{0,1\}^6 \to \N_{[0, 63]}$ be the decimal coding, then $\Tilde{p} = \iota^{-1}\circ p \circ\iota$ is a permutation of $\{0,1\}^6$, and $F \circ \Tilde{p}$ not a sum of bivariate functions. The claim can be verified by running the program in \cref{code:freelunch}.
    \end{example}

\subsection{Relaxation Principle}
\label{subsec:relaxation}%

The focus of this work are relaxation-based optimization algorithms. Relaxation is a transformation of an optimization problem $\min_{x \in \Omega} F(x)$ that generalizes the concept of a candidate $x \in \Omega$ to probability measures $\mathcal{M}_1^+(\Omega)$, as well as the objective value of a candidate from $F(x)$ to
\[
    \langle F, \mu \rangle \coloneqq \int_\Omega F \, \mathrm{d} \mu \,,
\]
and is applicable to integrable objectives. Relaxations are used in other domains of (applied) mathematics to generalize solution concepts, such as in game theory and PDE-analysis. Although the transformation is often not of practical interest, it associates the original optimization problem with a unique linear optimization problem on a convex and often compact domain.

For our purposes, relaxations serve us in generating insight into our problem as well as to derive practically useful problem formulations.\\

The following \cref{thm:relaxation} recalls that one can generally associate a relaxed optimization problem with the problem of minimizing an objective function, such that it has the same minimal value and any optimal measure evaluates to zero on the set of suboptimal candidates of the original problem. This is, in particular, true for sums of bivariates.

The deep insight of \cref{thm:relaxation} is that given a candidate $\mu \in \mathcal{M}_1^+(\Omega)$ of the relaxed (global) problem, it is sufficient to know for all $(i,j) \in \mathcal{E}$ its bivariate marginal measure
\begin{align*}
    \mu_{i,j} &\coloneqq \!\!\mu\big( \Omega_1 \!\times\! \cdots \!\times\! \Omega_{i-1} \!\times\! (\cdot) \!\times\! \Omega_{i+1} \!\times\! \cdots \!\times\! \Omega_{j-1} \!\times\! (\cdot) \!\times\! \Omega_{j+1} \!\times\! \cdots\!\times\! \Omega_n \big)
    \\
    &\,=
    \mu \circ \pi_{i,j}^{-1} \in \mathcal{M}_1^+(\Omega_i \times \Omega_j) \,,
\end{align*}
to determine its objective value. Here, the function $\pi_{i,j}: \Omega \to \Omega_i \times \Omega_j$ denotes the projection to the coordinates $i,j$. However, only if the graph that indexes the sum of bivariates is a tree, we associate optimal solutions of the relaxation with optimal solutions of the minimization of the sums of bivariates based solely on their bivariate marginal measures. By this we obtain a formulation for the optimization of tree-structured sums of bivariates that is already efficiently solvable and is an alternative to the dynamic programming approach presented in \cref{thm:dynprog}.

Additionally, in case the sum of bivariates is indexed by a star graph, one obtains an even simpler linear programming formulation, which we will later learn to admit a closed-form solution.

\begin{theorem}[Relaxation]
    If a sum of bivariates $F\equiv\sum_{(i,j)\in\mathcal{E}}f_{i,j}: \Omega \to \R$ has a unique minimum, then
    \begingroup
    \allowdisplaybreaks
    \begin{align*}
         &\argminmin_{x \in \Omega} F(x_1, \dots, x_n)
         \\
         &\quad=
         \argminmin_{\mu \in \mathcal{M}_1^+(\Omega)}\langle F, \mu \rangle \,.
    \tag*{\small{(global relaxation)}}
\intertext{Further, if $\mathcal{G}= (\mathcal{V}, \mathcal{E})$ is a tree, then}
        &\quad=
        \begin{cases}
        \displaystyle\argminmin_{\substack{\mu_{\cdot, \cdot} \in \mathcal{M}^+(\Omega_\cdot \times \Omega_\cdot)\\ \mu_\cdot \in \mathcal{M}_1^+(\Omega_\cdot)}} \sum_{(i,j) \in \mathcal{E}} \big\langle f_{i,j} , \mu_{i,j}\big\rangle
        \\
        \,\,\text{s.t.}\quad
        \mu_i = \mu_{i,j} \circ \pi_i^{-1} \,,\quad \forall (i,j) \in \mathcal{E}
        \\
        \phantom{\,\,\text{s.t.}\quad}
        \mu_j = \mu_{i,j} \circ \pi_j^{-1} \,,\quad \forall (i,j) \in \mathcal{E} \,.
        \end{cases}
    \tag*{\small{(tree relaxation)}}
\intertext{Further, if $\mathcal{G} = (\mathcal{V}, \mathcal{E})$ is a star graph with root $i^* \in \mathcal{V}$ and $\mathcal{E}=\mathcal{N}(i^*)$, then}
        &\quad=
        \begin{cases}
        \displaystyle\argminmin_{\substack{\mu_{i^*, \cdot} \in \mathcal{M}^+(\Omega_{i^*} \times \Omega_\cdot)\\ \mu_{i^*} \in \mathcal{M}_1^+(\Omega_{i^*})}} \sum_{(i^*,j) \in \mathcal{N}(i^*)} \big\langle f_{i^*,j} , \mu_{i^*,j}\big\rangle
        \\
        \,\,\text{s.t.}\quad
        \mu_{i^*} = \mu_{i^*,j} \circ \pi_{i^*}^{-1} \,,\quad \forall (i^*,j) \in \mathcal{N}(i^*) \,.
        \end{cases}
    \tag*{\small{(star relaxation)}}
    \end{align*}
    \endgroup
\label{thm:relaxation}
\end{theorem}
\begin{proof}\leavevmode
    \begin{enumerate}[label=\roman*.]
        \item The first claim---the equivalence of the minimization of the sum of bivariates and its global relaxation---is a generic result, which can be found in \cref{thm:globalrelaxation}.
        \item The second claim---the equivalence of the minimization of the global relaxation and the tree relaxation of the sum of bivariates---can be proven by induction over trees $\mathcal{G}$:\\
        The fact is trivially true if $\mathcal{G}$ has only one edge, as then the global measure and its $2$-marginals are equivalent.\\
        We prove now that if the second claim holds for all trees with $s \in \N$ edges, then it must also hold for all trees with $s+1$ edges. Let $\mathcal{G}$ now have $s+1$ edges, i.e.~$\abs{\mathcal{E}}=s+1$, and w.l.o.g.~let $(i^{j^*}, j^*)$ denote an edge, $j^*$ be a leaf, then
        \begingroup
        \allowdisplaybreaks
        \begin{align*}
            \argminmin_{\mu \in \mathcal{M}_1^+(\Omega)}\langle F, \mu \rangle
            &=
            \argminmin_{\mu \in \mathcal{M}_1^+(\Omega)} \bigg\langle f_{k^*, i^{j^*}} + f_{i^{j^*}, j^*}+ \sum_{\substack{(i,j) \in \mathcal{E}\\(i,j) \neq (i^{j^*}, j^*)\\ (i,j) \neq (k^*, i^{j^*})}} f_{i,j}, \mu \bigg\rangle
        \tag*{\small{(as $G$ is a tree, there exists some $(k^*, i^{j^*}) \in \mathcal{E}$)}}
            \\
            &=
            \argminmin_{\mu \in \mathcal{M}_1^+(\Omega)} \bigg\langle f_{k^*, i^{j^*}} + m_{i^{j^*}, j^*} + \sum_{\substack{(i,j) \in \mathcal{E}\\(i,j) \neq (i^{j^*}, j^*)\\ (i,j) \neq (k^*, i^{j^*})}} f_{i,j}, \mu \bigg\rangle
        \tag*{\small{(by $m_{i^{j^*}, j^*}(\cdot) \coloneqq \min_{x_{j^*} \in \Omega_{j^*}} f_{i^{j^*}, j^*}(\cdot, x_{j^*})$; $j^*$ is a leaf)}}
        \\
        &=
        \begin{cases}
        \displaystyle\argminmin_{\substack{\mu_{\cdot, \cdot} \in \mathcal{M}^+(\Omega_\cdot \times \Omega_\cdot)\\ \mu_\cdot \in \mathcal{M}_1^+(\Omega_\cdot)}} \big\langle f_{k^*, i^{j^*}} + m_{i^{j^*}, j^*}, \mu_{k^*,i^{j^*}} \big\rangle +\sum_{\substack{(i,j) \in \mathcal{E}\\(i,j) \neq (i^{j^*}, j^*)\\ (i,j) \neq (k^*, i^{j^*})}} \big\langle f_{i,j} , \mu_{i,j}\big\rangle
        \\
        \,\,\text{s.t.}\quad
        \mu_i = \mu_{i,j} \circ \pi_i^{-1} \,,\quad \forall (i,j) \in \mathcal{E}\setminus\{(i^{j^*}, j^*)\}
        \\
        \phantom{\,\,\text{s.t.}\quad}
        \mu_j = \mu_{i,j} \circ \pi_j^{-1} \,,\quad \forall (i,j) \in \mathcal{E}\setminus\{(i^{j^*}, j^*)\}
        \end{cases}
        \tag*{\small{(induction hypothesis; $m_{i^{j^*}, j^*}: \Omega_{i^{j^*}} \to \R$)}}
        \\
        &=
        \begin{cases}
        \displaystyle\argminmin_{\substack{\mu_{\cdot, \cdot} \in \mathcal{M}^+(\Omega_\cdot \times \Omega_\cdot)\\ \mu_\cdot \in \mathcal{M}_1^+(\Omega_\cdot)}} \big\langle m_{i^{j^*}, j^*}, \mu_{i^{j^*}} \big\rangle +\sum_{\substack{(i,j) \in \mathcal{E}\\(i,j) \neq (i^{j^*}, j^*)}} \big\langle f_{i,j} , \mu_{i,j}\big\rangle
        \\
        \,\,\text{s.t.}\quad
        \mu_i = \mu_{i,j} \circ \pi_i^{-1} \,,\quad \forall (i,j) \in \mathcal{E}\setminus\{(i^{j^*}, j^*)\}
        \\
        \phantom{\,\,\text{s.t.}\quad}
        \mu_j = \mu_{i,j} \circ \pi_j^{-1} \,,\quad \forall (i,j) \in \mathcal{E}\setminus\{(i^{j^*}, j^*)\}
        \end{cases}
        \tag*{\small{($m_{i^{j^*}, j^*}: \Omega_{i^{j^*}} \to \R$)}}
        \\
        &=
        \begin{cases}
        \displaystyle\argminmin_{\substack{\mu_{\cdot, \cdot} \in \mathcal{M}^+(\Omega_\cdot \times \Omega_\cdot)\\ \mu_\cdot \in \mathcal{M}_1^+(\Omega_\cdot)}} \big\langle m_{i^{j^*}, j^*}, \mu_{i^{j^*},j^*} \big\rangle +\sum_{\substack{(i,j) \in \mathcal{E}\\(i,j) \neq (i^{j^*}, j^*)}} \big\langle f_{i,j} , \mu_{i,j}\big\rangle
        \\
        \,\,\text{s.t.}\quad
        \mu_i = \mu_{i,j} \circ \pi_i^{-1} \,,\quad \forall (i,j) \in \mathcal{E}
        \\
        \phantom{\,\,\text{s.t.}\quad}
        \mu_j = \mu_{i,j} \circ \pi_j^{-1} \,,\quad \forall (i,j) \in \mathcal{E}
        \end{cases}
        \tag*{\small{(introduce $\mu_{j^*}\in\mathcal{M}_1^+(\Omega_{j^*})$, $\mu_{i^{j^*}, j^*} \in \mathcal{M}^+(\Omega_{i^{j^*}}\times\Omega_{j^*}): \mu_{i^{j^*}} = \mu_{i^{j^*}, j^*} \circ \pi_{i^{j^*}}^{-1}$, $\mu_{j^*} = \mu_{i^{j^*}, j^*} \circ \pi_{j^*}^{-1}$)}}
        \\
        &=
        \begin{cases}
        \displaystyle\argminmin_{\substack{\mu_{\cdot, \cdot} \in \mathcal{M}^+(\Omega_\cdot \times \Omega_\cdot)\\ \mu_\cdot \in \mathcal{M}_1^+(\Omega_\cdot)}}  \sum_{(i,j) \in \mathcal{E}} \big\langle f_{i,j} , \mu_{i,j}\big\rangle
        \\
        \,\,\text{s.t.}\quad
        \mu_i = \mu_{i,j} \circ \pi_i^{-1} \,,\quad \forall (i,j) \in \mathcal{E}
        \\
        \phantom{\,\,\text{s.t.}\quad}
        \mu_j = \mu_{i,j} \circ \pi_j^{-1} \,,\quad \forall (i,j) \in \mathcal{E} \,,
        \end{cases}
        \tag*{\small{(by $m_{i^{j^*}, j^*}(\cdot) \coloneqq \min_{x_{j^*} \in \Omega_{j^*}} f_{i^{j^*}, j^*}(\cdot, x_{j^*})$; $j^*$ is a leaf)}}
        \end{align*}
        \endgroup
        which proves the desired result by induction.

    \item The third claim---the equivalence of the tree relaxation and the star relaxation of the sum of bivariates---immediately follows from the fact that the unary measures of leafs are not part of the objective function and do not impose any non-trivial constraints.%
    \qedhere
    \end{enumerate}
\end{proof}

Given we have marginally-coupled bivariate measures $\mu_{i,j}$ on $\Omega_i \times \Omega_j, (i,j) \in \mathcal{E}$, such as the candidates of the tree relaxation in \cref{thm:relaxation}, we can ask the question of whether there exists a (global) measure on $\Omega$ that has these bivariate measures as marginals.
\cref{thm:rfm} proves existence of such a global measure for tree-indexed sets of bivariate measures and constructs one such global measure with maximum entropy.
\begin{theorem}[Reconstruction~From~Marginals]
\label{thm:rfm}
    If $\mathcal{G}= (\mathcal{V}, \mathcal{E})$ is an oriented tree with root $i^* \in \mathcal{V}$ and given
    \begin{align*}
    \begin{cases}
        \mu_{i,j} \in \mathcal{M}^+(\Omega_i \times \Omega_j)\,,\quad \forall (i,j) \in \mathcal{E}
        \\
        \mu_i \in \mathcal{M}^+_1(\Omega_i) \,,\quad \forall i \in \mathcal{V}
        \\
        \,\,\text{s.t.}\quad
        \mu_i = \mu_{i,j} \circ \pi_i^{-1} \,,\quad \forall (i,j) \in \mathcal{E}
        \\
        \phantom{\,\,\text{s.t.}\quad}
        \mu_j = \mu_{i,j} \circ \pi_j^{-1} \,,\quad \forall (i,j) \in \mathcal{E} \,,
    \end{cases}
    \end{align*}
    then the measure $\mu \in \mathcal{M}_1^+(\Omega)$ with
    \[
        \mu(x) \coloneqq \begin{cases}
            \mu_{i^*}(x_{i^*}) \prod_{(i,j) \in \mathcal{E}}  \frac{\mu_{i,j}(x_i, x_j)}{\mu_i(x_i)} & \text{if } \mu_i(x_i) >0 \quad \forall i \in \mathcal{V}
            \\
            0 & \text{else}
        \end{cases}
        \,,\,\,\forall (x_1, \dots, x_n) \in \Omega
    \]
    is a maximal-entropy measure, such that it has has marginals $\big(\mu_{i,j}\big),(i,j) \in \mathcal{E}$, i.e.
    \[
        \mu_{i,j} = \mu \circ \pi_{i,j}^{-1} \qquad \text{for all } (i,j) \in \mathcal{E} \,.
    \]
\end{theorem}
\begin{remark}
    The definition of $\mu$ in \cref{thm:rfm} does not depend on the choice of the root of a tree. See \cref{cor:rmer}.
\end{remark}
\begin{proof}\leavevmode
We assume w.l.o.g.~that for all $x \in \Omega$ we have $\mu(x) > 0$, otherwise restrict $\Omega$.
    \begin{enumerate}[label=\roman*.]
        \item\label{rfm:proof:i} First, we prove that $\mu$ is a probability measure. To this end, let $(i^{j^*}, j^*)$ denote an edge and $j^*$ be a leaf of $\mathcal{G}$. We have
        \begingroup
        \allowdisplaybreaks
        \begin{align*}
            \mu(\Omega)
            &=
            \sum_{x \in \Omega} \mu(x)
        \tag*{\small{(additivity of measures)}}
            \\
            &=
            \sum_{\substack{x_i \in \Omega_i\\i\neq j^*}} \sum_{x_{j^*} \in \Omega_{j^*}} \mu(x_1, \dots, x_n)
        \tag*{\small{(reorder summation)}}
            \\
            &=
            \sum_{\substack{x_i \in \Omega_i\\i\neq j^*}} \sum_{x_{j^*} \in \Omega_{j^*}} \mu_{i^*}(x_{i^*}) \prod_{(i,j) \in \mathcal{E}}  \frac{\mu_{i,j}(x_i, x_j)}{\mu_i(x_i)}
        \tag*{\small{(definition of $\mu$)}}
            \\
            &=
            \sum_{\substack{x_i \in \Omega_i\\i\neq j^*}} \mu_{i^*}(x_{i^*}) \bigg(\prod_{\substack{(i,j) \in \mathcal{E}\\j \neq j^*}}  \frac{\mu_{i,j}(x_i, x_j)}{\mu_i(x_i)} \bigg) \sum_{x_{j^*} \in \Omega_{j^*}} \frac{\mu_{i^{j^*},j^*}(x_{i^{j^*}}, x_{j^*})}{\mu_{i^{j^*}}(x_{i^{j^*}})}
        \tag*{\small{(distributivity; $j^*$ is leaf)}}
            \\
            &=
            \sum_{\substack{x_i \in \Omega_i\\i\neq j^*}} \mu_{i^*}(x_{i^*}) \bigg(\prod_{\substack{(i,j) \in \mathcal{E}\\j \neq j^*}}  \frac{\mu_{i,j}(x_i, x_j)}{\mu_i(x_i)} \bigg)
        \tag*{\small{(by $\mu_{i^{j^*}} = \mu_{i^{j^*}, j^*} \circ \pi_{i^{j^*}}^{-1}$)}}
            \\
            &=
            \sum_{x_{i^*} \in \Omega_{i^*}} \mu_{i^*}(x_{i^*})
        \tag*{\small{(repeated application of the last few steps for remaining leafs)}}
            \\
            &=
            1 \,.
        \tag*{\small{(by $\mu_{i^*} \in \mathcal{M}_1^*(\Omega_{i^*})$)}}
        \end{align*}
        \endgroup
    \item Similarly, we prove that $\mu \circ \pi_{i,j}^{-1} = \mu_{i,j}$ for all $(i,j) \in \mathcal{E}$. W.l.o.g., we assume that for an arbitrary edge $(i^{j^*}, j^*) \in \mathcal{E}$ the vertex $j^*$ is a leaf and that there are no other leafs---otherwise can repeatedly marginalize w.r.t.~leafs like in \cref{rfm:proof:i}. Denote with $(i^*, j^{i^*}) \in \mathcal{E}$ the edge that includes the root. We have for all $x_{i^{j^*}} \in \Omega_{i^{j^*}}, x_{j^*} \in \Omega_{j^*}$ that
    \begingroup
    \allowdisplaybreaks
    \begin{align*}
        &\mu \circ \pi_{i^{j^*},j^*}^{-1}(x_{i^{j^*}}, x_{j^*} )
    \\
        &\qquad=
        \sum_{\substack{x_i \in \Omega_i\\i \neq i^{j^*},j^*}}
        \mu_{i^*}(x_{i^*}) \prod_{(i,j) \in \mathcal{E}}  \frac{\mu_{i,j}(x_i, x_j)}{\mu_i(x_i)}
    \tag*{\small{(explicit summation; definition of $\mu$)}}
    \\
        &\qquad=
        \sum_{\substack{x_i \in \Omega_i\\i \neq i^*,i^{j^*},j^*}}
        \bigg(\sum_{\substack{x_{i^*} \in \Omega_{i^*}\\\mu_{i^*}(x_{i^*})>0}}
        \cancel{\mu_{i^*}(x_{i^*})} \frac{\mu_{i^*,j^{i^*}}(x_{i^*}, x_{j^{i^*}})}{\cancel{\mu_{i^*}(x_{i^*})}} \bigg)\prod_{\substack{(i,j) \in \mathcal{E}\\i\neq i^*}}  \frac{\mu_{i,j}(x_i, x_j)}{\mu_i(x_i)}
    \tag*{\small{(distributivity; $i^*$ is a root)}}
    \\
        &\qquad=
        \sum_{\substack{x_i \in \Omega_i\\i \neq i^*,i^{j^*},j^*}}
        \mu_{j^{i^*}}(x_{j^{i^*}}) \prod_{\substack{(i,j) \in \mathcal{E}\\i\neq i^*}}  \frac{\mu_{i,j}(x_i, x_j)}{\mu_i(x_i)}
    \tag*{\small{(cancellation; $\mu_{j^{i^*}} = \mu_{i^*, j^{i^*}} \circ \pi^{-1}_{j^{i^*}}$)}}
    \\
        &\qquad=
        \cancel{\mu_{i^{j^*}}(x_{i^{j^*}})}
        \frac{\mu_{i^{j^*},j^*}(x_{i^{j^*}}, x_{j^*})}{\cancel{\mu_{i^{j^*}}(x_{i^{j^*}})}} \,.
    \tag*{\small{(repeated application of the last few steps for remaining roots)}}
    \end{align*}
    \endgroup
    \item Finally, we prove that $\mu$ is a unique maximal-entropy measure, such that it has marginals $\big(\mu_{i,j}\big),(i,j) \in \mathcal{E}$.\\
    For an arbitrary measure $\lambda \in \mathcal{M}(\Omega)$ with $\mu_{ij} = \lambda \circ \pi^{-1}_{ij}$ for all $(i,j) \in \mathcal{E}$, we have
    \begingroup
    \allowdisplaybreaks
    \begin{align*}
        H(\lambda)
        &=
        H(\lambda, \mu) - D_{\mathrm{KL}}(\lambda, \mu)
    \tag*{\small{(cross-entropy formula)}}
        \\
        &=
        -\E_{\lambda}\big(\log \mu \big)
        - D_{\mathrm{KL}}(\lambda, \mu)
    \tag*{\small{(definition of cross-entropy)}}
        \\
        &=
        -\E_{\lambda}\bigg(\log \mu_{i^*}(\cdot_{i^*}) \prod_{\substack{(i,j) \in \mathcal{E}\\\mu_i(\cdot_i) > 0}}  \frac{\mu_{i,j}(\cdot_i, \cdot_j)}{\mu_i(\cdot_i)} \bigg)
        - D_{\mathrm{KL}}(\lambda, \mu)
    \tag*{\small{(definition of $\mu$)}}
        \\
        &=
        -\E_{\lambda}\big(\log \mu_{i^*}(\cdot_{i^*})\big)
        \\
        &\qquad- \bigg(\sum_{\substack{(i,j) \in \mathcal{E}\\\mu_i(\cdot_i) > 0}} \E_{\lambda}\big( \log\mu_{i,j}(\cdot_i, \cdot_j) \big) -   \E_{\lambda} \big(\log\mu_i(\cdot_i) \big) \bigg) 
        - D_{\mathrm{KL}}(\lambda, \mu)
    \tag*{\small{(properties of $\log$ and linearity of $\E$)}}
        \\
        &=
        -\E_{\lambda \circ \pi^{-1}_{i^*}}\big(\log \mu_{i^*}\big)
        \\
        &\qquad
        - \bigg(\sum_{\substack{(i,j) \in \mathcal{E}\\\mu_i(\cdot_i) > 0}} \E_{\lambda\circ\pi_{ij}^{-1}}\big( \log\mu_{i,j} \big) -   \E_{\lambda\circ\pi_i^{-1}} \big(\log\mu_i(\cdot_i) \big) \bigg) 
        - D_{\mathrm{KL}}(\lambda, \mu)
    \tag*{\small{(expect.~val.~of functions constant in an argument)}}
        \\
        &=
        -\E_{\mu_{i^*}}\big(\log \mu_{i^*}\big)
        \\
        &\qquad
        - \bigg(\sum_{\substack{(i,j) \in \mathcal{E}\\\mu_i(\cdot_i) > 0}} \E_{\mu_{ij}}\big( \log\mu_{i,j}(\cdot, \cdot) \big) -   \E_{\mu_i} \big(\log\mu_i(\cdot_i) \big) \bigg) 
        - D_{\mathrm{KL}}(\lambda, \mu)
    \tag*{\small{(assumption $\mu_{ij} = \lambda \circ \pi^{-1}_{ij}$ for all $(i,j) \in \mathcal{E}$)}}
        \\
        &=
        -\E_{\mu}\big(\log \mu\big) - D_{\mathrm{KL}}(\lambda, \mu)
    \tag*{\small{(by the reversed argument)}}
        \\
        &=
        H(\mu) - D_{\mathrm{KL}}(\lambda, \mu)
    \tag*{\small{(definition of entropy)}}
        \\
        &\leq
        H(\mu) \,.
    \tag*{\small{(property of the KL-divergence)}}
    \end{align*}
    \end{enumerate}
    \endgroup
\end{proof}

It turns out that the limits to the existence of a global measure that has given marginally-coupled bivariate measures as marginals are related to the consistency of the tree relaxation of \cref{thm:relaxation} in solving the original minimization of sums of bivariates.

We see in \cref{ex:counterrecon} that, in general, if we apply the tree relaxation to non-tree-structured sums of bivariates, we obtain a different problem. The problem we obtain differs both in its minimal value and in our ability to make sense of solutions by finding a global measure for a solution that is a set of marginally-coupled bivariate measures.

\begin{example}[Extension~Impossible]
\label{ex:counterrecon}
    Consider the functions $F: \{0,1\}^3 \to \R$, $f_{1,2}, f_{1,3}, f_{2,3}: \{0,1\}^2 \to \R$ with
    \[
        F(x_1, x_2, x_2) = f_{1,2}(x_1, x_2) + f_{1,3}(x_1, x_3) + f_{2,3}(x_2, x_3) \qquad \forall x \in \{0,1\}^3 \,,
    \]
    where
    \begin{align*}
        f_{1,2}(x_1, x_2) &\coloneqq \begin{cases}
            1 & \text{if } x_1=x_2=0
            \\
            0 & \text{if } x_1=x_2=1
            \\
            \infty & \text{else,}
        \end{cases}
        \\
        f_{1,3}(x_1, x_3) &\coloneqq \begin{cases}
            1 & \text{if } x_1=1-x_3=0
            \\
            0 & \text{if } x_1=x_3
            \\
            \infty & \text{else, and}
        \end{cases}
        \\
        f_{2,3}(x_2, x_3) &\coloneqq \begin{cases}
            1 & \text{if } x_2=1-x_3=0
            \\
            0 & \text{if } x_3=1-x_2=0
            \\
            \infty & \text{else.}
        \end{cases}
    \end{align*}
\begin{figure}%
\centering%
\begin{tikzpicture}[every node/.style={minimum size=1cm},on grid]%
\begin{scope}[every node/.append style={yslant=-0.5},yslant=-0.5]
  \shade[right color=blue!10, left color=blue!50] (0,0) rectangle +(2,2);
  \node at (0.5,1.5) {$0$};
  \node at (1.5,1.5) {$\infty$};
  \node at (0.5,0.5) {$\infty$};
  \node at (1.5,0.5) {$1$};
  \draw (0,0) grid (2,2);
\end{scope}
\begin{scope}[every node/.append style={yslant=0.5},yslant=0.5]
  \shade[right color=blue!50,left color=blue!10] (2,-2) rectangle +(2,2);
  \node at (3.5,-0.5) {$0$};
  \node at (2.5,-0.5) {$\infty$};
  \node at (3.5,-1.5) {$1$};
  \node at (2.5,-1.5) {$0$};
  \draw (2,-2) grid (4,0);
  \draw[->] (2.2,-2.2)   -- (4.2,-2.2) node[midway,below,shift={(0.,0.2)}] {$x_3$};
  \draw[->] (4.2,-2.0)   -- (4.2,0.0) node[midway,right,shift={(-.2,0.0)}] {$x_1$};
\end{scope}
\begin{scope}[every node/.append style={yslant=0.5,xslant=-1},yslant=0.5,xslant=-1]
  \shade[bottom color=blue!10, top color=blue!50] (4,2) rectangle +(-2,-2);
  \node at (3.5,1.5) {$\infty$};
  \node at (3.5,0.5) {$1$};
  \node at (2.5,1.5) {$0$};
  \node at (2.5,0.5) {$\infty$};
  \draw (2,0) grid (4,2);
  \draw[->] (4.2,0.2)   -- (4.2,2.2) node[midway,right,shift={(-0.2,0.0)}] {$x_2$};
\end{scope}%
\end{tikzpicture}%
\vspace*{-.5cm}
\caption{Geometric visualization of \cref{ex:counterrecon}. Measures $\mu_{1,2}, \mu_{1,3}$, and $\mu_{2,3}$ can be understood to be defined on the respective face of the cube. The relaxed objective function is the sum of the scalar product of each measure with the respective face.}%
\end{figure}
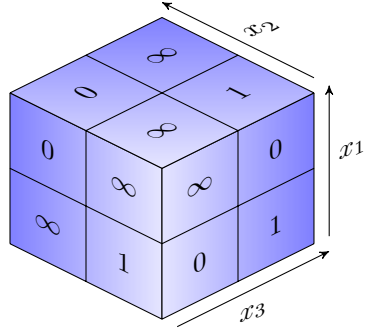%
    Then, we have
    \[
        F(x_1, x_2, x_3) = \begin{cases}
            3 & \text{if } x_1 = x_2 = 1-x_3 = 0
            \\
            \infty & \text{else}
        \end{cases}
    \]
    and the tree relaxation of \cref{thm:relaxation} (where in this example $\infty \cdot 0 \coloneqq 0$  and $\mathcal{G} \coloneqq(\mathcal{V}, \mathcal{E})\coloneqq(\{1,2,3\}, \{(1,2),(1,3),(2,3)\})$ is not a tree) yields
\begingroup
\allowdisplaybreaks
    \begin{align*}
        &\begin{cases}
        \displaystyle\argminmin_{\substack{\mu_{\cdot, \cdot} \in \mathcal{M}^+(\Omega_\cdot \times \Omega_\cdot)\\ \mu_\cdot \in \mathcal{M}_1^+(\Omega_\cdot)}} \sum_{(i,j) \in \mathcal{E}} \big\langle f_{i,j} , \mu_{i,j}\big\rangle
        \\
        \,\,\text{s.t.}\quad
        \mu_i = \mu_{i,j} \circ \pi_i^{-1} \,,\quad \forall (i,j) \in \mathcal{E}
        \\
        \phantom{\,\,\text{s.t.}\quad}
        \mu_j = \mu_{i,j} \circ \pi_j^{-1} \,,\quad \forall (i,j) \in \mathcal{E}
        \end{cases}
        \\
        &=
        \begin{cases}
        \displaystyle\argminmin_{\mu_{\cdot, \cdot} \in \mathcal{M}^+(\Omega_\cdot \times \Omega_\cdot)}  \big\langle f_{1,2} , \mu_{1,2}\big\rangle
        +
        \big\langle f_{1,3} , \mu_{1,3}\big\rangle
        +
        \big\langle f_{2,3} , \mu_{2,3}\big\rangle
        \\
        \,\,\text{s.t.}\quad
        \mu_{1,2}(\Omega_1 \times \Omega_2) = 1
        \\
        \phantom{\,\,\text{s.t.}\quad}
        \mu_{1,3}(\Omega_1 \times \Omega_3) = 1
        \\
        \phantom{\,\,\text{s.t.}\quad}
        \mu_{2,3}(\Omega_2 \times \Omega_3) = 1
        \\
        \phantom{\,\,\text{s.t.}\quad}
        \mu_{1,2} \circ \pi_1^{-1} = \mu_{1,3} \circ \pi_1^{-1}
        \\
        \phantom{\,\,\text{s.t.}\quad}
        \mu_{1,2} \circ \pi_2^{-1} = \mu_{2,3} \circ \pi_2^{-1}
        \\
        \phantom{\,\,\text{s.t.}\quad}
        \mu_{1,3} \circ \pi_3^{-1} = \mu_{2,3} \circ \pi_3^{-1}
        \end{cases}
    \tag*{\small{(definition of $\mathcal{E}$; eliminating marginal variables)}}
        \\
        &=
        \begin{cases}
        \displaystyle\argminmin_{\mu_{\cdot, \cdot} \in \mathcal{M}^+(\Omega_\cdot \times \Omega_\cdot)}  \mu_{1,2}(0,0)
        +
        \mu_{1,3}(0,1)
        +
        \mu_{2,3}(0,1)
        \\
        \,\,\text{s.t.}\quad
        \mu_{1,2}(0,0) + \mu_{1,2}(1,0)
        + \mu_{1,2}(0,1) + \mu_{1,2}(1,1) = 1
        \\
        \phantom{\,\,\text{s.t.}\quad}
        \mu_{1,3}(0,0) + \mu_{1,3}(1,0) + \mu_{1,3}(0,1) + \mu_{1,3}(1,1)= 1
        \\
        \phantom{\,\,\text{s.t.}\quad}
        \mu_{2,3}(0, 0) + \mu_{2,3}(1, 0) + \mu_{2,3}(0, 1) + \mu_{2,3}(1, 1) = 1
        \\
        \phantom{\,\,\text{s.t.}\quad}
        \mu_{1,2}(0,0)+ \mu_{1,2}(0,1) = \mu_{1,3}(0,0)+ \mu_{1,3}(0,1) 
        \\
        \phantom{\,\,\text{s.t.}\quad}
        \mu_{1,2}(1,0)+ \mu_{1,2}(1,1) = \mu_{1,3}(1,0)+ \mu_{1,3}(1,1) 
        \\
        \phantom{\,\,\text{s.t.}\quad}
        \mu_{1,2}(0,0)+ \mu_{1,2}(1,0) = \mu_{2,3}(0,0)+ \mu_{2,3}(0,1) 
        \\
        \phantom{\,\,\text{s.t.}\quad}
        \mu_{1,2}(0,1)+ \mu_{1,2}(1,1) = \mu_{2,3}(1,0)+ \mu_{2,3}(1,1) 
        \\
        \phantom{\,\,\text{s.t.}\quad}
        \mu_{1,3}(0,0)+ \mu_{1,3}(1,0) = \mu_{2,3}(0,0)+ \mu_{2,3}(1,0)
        \\
        \phantom{\,\,\text{s.t.}\quad}
        \mu_{1,3}(0,1)+ \mu_{1,3}(1,1) = \mu_{2,3}(0,1)+ \mu_{2,3}(1,1)
        \\
        \hline
        \phantom{\,\,\text{s.t.}\quad}
        \mu_{1,2}(0,1)=\mu_{1,2}(1,0)= \mu_{1,3}(1,0)= \mu_{2,3}(0,0) =\mu_{2,3}(1,1)=0
        \end{cases}
    \tag*{\small{(definitions of $f_{\cdot,\cdot}$ \& $\Omega_\cdot$)}}
        \\
        &=
        \begin{cases}
        \displaystyle\argminmin_{\mu_{\cdot, \cdot} \in \mathcal{M}^+(\Omega_\cdot \times \Omega_\cdot)}  \mu_{1,2}(0,0)
        +
        \mu_{1,3}(0,1)
        +
        \mu_{2,3}(0,1)
        \\
        \,\,\text{s.t.}\quad
        \mu_{1,2}(0,0) + \mu_{1,2}(1,1) = 1
        \\
        \phantom{\,\,\text{s.t.}\quad}
        \mu_{1,3}(0,0) + \mu_{1,3}(0,1) + \mu_{1,3}(1,1)= 1
        \\
        \phantom{\,\,\text{s.t.}\quad}
        \mu_{2,3}(1, 0) + \mu_{2,3}(0, 1) = 1
        \\
        \phantom{\,\,\text{s.t.}\quad}
        \mu_{1,2}(0,0) = \mu_{1,3}(0,0)+ \mu_{1,3}(0,1) 
        \\
        \phantom{\,\,\text{s.t.}\quad}
        \mu_{1,2}(1,1) = \mu_{1,3}(1,1) 
        \\
        \phantom{\,\,\text{s.t.}\quad}
        \mu_{1,2}(0,0) =\mu_{2,3}(0,1) 
        \\
        \phantom{\,\,\text{s.t.}\quad}
        \mu_{1,2}(1,1) = \mu_{2,3}(1,0)
        \\
        \phantom{\,\,\text{s.t.}\quad}
        \mu_{1,3}(0,0) = \mu_{2,3}(1,0)
        \\
        \phantom{\,\,\text{s.t.}\quad}
        \mu_{1,3}(0,1)+ \mu_{1,3}(1,1) = \mu_{2,3}(0,1)
        \\
        \hline
        \phantom{\,\,\text{s.t.}\quad}
        \mu_{1,2}(0,1)=\mu_{1,2}(1,0)= \mu_{1,3}(1,0)= \mu_{2,3}(0,0) =\mu_{2,3}(1,1)=0
    \tag*{\small{(plugging the last equation into all others)}}
        \end{cases}
        \\
        &=
        \begin{cases}
        \displaystyle\argminmin_{\mu_{\cdot, \cdot} \in \mathcal{M}^+(\Omega_\cdot \times \Omega_\cdot)}  2\mu_{1,2}(0,0)
        +
        \mu_{1,3}(0,1)
        \\
        \,\,\text{s.t.}\quad
        2\mu_{1,2}(1,1) + \mu_{1,3}(0,1) = 1
        \\
        \phantom{\,\,\text{s.t.}\quad}
        \mu_{1,2}(0,0) = \mu_{1,2}(1,1)+ \mu_{1,3}(0,1) 
        \\
        \phantom{\,\,\text{s.t.}\quad}
        \mu_{1,3}(0,1)+ \mu_{1,2}(1,1) = \mu_{1,2}(0,0)
        \\
        \hline
        \phantom{\,\,\text{s.t.}\quad}
        \mu_{1,2}(0,1)=\mu_{1,2}(1,0)= \mu_{1,3}(1,0)= \mu_{2,3}(0,0) =\mu_{2,3}(1,1)=0
        \\
        \phantom{\,\,\text{s.t.}\quad}
        \mu_{1,3}(1,1) = \mu_{1,2}(1,1)
        \\
        \phantom{\,\,\text{s.t.}\quad}
        \mu_{2,3}(0,1)  = \mu_{1,2}(0,0)
        \\
        \phantom{\,\,\text{s.t.}\quad}
        \mu_{1,3}(0,0) = \mu_{2,3}(1,0) = \mu_{1,2}(1,1)
        \\
        \phantom{\,\,\text{s.t.}\quad}
        \mu_{1,2}(1,1) = 1- \mu_{1,2}(0,0)
        \\
        \phantom{\,\,\text{s.t.}\quad}
        \mu_{2,3}(1, 0) = 1-\mu_{1,2}(0, 0)
        \end{cases}
    \tag*{\small{(replacing variables constrained by trivial equations)}}
        \\
        &=
        \begin{cases}
        \displaystyle\argminmin_{\mu_{\cdot, \cdot} \in \mathcal{M}^+(\Omega_\cdot \times \Omega_\cdot)}  1
        \\
        \,\,\text{s.t.}\quad
        \cancel{\mu_{1,2}(0,0) =  1-\mu_{1,2}(1,1)}
        \\
        \phantom{\,\,\text{s.t.}\quad}
        1-\mu_{1,2}(1,1) = \mu_{1,2}(0,0)
        \\
        \hline
        \phantom{\,\,\text{s.t.}\quad}
        \mu_{1,2}(0,1)=\mu_{1,2}(1,0)= \mu_{1,3}(1,0)= \mu_{2,3}(0,0) =\mu_{2,3}(1,1)=0
        \\
        \phantom{\,\,\text{s.t.}\quad}
        \mu_{1,3}(1,1) = \mu_{1,2}(1,1)
        \\
        \phantom{\,\,\text{s.t.}\quad}
        \mu_{2,3}(0,1)  = \mu_{1,2}(0,0)
        \\
        \phantom{\,\,\text{s.t.}\quad}
        \mu_{1,3}(0,0) = \mu_{2,3}(1,0) = \mu_{1,2}(1,1)
        \\
        \phantom{\,\,\text{s.t.}\quad}
        \mu_{1,2}(1,1) = 1- \mu_{1,2}(0,0)
        \\
        \phantom{\,\,\text{s.t.}\quad}
        \mu_{2,3}(1, 0) = 1-\mu_{1,2}(0, 0)
        \\
        \phantom{\,\,\text{s.t.}\quad}
        \mu_{1,3}(0,1) = 1 - 2\mu_{1,2}(1,1) \,.
        \end{cases}
    \tag*{\small{(replacing variables constrained by trivial equations)}}
    \end{align*}
    \endgroup
    Clearly, we have a minimal value of $1$ at a (non-unique) minimum of
    \[
        \mu_{1,2}(a,b) = \mu_{1,3}(a,b) = 1/2 - \mu_{2,3}(a,b) = \begin{cases}
            1/2 & \text{if } a=b
            \\
            0 & \text{else}
        \end{cases}
        \qquad \forall a,b \in \{0,1\} \,.
    \]
    Further, there exists no measure $\mu$ on $\{0,1\}^3$ such that it has marginals $\mu_{1,2}, \mu_{1,3}, \mu_{2,3}$, i.e.~there exists no $\mu$ such that
    \begin{gather*}
        \mu_{1,2} = \mu \circ \pi_{1,2}^{-1}\,, \qquad \mu_{1,3} = \mu \circ \pi_{1,3}^{-1} \,,\,\, \text{and} \qquad \mu_{2,3} = \mu \circ \pi_{2,3}^{-1} \,.
    \end{gather*}
    This is, as $\mu_{1,2}(1,0) = \mu_{1,2}(0,1) = \mu_{1,3}(1,0) = \mu_{1,3}(0,1) = 0$ would imply by marginalization that the support of $\mu$ is contained in $\{(0,0,0), (1,1,1)\}$. However, since $\mu_{2,3}(0,0) = \mu_{2,3}(1,1) = 0$, by marginalization, $\mu$ must be zero on $\{(\cdot, 0,0), (\cdot, 1,1)\}$. Therefore, $\mu$ cannot have the specified marginals.\\
    This example shows that neither the non-global relaxation of \cref{thm:relaxation} nor the reconstruction of \cref{thm:rfm} can be extended beyond trees in general. 
\end{example}

The following two postulates are inspired by and very similar to the results in \cite{ochs2024}.\\
The goal is to derive the Lagrangian dual of the star relaxation in \cref{thm:dual}, which restricts to star graph-indexed sums of bivariates. The dual of the tree relaxation would follow similarly and will be described briefly in \cref{thm:duallinprog} of \cref{sec:algorithms}. Later in our derivation, even the specific results for star graph-indexed sums of bivariates will prove useful, as we will be able to use them for a type of coordinate ascent minimization algorithm, similar and motivated to the results in \cite{ochs2024}.

Further, \cref{thm:dual} claims that the dual of the star relaxation is strong and involves fewer variables than its primal---the star relaxation.  This dual maximization can be understood as the objective to find lower bounds for the functions $\big(\min_{x_j \in \Omega_j} f_{i^*,j}(\cdot, x_j)\big) \in \R^{\Omega_{i^*}}$ for all $(i^*, j) \in \mathcal{N}(i^*)$, called \emph{$\min$-marginals}, such that the sum of the $\min$-marginals has the largest possible minimum.
\begin{postulate}[Star~Lagrangian~Duality]
    If $\mathcal{G}= (\mathcal{V}, \mathcal{E})$ is a star graph with root $i^* \in \mathcal{V}$ and $\mathcal{E}=\mathcal{N}(i^*)$, the Lagrangian dual of the star relaxation of \cref{thm:relaxation} is strong and
    \begin{align*}
    &\begin{cases}
    \displaystyle\min_{\substack{\mu_{i^*, \cdot} \in \mathcal{M}^+(\Omega_{i^*} \times \Omega_\cdot)\\ \mu_{i^*} \in \mathcal{M}_1^+(\Omega_{i^*})}} \sum_{(i^*,j) \in \mathcal{N}(i^*)} \big\langle f_{i^*,j} , \mu_{i^*,j}\big\rangle
        \\
        \,\,\text{s.t.}\quad
        \mu_{i^*} = \mu_{i^*,j} \circ \pi_{i^*}^{-1} \,,\quad \forall (i^*,j) \in \mathcal{N}(i^*)
        \end{cases}
    \tag*{\small{(star relaxation)}}
    \\
    &\qquad=
    \begin{cases}\displaystyle\max_{\rho_{i^*, \cdot} \in \R^{\Omega_{i^*}}}\bigg(\min_{x_{i^*} \in \Omega_{i^*}}
     \big(\textstyle\sum_{(i^*,j) \in \mathcal{N}(i^*)}  \rho_{i^*,j}\big)(x_{i^*}) \bigg)
    \\
     \,\,\text{s.t.}\quad \rho_{i^*, j} \leq \min_{x_j \in \Omega_j} f_{i^*, j}(\cdot, x_j)  \,,\quad \forall (i^*,j) \in \mathcal{N}(i^*) \,.
    \end{cases}
    \tag*{\small{(dual star relaxation)}}
    \end{align*}
\label{thm:dual}
\end{postulate}

We conclude the analysis of the star relaxation by deriving a closed-form solution to the dual star relaxation in \cref{thm:nscsdv} using \cref{lem:marginalize}. The solution is not unique, but is characterized by linear inequalities. Intuitively, optimal solutions exceed the minimum of the sum of the $\min$-marginals, while being individually bounded by them.

\begin{lemma}
\label{lem:marginalize}%
    In the setting of \cref{def:biv} and \cref{thm:dual}, we have for all $(i^*, j) \in \mathcal{N}(i^*)$ that
    \begin{align*}
        \min_{ \mu_{i^*, j} \in \mathcal{M}^+(\Omega_{i^*} \times \Omega_j)}\big\langle f_{i^*,j} - \rho_{i^*,j}, \mu_{i^*,j}\big\rangle
        &=
        \min_{ \nu_{i^*, j} \in \mathcal{M}^+(\Omega_{i^*})}
      \big\langle m_{i^*,j} - \rho_{i^*,j}, \nu_{i^*,j}\big\rangle
      \\
      &=
      \delta_{\rho_{i^*, j} \leq m_{i^*, j}} \coloneqq
      \begin{cases}
          0 & \text{if } \rho_{i^*, j} \leq m_{i^*, j}
          \\
          -\infty & \text{else}\,,
      \end{cases}
    \end{align*}
    where $m_{i^*, j}(\cdot) \coloneqq \min_{x_j \in \Omega_j} f_{i^*, j}(\cdot, x_j)$ and $\rho_{i^*, j} \in \R^{\Omega_{i^*}}$.
\end{lemma}
\begin{proof}
    Clear.
\end{proof}

\begin{postulate}[Necessary and Sufficient Conditions for Star Dual Variables]
\label{thm:nscsdv}
 If $\mathcal{G}= (\mathcal{V}, \mathcal{E})$ is a star graph with root $i^* \in \mathcal{V}$ and $\mathcal{E}=\mathcal{N}(i^*)$, the variables of the Lagrangian dual of the star relaxation of \cref{thm:dual} are optimal, i.e.
    \begin{align*}
    &\rho^*_{i^*, \cdot} \in \begin{cases}\displaystyle\argmax_{\rho_{i^*, \cdot} \in \R^{\Omega_{i^*}}}\bigg(\min_{        x_{i^*} \in \Omega_{i^*}}
     \big(\textstyle\sum_{(i^*,j) \in \mathcal{N}(i^*)}  \rho_{i^*,j}\big)(x_{i^*}) \bigg)
    \\
     \,\,\text{s.t.}\quad \rho_{i^*, j} \leq  m_{i^*, j}  \,,\quad \forall (i^*,j) \in \mathcal{N}(i^*)
    \end{cases}
    \tag*{\small{(dual star relaxation)}}
      \\
      & \iff
      \begin{cases}\displaystyle \min_{x_{i^*} \in \Omega_{i^*}} \big(\textstyle\sum_{(i^*, j) \in \mathcal{N}(i^*)} m_{i^*, j}\big)(x_{i^*}) \leq \sum_{(i^*,j) \in \mathcal{N}(i^*)} \rho_{i^*, j}
      \\      
      \,\,\text{s.t.}\quad\rho_{i^*, j} \leq m_{i^*, j}  \,,\quad \forall (i^*,j) \in \mathcal{N}(i^*) \,,
      \end{cases}
    \tag*{\small{(NSCSDV)}}
    \end{align*}
    where $m_{i^*, j}(\cdot) \coloneqq \min_{x_j \in \Omega_j} f_{i^*, j}(\cdot, x_j)$ for all $(i^*, j) \in \mathcal{N}(i^*)$.
\end{postulate}

\subsection{Entropy-Regularized Relaxation Principle}
\label{subsec:entropyrelaxation}%

As solutions to the dual star relaxation from \cref{thm:nscsdv} are not unique, and with applications of the results to non-tree-structured sums of bivariates in mind, we aim to derive an \emph{entropy-regularized star relaxation}. The entropy-regularization will, in fact, turn out to help us in deriving a unique solution, and the new objective that we will construct will uniformly approximate our known star relaxation. The heuristic hope is that an entropy-regularization would also help in \enquote{distributing} measure by preferring higher entropies in coupled bivariate marginal measures, such that the possibility of suboptimal convergence may be partly prevented.

To this end, we formulate a \enquote{parameter-free} version of the maximal-entropy measure that has given coupled bivariate marginal measures, known from \cref{thm:rfm}, in \cref{cor:rmer}. We also find a simple expression for its entropy that is a linear combination of the entropies of its marginals.

\begin{corollary}[Representation of the Maximal-Entropy Reconstruction]
\label{cor:rmer}
    If the graph $\mathcal{G}= (\mathcal{V}, \mathcal{E})$ is a tree and given
    \begin{align*}
    \begin{cases}
        \mu_{i,j} \in \mathcal{M}^+(\Omega_i \times \Omega_j)\,, &\forall (i,j) \in \mathcal{E}
        \\
        \mu_i \in \mathcal{M}^+_1(\Omega_i) \,,&\forall i \in \mathcal{V}
        \\
        \,\,\text{s.t.}\quad
        \mu_i = \mu_{i,j} \circ \pi_i^{-1} \,,&\forall (i,j) \in \mathcal{E}
        \\
        \phantom{\,\,\text{s.t.}\quad}
        \mu_j = \mu_{i,j} \circ \pi_j^{-1} \,,&\forall (i,j) \in \mathcal{E} \,,
    \end{cases}
    \end{align*}
    then the maximal-entropy measure $\mu \in \mathcal{M}_1^+(\Omega)$ with marginals $\big(\mu_{i,j}\big),(i,j) \in \mathcal{E}$, i.e.
    \[
        \mu_{i,j} = \mu \circ \pi_{i,j}^{-1} \qquad \forall (i,j) \in \mathcal{E} \,,
    \]
    known from \cref{thm:rfm}, has a representation
    \[
        \mu(x) = \begin{cases}
            \bigg(\prod_{(i,j) \in \mathcal{E}} \mu_{i,j}(x_i, x_j) \bigg)\!\!\bigg/\!\! \bigg( \prod_{i \in \mathcal{V}} \mu_i(x_i)^{\abs{\widetilde{\mathcal{N}}(i)}-1} \bigg)
            & \text{if } \mu_i(x_i) >0  \,\,\, \forall i \in \mathcal{V}
            \\
            0 & \text{else}
        \end{cases}
        \,,\,\, \forall x \in \Omega \,.
    \]
    Further, we have for its entropy
    \[
        H(\mu) = \sum_{(i,j) \in \mathcal{E}} H(\mu_{i,j}) - \sum_{i \in \mathcal{V}} (\abs{\widetilde{\mathcal{N}}(i)}-1)\cdot H(\mu_i) \,.
    \]
\end{corollary}
\begin{proof}\leavevmode
W.l.o.g.,~we assume $\mathcal{G}$ to be oriented. Otherwise, we orient it.\\
Note, that we now have
\[
    \abs{\widetilde{\mathcal{N}}(i)}-1 = \abs{\mathcal{N}(i)}  \,,\,\, \forall i \in \mathcal{V}\setminus\{i^*\}\,,
\]
where the root of the oriented tree is denoted by $i^* \in \mathcal{V}$.
\begin{enumerate}[label=\roman*., listparindent=1.5em]
    \item We have to prove that $\mu$ as defined in \cref{thm:rfm} equals the representation given. Clearly, this is the case if there exists $x \in \Omega$ and $i \in \mathcal{V}$ with $\mu_i(x_i) = 0$ by the \enquote{else}-condition of the definitions. Therefore, we assume that for all $x \in \Omega$, we have $\mu_i(x_i) > 0$ for all $i \in \mathcal{V}$.\\
    \indent We prove the equality by induction over trees. The fact is trivially true for a tree with the root as the only parent vertex. We assume, therefore, that the equality holds for trees of a fixed depth. Denote by $(i^{j^*}, j^*) \in \mathcal{E}$ an arbitrary leaf edge of $\mathcal{G}$. Then,
\begingroup%
\allowdisplaybreaks%
    \begin{align*}
        \mu(x)
        &=
        \mu_{i^*}(x_{i^*}) \prod_{(i,j) \in \mathcal{E}}  \frac{\mu_{i,j}(x_i, x_j)}{\mu_i(x_i)}
    \tag*{\small{(definition in \cref{thm:rfm})}}
    \\
    &=
    \bigg(\mu_{i^*}(x_{i^*}) \prod_{(i,j) \in \mathcal{E}\setminus \{(i^{j^*}, j^*)\}}  \frac{\mu_{i,j}(x_i, x_j)}{\mu_i(x_i)} \bigg) \frac{\mu_{i^{j^*},j^*}(x_{i^{j^*}},x_{j^*})}{\mu_{i^{j^*}}(x_{i^{j^*}})}
    \tag*{\small{(splitting the product)}}
    \\
    &=
    \bigg(\bigg(\prod_{(i,j) \in \mathcal{E}\setminus\{(i^{j^*}, j^*)\}} \mu_{i,j}(x_i,x_j) \bigg)\bigg/ \bigg( \prod_{i \in (\mathcal{V}\setminus\{j^*\})} \mu_i(x_i)^{\abs{\mathcal{N}^{j^*}(i)}} \bigg) \bigg) \frac{\mu_{i^{j^*},j^*}(x_{i^{j^*}},x_{j^*})}{\mu_{i^{j^*}}(x_{i^{j^*}})}
    \tag*{\small{(induction hypothesis; $N^{j^*}(i) \coloneqq \{(i,j) \in \mathcal{E} \setminus \{(i^{j^*}, j^*)\}$)}}
    \\
    &=
    \bigg(\prod_{(i,j) \in \mathcal{E}} \mu_{i,j}(x_i, x_j) \bigg)\bigg/ \bigg( \prod_{i \in \mathcal{V}} \mu_i(x_i)^{\abs{\mathcal{N}^{j^*}(i)}}  \bigg)
    \tag*{\small{(by$\quad\abs{\mathcal{N}^{j^*}(i)} = \abs{\mathcal{N}(i)}$, if $i\neq i^{j^*}\,;\quad\abs{\mathcal{N}^{j^*}(i^{j^*})} + 1 = \abs{\mathcal{N}(i^{j^*})}\,;\quad\abs{\mathcal{N}(j^*)} = 0$)}} \,.
    \end{align*}
    \endgroup
    \item W.l.o.g.~assume $\mu(x) > 0$ for all $x \in \Omega$. Otherweise, restrict $\Omega$. We have
    \begingroup%
    \allowdisplaybreaks%
    \begin{align*}
    H(\mu)
    &=
    -\E_\mu\big( \log \mu(\cdot ) \big)
    \tag*{\small{(definition of the entropy $H$)}}
    \\
    &=
    -\E_\mu\bigg( \log \bigg(\prod_{(i,j) \in \mathcal{E}} \mu_{i,j}(\cdot_i,\cdot_j) \bigg)\bigg/ \bigg( \prod_{i \in \mathcal{V}} \mu_i\big(\cdot_i)^{\abs{\widetilde{\mathcal{N}}(i)}-1} \bigg) \bigg)
    \tag*{\small{(representation of $\mu$)}}
    \\
    &=
    -\sum_{(i,j) \in \mathcal{E}} \E_\mu\bigg( \log \mu_{i,j}(\cdot_i, \cdot_j) \bigg) +\sum_{i \in \mathcal{V}}(\abs{\widetilde{\mathcal{N}}(i)}-1) \cdot \E_\mu\bigg( \log \mu_i\big(\cdot_i\big) \bigg)
    \tag*{\small{(properties of the logarithm; linearity of $\E$)}}
    \\
    &=
    -\sum_{(i,j) \in \mathcal{E}} \E_{\mu\circ \pi_{ij}^{-1}}\bigg( \log \mu_{i,j} \bigg) +   \sum_{i \in \mathcal{V}}(\abs{\widetilde{\mathcal{N}}(i)}-1) \cdot \E_{\mu\circ \pi_i^{-1}}\bigg( \log \mu_i \bigg)
    \tag*{\small{(expect.~val.~of functions constant in an argument)}}
    \\
    &=
    -\sum_{(i,j) \in \mathcal{E}} \E_{\mu_{ij}}\bigg( \log \mu_{i,j} \bigg) +   \sum_{i \in \mathcal{V}}(\abs{\widetilde{\mathcal{N}}(i)}-1) \cdot \E_{\mu_i}\bigg( \log \mu_i \bigg)
    \tag*{\small{(assumption $\mu_{ij} = \lambda \circ \pi^{-1}_{ij}$ for all $(i,j) \in \mathcal{E}$)}}
    \\
    &=
    \sum_{(i,j) \in \mathcal{E}} H(\mu_{ij}) -   \sum_{i \in \mathcal{V}}(\abs{\widetilde{\mathcal{N}}(i)}-1) \cdot H(\mu_i) \,.
    \tag*{\small{(definition of the entropy $H$)}}
    \end{align*}
\endgroup
\end{enumerate}
\end{proof}

Next, we prove the concavity of the function that, given a family of bivariate measures, assigns the entropy of a global measure that has the family as marginals. This result will be essential to prove the strong Lagrangian duality in \cref{thm:entropydual}.

\begin{lemma}
\label{lem:concavity}
    If $\mathcal{G}= (\mathcal{V}, \mathcal{E})$ is a star graph with root $i^* \in \mathcal{V}$ and $\mathcal{E}=\mathcal{N}(i^*)$, the function
    \[
        \mu_{i^*,j} \in \mathcal{M}_1^+(\Omega_{i^*}\times \Omega_j), (i^*,j) \in \mathcal{E} \longmapsto \bigg(\sum_{(i^*,j) \in \mathcal{E}} H(\mu_{i^*,j})\bigg) -  \big(\abs{\widetilde{\mathcal{N}}(i^*)}-1\big) \cdot H\big(\mu_{i^*,j}\circ\pi_{i^*}^{-1}\big)
    \]
    is concave.
\end{lemma}
\begin{proof}
    Let $\mu_{i^*,j},\lambda_{i^*,j} \in \mathcal{M}_1^+(\Omega_{i^*}\times \Omega_j)$, $t \in [0,1]$, and define $z_{i^*,j} \coloneqq t \mu_{i^*,j} + (1-t)\lambda_{i^*,j}$.
    Due to concavity of the entropy and $\abs{\mathcal{E}} = \abs{\widetilde{\mathcal{N}}(i^*)}$, it is sufficient to show the concavity of
\begingroup%
\allowdisplaybreaks%
    \begin{align*}
         &\sum_{(i^*,j) \in \mathcal{E}} H(z_{i^*,j}) -  H\big(z_{i^*,j}\circ\pi_{i^*}^{-1}\big)
        \\
        &\quad=
        \sum_{(i^*,j) \in \mathcal{E}}\log(\abs{\Omega_j}) - D_{\mathrm{KL}}\big(z_{i^*, j}, z_{i^*,j}\circ\pi_{i^*}^{-1}(\cdot_{i^*})/\abs{\Omega_j}\big)
    \tag*{\small{(see below argument)}}
        \\
        &\quad=
        \sum_{(i^*,j) \in \mathcal{E}}
        \log(\abs{\Omega_j})
        \\
        &\quad\qquad\quad
        -D_{\mathrm{KL}}\bigg(t \mu_{i^*,j} + (1-t)\lambda_{i^*,j}\,,\,\, t\mu_{i^*,j}\circ\pi_{i^*}^{-1}(\cdot_{i^*})/\abs{\Omega_j} + (1-t)\lambda_{i^*,j}\circ\pi_{i^*}^{-1}(\cdot_{i^*})/\abs{\Omega_j}\bigg)
    \tag*{\small{(definition of $z_{i^*,j}$ \& linearity of marginalization)}}
        \\
        &\quad\geq
        \sum_{(i^*,j) \in \mathcal{E}} t \bigg(\log(\abs{\Omega_j})- D_{\mathrm{KL}}\bigg( \mu_{i^*,j} \,,\,\, \mu_{i^*,j}\circ\pi_{i^*}^{-1}(\cdot_{i^*})/\abs{\Omega_j} \bigg)\bigg)
        \\
        &\quad\qquad\quad
        +
        (1-t) \bigg(\log(\abs{\Omega_j})- D_{\mathrm{KL}}\bigg( \lambda_{i^*,j} \,,\,\, \lambda_{i^*,j}\circ\pi_{i^*}^{-1}(\cdot_{i^*})/\abs{\Omega_j} \bigg)\bigg)
    \tag*{\small{(convexity of $\mathrm{KL}$-divergence)}}
        \\
        &\quad=
        t \bigg(\sum_{(i^*,j) \in \mathcal{E}} H(\mu_{i^*,j}) -  H\big(\mu_{i^*,j}\circ\pi_{i^*}^{-1}\big)\bigg)
        +
        (1-t) \bigg(\sum_{(i^*,j) \in \mathcal{E}} H(\lambda_{i^*,j}) -  H\big(\lambda_{i^*,j}\circ\pi_{i^*}^{-1}\big)\bigg) \,.
    \tag*{\small{(see below argument; linearity)}}
    \end{align*}
\endgroup
    It follows a similar argument as in \cite[Appendix~A]{globerson07a}. We have a representation of what is know as the conditional entropy of $\mu_{i^*, j}$, as
\begingroup%
\allowdisplaybreaks%
    \begin{align*}
        &H(\mu_{i^*,j}) -  H\big(\mu_{i^*,j}\circ\pi_{i^*}^{-1}\big)
        \\
        &\qquad=
        -\E_{\mu_{i^*,j}}\big(\log \mu_{i^*,j}\big) + \E_{\mu_{i^*,j}\circ\pi_{i^*}^{-1}}\big(\log \mu_{i^*,j}\circ\pi_{i^*}^{-1}\big)
    \tag*{\small{(definition entropy)}}
        \\
        &\qquad=
        -\E_{\mu_{i^*,j}}\bigg(\log \mu_{i^*,j} - \log \mu_{i^*,j}\circ\pi_{i^*}^{-1}(\cdot_{i^*})\bigg)
    \tag*{\small{(constant integration)}}
        \\
        &\qquad=
        -\E_{\mu_{i^*,j}}\bigg(\log \frac{ \mu_{i^*,j}}{\mu_{i^*,j}\circ\pi_{i^*}^{-1}(\cdot_{i^*})}\bigg)
    \tag*{\small{(properties of $\log$)}}
        \\
        &\qquad=
        \log(\abs{\Omega_j}) -\log(\abs{\Omega_j}) -\E_{\mu_{i^*,j}}\bigg(\log \frac{ \mu_{i^*,j}}{\mu_{i^*,j}\circ\pi_{i^*}^{-1}(\cdot_{i^*})}\bigg)
    \tag*{\small{(constant addition)}}
        \\
        &\qquad=
        \log(\abs{\Omega_j})-\E_{\mu_{i^*,j}}\bigg(\log \frac{ \mu_{i^*,j}(\cdot, \cdot)}{\mu_{i^*,j}\circ\pi_{i^*}^{-1}(\cdot_{i^*})/\abs{\Omega_j}} \bigg)
    \tag*{\small{(linearity \& properties of $\log$)}}
        \\
        &\qquad=
        \log(\abs{\Omega_j}) - D_{\mathrm{KL}}\big(\mu_{i^*, j}, \mu_{i^*,j}\circ\pi_{i^*}^{-1}(\cdot_{i^*})/\abs{\Omega_j}\big) \,.
    \tag*{\small{(definition of $\mathrm{KL}$-divergence; marginalized distribution is constant in marginalized coordinate)}}
    \end{align*}
\endgroup
\end{proof}

By subtracting a scaled version of the entropy function of \cref{lem:concavity} from the star relaxation of \cref{thm:relaxation}, one obtains a new regularized objective, termed \emph{entropy-regularized star relexation}, which, due to the boundedness of the entropy, uniformly approximates the star relaxation. \cref{thm:entropydual} shows that the dual of this new objective now has a unique solution and derives its representation based on the \emph{$\varepsilon$-$\mathrm{LogSumExp}$}.

Similarly to the results in \cref{subsec:relaxation}, the theorem will also aid in deriving optimization algorithms that are applicable to generic sums of bivariates.

\begin{theorem}[Entropy~Star~Lagrangian~Duality]
    If $\mathcal{G}= (\mathcal{V}, \mathcal{E})$ is a star graph with root $i^* \in \mathcal{V}$ and $\mathcal{E}=\mathcal{N}(i^*)$, the Lagrangian dual of the following entropy-regularized star relaxation of \cref{thm:relaxation} is strong, has a unique solution, and
\begingroup%
\allowdisplaybreaks%
    \begin{align*}
    &\begin{cases}
    \displaystyle\min_{\substack{\mu_{i^*, \cdot} \in \mathcal{M}^+(\Omega_{i^*} \times \Omega_\cdot)\\ \mu_{i^*} \in \mathcal{M}_1^+(\Omega_{i^*})}} \bigg(\sum_{(i^*,j) \in \mathcal{N}(i^*)} \big\langle f_{i^*,j} , \mu_{i^*,j}\big\rangle - \varepsilon H(\mu_{i^*,j})\bigg)+ \varepsilon(\abs{\mathcal{N}(i^*)}-1) \cdot H(\mu_{i^*}) 
    \\
    \,\,\text{s.t.}\quad
    \mu_{i^*} = \mu_{i^*,j} \circ \pi_{i^*}^{-1} \,,\quad \forall (i^*,j) \in \mathcal{N}(i^*)\qquad\qquad\qquad
    \text{\small{\emph{(entropy-regularized star relaxation)}}}%
    \end{cases}%
    \\
    &\quad=
     \displaystyle\max_{\rho_{i^*, \cdot} \in \R^{\Omega_{i^*}}}\bigg(\bigg(\min_{
        x_{i^*}\in  \Omega_{i^*}}\big(\textstyle\sum_{(i^*,j) \in \mathcal{N}(i^*) } \rho_{i^*, j}\big)(x_{i^*}) \bigg)
    \\
    &\qquad\qquad\qquad\qquad\qquad
    -\varepsilon\!\!\!\displaystyle\sum_{(i^*,j) \in \mathcal{N}(i^*) }\sum_{x_{i^*} \in \Omega_{i^*}} \sum_{x_j \in \Omega_j}
        \exp \big(\big(\rho_{i^*,j}(x_{i^*})-f_{i^*,j}(x_{i^*},x_j)\big)/\varepsilon -1 \big)\bigg)
    \tag*{\small{(entropy dual star relaxation)}}
    \\
    &\quad=
    \begin{cases}
     \rho - \abs{\mathcal{N}(i^*)}
    \\
    \text{where }
    \\
    \rho_{i^*,j}(x_{i^*})
     = \mathrm{lse}_\varepsilon\big(f_{i^*,j}(x_{i^*}, \cdot)\big) + \tfrac{1}{\abs{\mathcal{N}(i^*)}} \bigg(\rho
     -\!\!\!\!\!\!\!\!\!\!\displaystyle\sum_{(i^*,k)\in \mathcal{N}(i^*)} \!\!\!\!\!\!\!\!\!\mathrm{lse}_\varepsilon\big(f_{i^*,k}(x_{i^*}, \cdot)\big)\bigg)\,
     \\
     \qquad\qquad\forall x_{i^*} \in \Omega_{i^*} \,,\,\, (i^*, j) \in \mathcal{N}(i^*)
    \\
    \rho = \varepsilon\abs{\mathcal{N}(i^*)} +\varepsilon\abs{\mathcal{N}(i^*)}\log\big(\abs{\mathcal{N}(i^*)}/\varepsilon\big)
    \\
    \qquad
    - \abs{\mathcal{N}(i^*)}\mathrm{lse}_\varepsilon^{x_{i^*}, j}\bigg(\mathrm{lse}^{x_j}_\varepsilon\big(f_{i^*,j}(x_{i^*}, x_j)\big) +\mathrm{lse}^{x_j}_\varepsilon\big(-f_{i^*, j}(x_{i^*}, x_j)\big)\\
    \qquad\qquad\qquad\qquad\qquad
    -\tfrac{1}{\abs{\mathcal{N}(i^*)}}\!\!\!\!\!\!\!\!\!\!\displaystyle\sum_{(i^*,k)\in \mathcal{N}(i^*)} \!\!\!\!\!\!\!\!\!\mathrm{lse}^{x_k}_\varepsilon\big(f_{i^*,k}(x_{i^*}, x_k)\big)\bigg) \,.
    \end{cases}
    \end{align*}%
\endgroup%
\label{thm:entropydual}%
\end{theorem}%
\newpage%
\begin{proof}%
The equivalence of the value of Lagrangian primal and dual, i.e.~the strong duality, as well as the uniqueness of the solution can be derived directly, as
\begingroup%
\allowdisplaybreaks%
\begin{align*}%
    &
    \begin{cases}
    \displaystyle\min_{\substack{\mu_{i^*, \cdot} \in \mathcal{M}^+(\Omega_{i^*} \times \Omega_\cdot)\\ \mu_{i^*} \in \mathcal{M}_1^+(\Omega_{i^*})}} \bigg(\sum_{(i^*,j) \in \mathcal{N}(i^*)}\big\langle f_{i^*,j} , \mu_{i^*,j}\big\rangle -  \varepsilon H(\mu_{i^*,j}) \bigg)+ \varepsilon(\abs{\mathcal{N}(i^*)}-1) \cdot H(\mu_{i^*})
        \\
        \,\,\text{s.t.}\quad
        \mu_{i^*} = \mu_{i^*,j} \circ \pi_{i^*}^{-1} \,,\quad \forall (i^*,j) \in \mathcal{N}(i^*) \,.
        \end{cases}
    \\ 
    &\quad=
    \displaystyle\max_{\rho_{i^*, \cdot} \in \R^{\Omega_{i^*}}}\min_{\substack{ \mu_{i^*, \cdot} \in \mathcal{M}^+(\Omega_{i^*} \times \Omega_\cdot)\\
        \mu_{i^*}\in \mathcal{M}_1^+(\Omega_{i^*})}} 
     \bigg(\sum_{(i^*,j) \in \mathcal{N}(i^*)}  \big\langle f_{i^*,j} , \mu_{i^*,j}\big\rangle - \varepsilon H(\mu_{i^*,j}) \bigg)
     \\
     &\qquad\qquad\qquad\qquad\qquad\qquad\qquad\qquad
     + \varepsilon(\abs{\mathcal{N}(i^*)}-1) \cdot H(\mu_{i^*})
     \\
     &\qquad\qquad\qquad\qquad\qquad\qquad\qquad\qquad+ \sum_{(i^*,j) \in \mathcal{N}(i^*) } \big\langle \rho_{i^*,j}, \mu_{i^*} - \mu_{i^*,j} \circ \pi_{i^*}^{-1} \big\rangle
\tag*{\small{(strong Lagrangian duality by \cref{thm:slater}; \cref{lem:concavity})}}
    \\
    &\quad=
    \displaystyle\max_{\rho_{i^*, \cdot} \in \R^{\Omega_{i^*}}}\min_{\substack{ \mu_{i^*, \cdot} \in \mathcal{M}^+(\Omega_{i^*} \times \Omega_\cdot)\\
        \mu_{i^*}\in \mathcal{M}_1^+(\Omega_{i^*})}} 
     \bigg\langle -\varepsilon (\abs{\mathcal{N}(i^*)}-1) \log\mu_{i^*} + \sum_{(i^*,j) \in \mathcal{N}(i^*) } \rho_{i^*, j}, \mu_{i^*} \bigg\rangle
     \\
     &\qquad\qquad\qquad\qquad\qquad\qquad\qquad\qquad+
     \sum_{(i^*, j ) \in \mathcal{N}(i^*)} \bigg\langle \varepsilon \log \mu_{i^*,j} + f_{i^*,j} -\rho_{i^*,j}, \mu_{i^*,j}\bigg\rangle
\tag*{\small{(definition entropy; linearity; constant integration)}}
    \\
    &\quad=
    \displaystyle\max_{\rho_{i^*, \cdot} \in \R^{\Omega_{i^*}}} \bigg(\min_{
        \mu_{i^*}\in \mathcal{M}_1^+(\Omega_{i^*})} 
     \bigg\langle -\varepsilon (\abs{\mathcal{N}(i^*)}-1) \log\mu_{i^*} + \sum_{(i^*,j) \in \mathcal{N}(i^*) } \rho_{i^*, j}, \mu_{i^*} \bigg\rangle \bigg)
     \\
     &\qquad\qquad\qquad\qquad\qquad+
     \bigg(\sum_{(i^*, j )\in \mathcal{N}(i^*)}  \min_{\mu_{i^*, j} \in \mathcal{M}^+(\Omega_{i^*} \times \Omega_j)}\bigg\langle \varepsilon \log \mu_{i^*,j} + f_{i^*,j} -\rho_{i^*,j}, \mu_{i^*,j}\bigg\rangle\bigg)
\tag*{\small{(minimize separable problems)}}
     \\
     &\quad=
     \displaystyle\max_{\rho_{i^*, \cdot} \in \R^{\Omega_{i^*}}}\bigg(\bigg(\min_{
        x_{i^*}\in  \Omega_{i^*}}\big(\textstyle\sum_{(i^*,j) \in \mathcal{N}(i^*) } \rho_{i^*, j}\big)(x_{i^*}) \bigg)
    \\
    &\qquad\qquad\qquad\qquad\qquad
    -\varepsilon\!\!\!\displaystyle\sum_{(i^*,j) \in \mathcal{N}(i^*) }\sum_{x_{i^*} \in \Omega_{i^*}} \sum_{x_j \in \Omega_j}
        \exp \bigg(\big(\rho_{i^*,j}(x_{i^*}) -f_{i^*,j}(x_{i^*},x_j)\big)/\varepsilon -1 \bigg)\bigg)
\tag*{\small{(by \cref{lem:marginalize_entropy1} \& \cref{lem:marginalize_entropy2})}}
     \\
     &\quad=
     \displaystyle\max_{\rho_{i^*, \cdot} \in \R^{\Omega_{i^*}}}\bigg(\bigg(\min_{
        x_{i^*}\in  \Omega_{i^*}}\big(\textstyle\sum_{(i^*,j) \in \mathcal{N}(i^*) } \rho_{i^*, j}\big)(x_{i^*}) \bigg)
    \\
    &\qquad\qquad\qquad\quad
    -\varepsilon\!\!\!\displaystyle\sum_{x_{i^*} \in \Omega_{i^*}} \sum_{(i^*,j) \in \mathcal{N}(i^*) }\exp\big(\rho_{i^*,j}(x_{i^*})/\varepsilon\big)\sum_{x_j \in \Omega_j}
        \exp \big(-f_{i^*,j}(x_{i^*},x_j)/\varepsilon -1 \big)\bigg)
\tag*{\small{(distributivity \& order of summation)}}
    \\
    &\quad=
    \begin{cases}
     \rho - \abs{\mathcal{N}(i^*)}
    \\
    \text{where }
    \\
    \rho_{i^*,j}(x_{i^*})
     = \mathrm{lse}_\varepsilon\big(f_{i^*,j}(x_{i^*}, \cdot)\big) + \tfrac{1}{\abs{\mathcal{N}(i^*)}} \bigg(\rho
     -\!\!\!\!\!\!\!\!\!\!\displaystyle\sum_{(i^*,k)\in \mathcal{N}(i^*)} \!\!\!\!\!\!\!\!\!\mathrm{lse}_\varepsilon\big(f_{i^*,k}(x_{i^*}, \cdot)\big)\bigg)\,,
     \\
     \qquad\qquad\forall x_{i^*} \in \Omega_{i^*} \,,\,\, (i^*, j) \in \mathcal{N}(i^*)
    \\
    \rho = \varepsilon\abs{\mathcal{N}(i^*)} +\varepsilon\abs{\mathcal{N}(i^*)}\log\big(\abs{\mathcal{N}(i^*)}/\varepsilon\big)
    \\
    \qquad
    - \abs{\mathcal{N}(i^*)}\mathrm{lse}_\varepsilon^{x_{i^*}, j}\bigg(\mathrm{lse}^{x_j}_\varepsilon\big(f_{i^*,j}(x_{i^*}, x_j)\big) +\mathrm{lse}^{x_j}_\varepsilon\big(-f_{i^*, j}(x_{i^*}, x_j)\big)\\
    \qquad\qquad\qquad\qquad\qquad
    -\tfrac{1}{\abs{\mathcal{N}(i^*)}}\!\!\!\!\!\!\!\!\!\!\displaystyle\sum_{(i^*,k)\in \mathcal{N}(i^*)} \!\!\!\!\!\!\!\!\!\mathrm{lse}^{x_k}_\varepsilon\big(f_{i^*,k}(x_{i^*}, x_k)\big)\bigg) \,.
    \end{cases}
\tag*{\small{(by \cref{lem:entropy_closed_form})}}
\end{align*}
\endgroup
\end{proof}

\subsection{Recovering Primal Solutions from Dual Solutions}
\label{subsec:dual2primal}%

Finally, we describe how to obtain a minimum of the sum of bivariates from a dual solution, e.g., in the problem formulations of \cref{thm:dual,thm:nscsdv,thm:entropydual,thm:duallinprog,rem:reprdetr}. Since we apply the method beyond star graphs, as opposed to much of the theory of this section, we prove that the result of \cref{thm:dual2primal} holds for tree-indexed sums of bivariates. However, the result will be heuristically applicable beyond tree-indexed sums of bivariates.\\

As we have not introduced the dual of the tree relaxation explicitly for this case yet---mainly because closed-form solutions do not exist in this generality---we forward-reference to the description in \cref{thm:duallinprog} to avoid redundancy.

Intuitively, the method in \cref{thm:dual2primal} determines a solution by sequentially minimizing a sum which fixes determined arguments in bivariates, ignores bivariates that are constant in the active argument, and replaces all bivariates with two undertermined arguments by their $\min$-marginals. 

\begin{theorem}[Recovering~Primal~Solutions]
\label{thm:dual2primal}%
    In the setting of \cref{thm:duallinprog}, let
    \begin{itemize}
    \item $\mathcal{G} = (\mathcal{V}, \mathcal{E})$ be an oriented tree,
    \item $\rho_{i,j} \in \R^{\Omega_i},\rho_{j,i} \in \R^{\Omega_j}; (i,j)\in \mathcal{E}$ be a solution to the \emph{dual linear program},
    \item $i_j \in \mathcal{V}$ denote the parent node for all non-root $j \in \mathcal{V}$, and%
    \vspace*{-1mm}%
    \[
        m_{i,j}(\cdot) \coloneqq \min_{x_j \in \Omega_j} f_{i,j}(\cdot, x_j) - \rho_{j,i}(x_j)
        \qquad \forall (i,j) \in \mathcal{E} \,.
    \]\vspace*{-1mm}%
    \end{itemize}%
    Then, we know that%
    \vspace*{-1mm}%
    \[
        x_j^* \in \argmin_{x_j \in \Omega_{i_t}}
        \begin{cases}
            \sum_{k \in \widetilde{\mathcal{N}}(j)\setminus\{i_j\}} m_{j, k}(x_{j})
            & \text{if $j$ is the root}
            \\
            f_{i_j, j}(x_{i_j}^*, x_{j}) + \sum_{k \in \widetilde{\mathcal{N}}(j)\setminus\{i_j\}}m_{j, k}(x_{j})
            & \text{else}
        \end{cases} 
    \]
    minimizes the sum of bivariates $F\equiv\sum_{(i,j)\in\mathcal{E}}f_{i,j}: \Omega \to \R$.
\end{theorem}
\begin{proof}
    Define $\rho \coloneqq \sum_{i\in \mathcal{V}} \rho_i \in \R$. By strong duality (similar proof as \cref{thm:dual}), we know $\min F - \rho = 0$ and that
    \[
        F-\rho = \sum_{(i,j) \in \mathcal{E}} f_{i,j} - \rho_{i,j} - \rho_{j,i} \,.
    \]
    By the constraint of the \emph{dual linear program} of \cref{thm:duallinprog}, we know that each summand $f_{i,j} - \rho_{i,j} - \rho_{j,i}; (i,j) \in \mathcal{E}$ is non-negative, and, therefore, that
    \[
        x^* \in \argmin_{x \in \Omega} F(x)-\rho \iff \bigg( f_{i,j}(x_i^*,x_j^*) - \rho_{i,j}(x_i^*) - \rho_{j,i}(x_j^*) = 0 \quad \forall (i,j) \in \mathcal{E} \bigg) \,.
    \]
    We inductively prove over the tree height that each summand is $0$ for our candidate $x^* \in \Omega$.
    \begin{enumerate}
        \item Let $i^*$ be the root. We show that
        \begin{align*}
            &\sum_{j \in \widetilde{\mathcal{N}}(i^*)} \min_{x_j \in \Omega_j} f_{i^*, j}(x_{i^*}^*, x_j) - \rho_{i^*, j}(x_{i^*}^*) - \rho_{j,i^*}(x_j)
            \\
            &\qquad
            =
            - \rho_{i^*} + \sum_{j \in \widetilde{\mathcal{N}}(i^*)} \min_{x_j \in \Omega_j} f_{i^*, j}(x_{i^*}^*, x_j) - \rho_{j,i^*}(x_j)
        \tag*{\small{($i_1$ is a root; definition of $\rho_{i^*}$ in \cref{thm:duallinprog})}}
            \\
            &\qquad
            =
            - \rho_{i^*} + \min_{x_{i^*} \in \Omega_{i^*}}\sum_{j \in \widetilde{\mathcal{N}}(i^*)} m_{i^*,j}(x_{i^*})
        \tag*{\small{(definition of $m_{i^*,j}$; definition of $x_{i^*}^*$)}}
            \\
            &\qquad=
            0 \,,
\tag*{\small{(by \cref{thm:nscsdv}; definition of $\rho_{i^*}$ in \cref{thm:duallinprog})}}
        \end{align*}
        Which implies that $x_j \in \Omega_j$ can be selected for all vertices $j$ of depth one, such that the the associated summand is zero, i.e.
        \[
        0 = \min_{x_j \in \Omega_j} f_{i^*, j}(x_{i^*}^*, x_j) - \rho_{i^*, j}(x_{i^*}^*) - \rho_{j,i^*}(x_j) \qquad \forall j \in \widetilde{\mathcal{N}}(i^*) \,.
        \]
        \item Given the presented method actually selects a successive value $x_j^* \in \Omega_j$ in this way, i.e.
        \[
            0 = f_{i, j}(x_{i}^*, x_{j}^*) - \rho_{i, j}(x_{i}^*) - \rho_{j,i}(x_j^*)
        \]
        for some $(i,j) \in \mathcal{E}$. We know that
        \begingroup%
        \allowdisplaybreaks%
        \begin{align*}
            &\sum_{k \in \widetilde{\mathcal{N}}(j)\setminus\{i\}} \min_{x_k \in \Omega_k} f_{j, k}(x_{j}^*, x_k) - \rho_{j, k}(x_j^*) - \rho_{k,j}(x_k)
            \\
            &\qquad=
             \bigg(f_{i, j}(x_{i}^*, x_{j}^*) - \rho_{i, j}(x_{i}^*) - \rho_{j,i}(x_j^*)\bigg)
             \\
             &\qquad\qquad
             + \bigg(\sum_{k \in \widetilde{\mathcal{N}}(j)\setminus\{i\}} \min_{x_k \in \Omega_k} f_{j, k}(x_{j}^*, x_k) - \rho_{j, k}(x_j^*) - \rho_{k,j}(x_k) \bigg)
        \tag*{\small{(assumption; adding $0$)}}
            \\
            &\qquad=
             \bigg(\min_{x_i \in \Omega_i} f_{i, j}(x_{i}, x_{j}^*) - \rho_{i, j}(x_{i}) - \rho_{j,i}(x_j^*)\bigg)
             \\
             &\qquad\qquad
             + \bigg(\!\sum_{k \in \widetilde{\mathcal{N}}(j)\setminus\{i\}} \min_{x_k \in \Omega_k} f_{j, k}(x_{j}^*, x_k) - \rho_{j, k}(x_j^*) - \rho_{k,j}(x_k) \bigg)
        \tag*{\small{(by $f_{i,j} - \rho_{i,j}-\rho_{j,i} \geq 0$)}}
            \\
            &\qquad=
             -\rho_j + \bigg(\min_{x_i \in \Omega_i} f_{i, j}(x_{i}, x_{j}^*) - \rho_{i, j}(x_{i})+ \sum_{k \in \widetilde{\mathcal{N}}(j)\setminus\{i\}} m_{j,k}(x_j^*)\bigg)
        \tag*{\small{(definition of $m_{j,k}$ and $\rho_j$ in \cref{thm:duallinprog})}}
            \\
            &\qquad=
             -\rho_j + \min_{x_j \in \Omega_i}\bigg( \min_{x_ji \in \Omega} f_{i, j}(x_{i}, x_{j}) - \rho_{i, j}(x_{i})+ \sum_{k \in \widetilde{\mathcal{N}}(j)\setminus\{i\}} m_{j,k}(x_j)\bigg)
        \tag*{\small{(definition of $x_j^*$; swapping minimization)}}
            \\
            &\qquad=
            0 \,.
    \tag*{\small{(by \cref{thm:nscsdv}; definition of $\rho_{i^*}$ in \cref{thm:duallinprog})}}
        \end{align*}
        \endgroup%
        This implies that all values $x_k \in \Omega_k$ of child vertices $k$ of $j$, can also be selected, such that their associated summand is zero, i.e.
        \[
            0 = \min_{x_k \in \Omega_k} f_{j, k}(x_{j}^*, x_k) - \rho_{j, k}(x_j^*) - \rho_{k,j}(x_k) \qquad \forall k \in \widetilde{\mathcal{N}}(j)\setminus\{i\} \,.
        \]\vspace*{-2mm}%
    \item The method, in fact, selects successive values such that summands are minimized to be zero, i.e.
    \vspace*{-2mm}%
    \begingroup%
    \allowdisplaybreaks%
    \begin{align*}
         & \argmin_{x_j \in \Omega_j} f_{i_j, j}(x_{i_j}^*, x_{j}) +\!\!\!\! \sum_{k \in \widetilde{\mathcal{N}}(j)\setminus\{i_j\}}m_{j, k}(x_{j})
         \\
         &\qquad=
         \argmin_{x_j \in \Omega_j} f_{i_j, j}(x_{i_j}^*, x_{j}) -\rho_{i_j,j}(x_i^*) - \rho_j + \!\!\!\!\sum_{k \in \widetilde{\mathcal{N}}(j)\setminus\{i_j\}} \min_{x_k \in \Omega_k} f_{j,k}(x_j, x_k) - \rho_{k,j}(x_k)
    \tag*{\small{(definition of $m_{j,k}$; $\rho_j, \rho_{i_j,j}(x_i^*) \in \R$)}}
         \\
         &\qquad=
         \argmin_{x_j \in \Omega_j} f_{i_j, j}(x_{i_j}^*, x_{j}) - \rho_{i_j,j}(x_i^*) - \rho_{j,i_j}(x_j)
         \\
         &\qquad\qquad\qquad+ \!\!\!\! \sum_{k \in \widetilde{\mathcal{N}}(j)\setminus\{i_j\}}\!\!\!\! \min_{x_k \in \Omega_k} f_{j,k}(x_j, x_k) - \rho_{j,k}(x_j) -\rho_{k,j}(x_k)
    \tag*{\small{(definition of $\rho_j$ in \cref{thm:duallinprog})}}
         \\
         &\qquad\subseteq
         \argmin_{x_j \in \Omega_j} f_{i_j, j}(x_{i_j}^*, x_{j}) - \rho_{i_j,j}(x_i^*) - \rho_{j,i_j}(x_j) \,.
    \tag*{\small{(previous result; non-negativity of summands)}}%
    \end{align*}%
    \endgroup%
    \end{enumerate}%
    \vspace{-3mm}\qedhere%
\end{proof}%

\clearpage%
\section{Algorithms}
\label{sec:algorithms}

We have developed the necessary results to derive several relevant types of optimization algorithms for sums of bivariates in this section. The optimization algorithms that we present differ in the objective function they operate on and their degree of approximating it.
For example, the first algorithm we introduce, coordinate descent (\hyperref[alg:cd]{\texttt{CD}}), operates directly on the sum of bivariates, which is the problem we ultimately aim to solve. All other algorithms, in contrast, operate on some form of dual relaxation of the sum of bivariates that we covered in the previous section.

While we derived insight into tractable optimization problems in the previous section, we will use these results to derive optimization algorithms for more general but similar problem formulations.

\subsection{Coordinate Descent for the Sum of Bivariates}
As a baseline algorithm for comparison in experiments, we define a simple coordinate descent procedure to heuristically minimize a sum of bivariates. Since we have shown the problem to be \text{NP}-hard, we do not expect efficiency nor convergence to the optimal solution, in general.\\

In contrast to the coordinate optimization approaches presented below, the following coordinate descent algorithm (\hyperref[alg:cd]{\texttt{CD}}) sequentially optimizes along single coordinates. The reason is that optimizing over neighborhoods of coordinates could yield an objective that is a partial sum of bivariates that does not have tree structure.

\begin{algorithm}[h]
\caption{Coordinate Descent for the Sum of Bivariates (\texttt{CD})%
}
\begin{algorithmic}[1]
    \INPUT vertex set $\mathcal{V} \coloneqq \N_{\leq n}$, $n \in \N$;
    \SUBINPUT edge set $\mathcal{E} \subseteq \{(i,j) \in \mathcal{V} \times \mathcal{V} \mid i < j\}$;
    \SUBINPUT candidate sets $\Omega_i \subset \R$, where $\abs{\Omega_i} \in \N$ and $i \in \mathcal{V}$;
    \SUBINPUT functions $f_{i,j}: \Omega_i \times \Omega_j \to \R$, where $(i,j) \in \mathcal{E}$;
    \SUBINPUT bijection $c: \mathcal{V} \to \N_{\leq n}$;
    \SUBINPUT budget $B \in \N$.
    \OUTPUT Element $x^* \in \Omega_1 \times \dots \times \Omega_n$ with \enquote{low} function value $\sum_{(i,j) \in \mathcal{E}} f_{i,j}(x_i, x_j)$.
    \INIT $\Tilde{\mathcal{E}} \coloneqq \big\{ \{ i,j\} \mid (i,j) \in \mathcal{E}\}$;\quad $\widetilde{\mathcal{N}}(i) =\{j \mid \{i, j\} \in \widetilde{\mathcal{E}}\},i \in \mathcal{V}$.
    \AlgBlankLine 
    \FORALL{$i \in (c^{-1}(1), \dots, c^{-1}(n))$}
    \STATE $x_i^* \in \argmin_{x_i\in\Omega_i} \sum_{j \in \widetilde{\mathcal{N}}(i)} f_{\min\{i,j\}, \max\{i,j\}}(x_{\min\{i,j\}}, x_{\max\{i,j\}})$
    \ENDFOR
    \IF{$B\leq t$} \algcomment{stopping criterion}
        \STATE \texttt{return} $x^*$
    \ENDIF
    \STATE $t \leftarrow t+1$
    \STATE \texttt{go to line 1}
    \algcomment{repeat the loop}
  \end{algorithmic}
\label{alg:cd}
\end{algorithm}

\subsection{Linear Programming for the Dual Linear Program}%
\label{subsec:duallinprog}%
Relaxation is the main alternative to an elementary approach that operates directly on the sum of bivariates. As in the previous section, due to strong duality and fewer variables, the dual of the relaxation replaces the relaxation.

The variables of the dual can be interpreted as constant sums of bivariates, such that each bivariate lower-bounds the respective bivariate of the objective function (note the constraint in \cref{thm:duallinprog}). In the case of \cref{thm:duallinprog}, we introduce the dual from a primal perspective, as it arguably better motivates constant univariate dual variables, namely as a type of local lower bound. Specifically, the problem we present in \cref{thm:duallinprog}, generalizes the dual relaxation of \cref{thm:dual} to non-star graphs, while \cref{rem:duallinprog} explains when it is still exact.

\begin{theorem}[Dual~Linear~Program]
\label{thm:duallinprog}
    In the setting of \cref{def:biv}, let a sum of bivariates $F\equiv\sum_{(i,j)\in\mathcal{E}}f_{i,j}: \Omega \to \R$ be given, then
        \begingroup%
    \allowdisplaybreaks%
    \begin{align*}%
&\begin{cases}
        \max_{P \in \R^\Omega} \min_{x \in \Omega} P(x)
        \\
        \text{\,\,s.t.\quad} P \text{ is separable}
        \\
        \phantom{\text{\,\,s.t.\quad}} P \equiv \sum_{(i,j) \in \mathcal{E}}  \tilde{\rho}_{i,j} 
        \\
        \phantom{\text{\,\,s.t.\quad}}\tilde{\rho}_{i,j} \in \R^{\Omega_i \times \Omega_j}
            \\
        \phantom{\text{\,\,s.t.\quad}}
            \tilde{\rho}_{i,j} \leq f_{i,j} \quad \forall (i,j) \in \mathcal{E}
        \end{cases}
    \tag*{\small{(Separable Local Lower Bound)}}
        \\
        &\qquad=
        \begin{cases}
        \max_{\rho_{\cdot,\cdot} \in \R^{\Omega_\cdot}}   \sum_{i \in \mathcal{V}} \min_{x_i\in \Omega_i}\sum_{j \in \widetilde{\mathcal{N}}(i)} \rho_{i,j}(x_i)
            \\        \text{\,\,s.t.\quad}
            \rho_{i,j} + \rho_{j,i} \leq f_{i,j} \quad \forall (i,j) \in \mathcal{E}
        \end{cases}
    \tag*{\small{(Dual Tree Relaxation)}}
        \\
        &\qquad=
                \begin{cases}
        \max_{\substack{\rho_{\cdot,\cdot} \in \R^{\Omega_\cdot} \\ \rho_\cdot \in \R}}  \sum_{i \in \mathcal{V}} \rho_i
            \\        \text{\,\,s.t.\quad}
            \rho_{i,j} + \rho_{j,i} \leq f_{i,j} \quad \forall (i,j) \in \mathcal{E}
            \\
            \phantom{\text{\,\,s.t.\quad}}
            \sum_{j \in \widetilde{\mathcal{N}}(i)} \rho_{i,j} \equiv \rho_i \quad \forall i \in \mathcal{V} \,.
        \end{cases}
    \tag*{\small{(Dual Linear Program)}}
    \end{align*}%
    \endgroup%
\end{theorem}%
\begin{remark}%
\label{rem:duallinprog}
    By \cref{thm:relaxation} and \cref{thm:dual}, we know that for star graphs the minimal value of the Dual Linear Program is the minimal value of the sum of bivariates. By \cref{ex:counterrecon}, we know that this does not apply beyond trees in general.
\end{remark}
\begin{proof}
We prove the equality of the first and third term
        \begingroup%
    \allowdisplaybreaks%
    \begin{align*}%
        &\begin{cases}
        \max_{P \in \R^\Omega} \min_{x \in \Omega} P(x)
        \\
        \text{\,\,s.t.\quad} P \text{ is separable}
        \\
        \phantom{\text{\,\,s.t.\quad}} P \equiv \sum_{(i,j) \in \mathcal{E}}  \tilde{\rho}_{i,j} 
        \\
        \phantom{\text{\,\,s.t.\quad}}\tilde{\rho}_{i,j} \in \R^{\Omega_i \times \Omega_j}
            \\
        \phantom{\text{\,\,s.t.\quad}}
            \tilde{\rho}_{i,j} \leq f_{i,j} \quad \forall (i,j) \in \mathcal{E}
        \end{cases}
        \\
        &\qquad=
                \begin{cases}
        \max_{P \in \R^\Omega} \min_{x \in \Omega} P(x)
        \\
        \text{\,\,s.t.\quad} P \text{ is constant}
        \\
        \phantom{\text{\,\,s.t.\quad}} P \equiv \sum_{(i,j) \in \mathcal{E}}  \tilde{\rho}_{i,j}
            \\
        \phantom{\text{\,\,s.t.\quad}}
        \tilde{\rho}_{i,j} \in \R^{\Omega_i \times \Omega_j}
            \\        \phantom{\text{\,\,s.t.\quad}}
            \tilde{\rho}_{i,j} \leq f_{i,j} \quad \forall (i,j) \in \mathcal{E}
        \end{cases}
    \tag*{\small{(by projection $P \mapsto ( x \mapsto \min P)$)}}
        \\
        &\qquad=
                \begin{cases}
        \max_{P \in \R^\Omega} \min_{x \in \Omega} P(x)
        \\
        \text{\,\,s.t.\quad} P \text{ is constant}
        \\
        \phantom{\text{\,\,s.t.\quad}} P \equiv \sum_{(i,j) \in \mathcal{E}}  \rho_{i,j} + \rho_{j,i}
            \\
        \phantom{\text{\,\,s.t.\quad}}
        \rho_{i,j} \in \R^{\Omega_i} \,,\,\, \rho_{j,i} \in \R^{\Omega_j}
            \\        \phantom{\text{\,\,s.t.\quad}}
            \rho_{i,j} + \rho_{j,i} \leq f_{i,j} \quad \forall (i,j) \in \mathcal{E}
        \end{cases}
    \tag*{\small{(\cref{cor:dv})}}
        \\
        &\qquad=
                \begin{cases}
        \max_{P \in \R^\Omega} \min_{x \in \Omega} \sum_{i \in \mathcal{V}} \sum_{j \in \widetilde{\mathcal{N}}(i)} \rho_{i,j}(x_i)
        \\
        \text{\,\,s.t.\quad} P \text{ is constant}
        \\
        \phantom{\text{\,\,s.t.\quad}} P \equiv \sum_{(i,j) \in \mathcal{E}}  \rho_{i,j} + \rho_{j,i}
            \\
        \phantom{\text{\,\,s.t.\quad}}
        \rho_{i,j} \in \R^{\Omega_i} \,,\,\, \rho_{j,i} \in \R^{\Omega_j}
            \\        \phantom{\text{\,\,s.t.\quad}}
            \rho_{i,j} + \rho_{j,i} \leq f_{i,j} \quad \forall (i,j) \in \mathcal{E}
        \end{cases}
    \tag*{\small{(rewrite objective)}}
        \\
        &\qquad=
                \begin{cases}
        \max_{P \in \R^\Omega}  \sum_{i \in \mathcal{V}} \min_{x_i \in \Omega_i}\sum_{j \in \widetilde{\mathcal{N}}(i)} \rho_{i,j}(x_i)
        \\
        \text{\,\,s.t.\quad} P \text{ is constant}
        \\
        \phantom{\text{\,\,s.t.\quad}} P \equiv \sum_{(i,j) \in \mathcal{E}}  \rho_{i,j} + \rho_{j,i}
            \\
        \phantom{\text{\,\,s.t.\quad}}
        \rho_{i,j} \in \R^{\Omega_i} \,,\,\, \rho_{j,i} \in \R^{\Omega_j} \,,\,\, \rho_i, \rho_j \in \R
            \\        \phantom{\text{\,\,s.t.\quad}}
            \rho_{i,j} + \rho_{j,i} \leq f_{i,j} \quad \forall (i,j) \in \mathcal{E}
            \\
            \phantom{\text{\,\,s.t.\quad}}
            \sum_{j \in \widetilde{\mathcal{N}}(i)} \rho_{i,j} \equiv \rho_i \quad \forall i \in \mathcal{V}
        \end{cases}
    \tag*{\small{($P$ constant $\iff$ $\sum_{j \in \widetilde{\mathcal{N}}(i)} \rho_{i,j}$ constant, $i \in \mathcal{V}$)}}
        \\
        &\qquad=
                \begin{cases}
        \max_{\substack{\rho_{\cdot,\cdot} \in \R^{\Omega_\cdot} \\ \rho_\cdot \in \R}}  \sum_{i \in \mathcal{V}} \rho_i
            \\        \text{\,\,s.t.\quad}
            \rho_{i,j} + \rho_{j,i} \leq f_{i,j} \quad \forall (i,j) \in \mathcal{E}
            \\
            \phantom{\text{\,\,s.t.\quad}}
            \sum_{j \in \widetilde{\mathcal{N}}(i)} \rho_{i,j} \equiv \rho_i \quad \forall i \in \mathcal{V} \,.
        \end{cases}
    \tag*{\small{(eliminate variables and rewrite objective)}}
    \end{align*}%
    \endgroup%
    The first equality becomes clear, by replacing $\rho_i, i \in \mathcal{V}$ in the objective and optimizing over non-constant functions $\sum_{j \in \widetilde{\mathcal{N}}(i)} \rho_{i,j}$ instead.
\end{proof}

In algorithm \hyperref[alg:lpdlp]{\texttt{LPDLP}}, we minimize a sum of bivariates by using the dual linear program formulation of \cref{thm:duallinprog} to generate a dual solution and then deriving a solution candidate to the sum of bivariates using the method described in \cref{thm:dual2primal}. The algorithm is proven to be exact for tree-indexed sums of bivariates, as the tree relaxation of \cref{thm:relaxation} is exact, its dual is strong via a similar argument as in \cref{thm:dual}, and the recovery of a primal solution is exact via \cref{thm:dual2primal}.

\begin{algorithm}[h]
\caption{Linear Programming for the Dual Linear Program (\texttt{LPDLP})%
}
\begin{algorithmic}[1]
    \INPUT vertex set $\mathcal{V} \coloneqq \N_{\leq n}$, $n \in \N$;
    \SUBINPUT edge set $\mathcal{E} \subseteq \{(i,j) \in \mathcal{V} \times \mathcal{V} \mid i < j\}$;
    \SUBINPUT candidate sets $\Omega_i \subset \R$, where $\abs{\Omega_i} \in \N$ and $i \in \mathcal{V}$;
    \SUBINPUT functions $f_{i,j}: \Omega_i \times \Omega_j \to \R$, where $(i,j) \in \mathcal{E}$.
    \SUBINPUT bijection $c: \mathcal{V} \to \N_{\leq n}$;
    \OUTPUT Element $x^* \in \Omega_1 \times \dots \times \Omega_n$ with \enquote{low} function value $\sum_{(i,j) \in \mathcal{E}} f_{i,j}(x_i, x_j)$.
    \INIT $\Tilde{\mathcal{E}} \coloneqq \big\{ \{ i,j\} \mid (i,j) \in \mathcal{E}\}$;\quad $\widetilde{\mathcal{N}}(i) =\{j \mid \{i, j\} \in \widetilde{\mathcal{E}}\},i \in \mathcal{V}$;
    \SUBINPUT \quad\,\, $\rho_{i,j}, m_{i,j} \coloneqq 0 \in \R^{\Omega_i}, \rho_{j,i},m_{j,i} \coloneqq 0 \in \R^{\Omega_j}, (i,j) \in \mathcal{E}$;\quad $\rho_i \coloneqq 0 \in \R, i \in \mathcal{V}$.
    \AlgBlankLine 
    \STATE \texttt{solve for}
    \algcomment{solving the \emph{dual linear program}; see \cref{thm:duallinprog}}
        \[\begin{cases}
        (\rho^*_{\cdot, \cdot}, \rho^*_\cdot) \in \argmax_{\substack{\rho_{\cdot,\cdot} \in \R^{\Omega_\cdot} \\ \rho_\cdot \in \R}}  \sum_{i \in \mathcal{V}} \rho_i
            \\        \text{\,\,s.t.\quad}
            \rho_{i,j}(x_i) + \rho_{j,i}(x_j) \leq f_{i,j}(x_i, x_j) \quad \forall x_i \in \Omega_i \,\, \forall x_j \in \Omega_j \,\, \forall (i,j) \in \mathcal{E}
            \\
            \phantom{\text{\,\,s.t.\quad}}
            \sum_{j \in \widetilde{\mathcal{N}}(i)} \rho_{i,j}(x_i) = \rho_i \quad \forall x_i \in \Omega_i \,\, \forall i \in \mathcal{V}
        \end{cases}\]
        \STATE $m_{i,j} \leftarrow \min_{x_j \in \Omega_j} f_{i,j}(\cdot, x_j) - \rho^*_{j,i}(x_j), (i,j) \in \mathcal{E}$
         \algcomment{dual to primal; see \cref{thm:dual2primal}}
        \STATE $m_{j,i} \leftarrow \min_{x_i \in \Omega_i} f_{i,j}(x_i, \cdot) - \rho^*_{i,j}(x_i), (i,j) \in \mathcal{E}$
        \FOR{$i \in (c^{-1}(1), \dots, c^{-1}(n))$}
        \STATE $x^*_i \!\!\in \!\argmin_{x_i \in \Omega_i} \!\sum_{\!\!\!\substack{\{i,j\} \in \widetilde{\mathcal{E}}\\c(j)<c(i)}} \!\!\! f_{\min\{i,j\},\max\{i,j\}}(x_{\min\{i,j\}}, x_{\max\{i,j\}}) + \!\sum_{\substack{\{i,j\}\in\widetilde{\mathcal{E}}\\c(i)<c(j)}} \!\!m_{i,j}(x_i)$
        \ENDFOR
    \STATE \texttt{return} $x^*$
  \end{algorithmic}
\label{alg:lpdlp}
\end{algorithm}

\subsection{Block Coordinate Ascent for the Dual Tree Relaxation}
\label{subsec:coorddesc}
Consider restricting the minimization of the dual tree relexation from \cref{thm:duallinprog} to $\rho_{i^*,j} \in \R^{\Omega_i}, j \in \widetilde{\mathcal{N}}(i^*)$ for some coordinate $i^* \in \mathcal{V}$. This yields a new problem on what is called a \emph{block coordinate} of the objective function. We have
\begin{align*}
    &\begin{cases}
        \max_{\rho_{i^*,\cdot} \in \R^{\Omega_{i^*}}}   \sum_{i \in \mathcal{V}} \min_{x_i\in \Omega_i}\sum_{j \in \widetilde{\mathcal{N}}(i)} \rho_{i,j}(x_i)
            \\        \text{\,\,s.t.\quad}
            \rho_{i,j} + \rho_{j,i} \leq f_{i,j} \quad \forall (i,j) \in \mathcal{E}
    \end{cases}
    \\
    &\qquad=
    \begin{cases}
        \max_{\rho_{i^*,\cdot} \in \R^{\Omega_{i^*}}}   \min_{x_{i^*}\in \Omega_{i^*}}\sum_{j \in \widetilde{\mathcal{N}}(i^*)} \rho_{i^*,j}(x_{i^*})
            \\        \text{\,\,s.t.\quad}
            \rho_{i^*,j} \leq \min_{x_j \in \Omega_j} f_{i^*,j}(\cdot, x_j) - \rho_{j,i^*}(x_j) \quad \forall j \in \widetilde{\mathcal{N}}(i^*) \,.
    \end{cases}
\tag*{\small{(possibly transposing $f_{\cdot,\cdot}$; dropping constant trerms)}}
\end{align*}

 However, by \cref{thm:nscsdv}, we know how to solve this problem. Repeatedly replacing variables by their block coordinate optima in some given order yields our version of what is called \emph{(block) coordinate ascent} in \hyperref[alg:bcadtr]{\texttt{BCADTR}}.\\
 
 Some authors provide counterexamples even for convergence to an optimal solution of the dual linear program of \cref{thm:duallinprog} for methods of this class \cite[Sec.~4.4]{Dlask2022}. Others make similar claims \cite{werner2007,prusa2016,tourani2018,savchynskyy2019discrete}. Evidence for the hypothesis that \hyperref[alg:bcadtr]{\texttt{BCADTR}} does not, in general, converge to an optimal solution of the dual linear program can also be observed in \cref{sec:firstexp}. To obtain a solution candidate for the minimization of the sum of bivariates, we apply the method described in \cref{thm:dual2primal} once again.

\begin{algorithm}[h]
\caption{Block Coordinate Ascent for the Dual Tree Relaxation (\texttt{BCADTR})%
}
\begin{algorithmic}[1]
    \INPUT vertex set $\mathcal{V} \coloneqq \N_{\leq n}$, $n \in \N$;
    \SUBINPUT edge set $\mathcal{E} \subseteq \{(i,j) \in \mathcal{V} \times \mathcal{V} \mid i < j\}$;
    \SUBINPUT candidate sets $\Omega_i \subset \R$, where $\abs{\Omega_i} \in \N$ and $i \in \mathcal{V}$;
    \SUBINPUT functions $f_{i,j}: \Omega_i \times \Omega_j \to \R$, where $(i,j) \in \mathcal{E}$;
    \SUBINPUT bijection $c: \mathcal{V} \to \N_{\leq n}$;
    \SUBINPUT weights $w_{i,j;t} \in \R$, $j \in \widetilde{\mathcal{N}}(i)$, $i \in \mathcal{V}$ with $\sum_{j \in \widetilde{\mathcal{N}}(i)} w_{i,j; t} =1$, $t \in \N_{\leq B}$;
    \SUBINPUT budget $B \in \N$.
    \OUTPUT Element $x^* \in \Omega_1 \times \dots \times \Omega_n$ with \enquote{low} function value $\sum_{(i,j) \in \mathcal{E}} f_{i,j}(x_i, x_j)$.
    \INIT $t\coloneqq 1 \in \N$;\quad $\widetilde{\mathcal{E}} \coloneqq \big\{ \{ i,j\} \mid (i,j) \in \mathcal{E}\}$;
    \SUBINPUT \quad\,\, $\widetilde{\mathcal{N}}(i) =\{j \mid \{i, j\} \in \widetilde{\mathcal{E}}\}, m_i \coloneqq 0 \in \R^{\Omega_i},i \in \mathcal{V}$;
    \SUBINPUT \quad\,\, $m_{i,j},\rho_{i,j} \coloneqq 0 \in \R^{\Omega_i}, m_{j,i},\rho_{j,i} \coloneqq 0 \in \R^{\Omega_j}, (i,j) \in \mathcal{E}$;
    \SUBINPUT \quad\,\, $F \equiv \sum_{(i,j) \in \mathcal{E}}f_{i,j}$;
    \SUBINPUT \quad\,\, $x^* \in \Omega_1\times \dots \times \Omega_n$ uniformly at random.
    \AlgBlankLine
    \FORALL{$i \in (c^{-1}(1), \dots, c^{-1}(n))$}           \algcomment{one loop in order defined by $c$}
        \FORALL{$\ell \in \widetilde{\mathcal{N}}(i)$}\algcomment{block coordinate maximization; see \cref{thm:nscsdv}}
        \STATE $m_{i,j}(x_i) \leftarrow \min_{x_j \in \Omega_j} f_{\min\{i,j\}, \max\{i,j\}}(x_{\min\{i,j\}}, x_{\max\{i,j\}}) - \rho_{j,i}(x_j), \forall x_i \in \Omega_i$
        \ENDFOR
        \STATE $m_i \leftarrow \sum_{j \in \widetilde{\mathcal{N}}(i)} m_{i,j}$
        \STATE $\rho_{i,j}  
     \leftarrow m_{i,j} - w_{i,j;t} \big(m_i - \min_{x_i \in \Omega_i} m_i(x_i)\big), j \in \widetilde{\mathcal{N}}(i)$
    \ENDFOR
    \FOR{$i \in (c^{-1}(1), \dots, c^{-1}(n))$}
         \algcomment{dual to primal; see \cref{thm:dual2primal}}
        \STATE $y^*_i \!\!\in \!\argmin_{x_i \in \Omega_i} \!\sum_{\!\!\!\substack{\{i,j\} \in \widetilde{\mathcal{E}}\\c(j)<c(i)}} \!\!\! f_{\min\{i,j\},\max\{i,j\}}(x_{\min\{i,j\}}, x_{\max\{i,j\}}) + \!\sum_{\substack{\{i,j\}\in\widetilde{\mathcal{E}}\\c(i)<c(j)}} \!\!m_{i,j}(x_i)$
    \ENDFOR
    \IF{$F(y^*)<F(x^*)$} \algcomment{keep best solution candidate}
        \STATE $x^* \leftarrow y^*$
    \ENDIF
    \IF{$B\leq t$} \algcomment{stopping criterion}
        \STATE \texttt{return} $x^*$
    \ENDIF
    \STATE $t \leftarrow t+1$
    \STATE \texttt{go to line 1}
    \algcomment{repeat the loop}
  \end{algorithmic}
\label{alg:bcadtr}
\end{algorithm}

\subsection{Block Coordinate Ascent for the Dual Entropy-Regularized Tree Relaxation}

This algorithm is inspired by the entropy-based algorithm of \cite{ochs2024}. We showed by \cref{ex:hamilton,ex:counterrecon} that, in case the graph $\mathcal
 {G} = (\mathcal{V}, \mathcal{E})$ is not a tree,
 translating solutions of the \emph{tree relaxation} presented in \cref{thm:relaxation} to solutions to the minimization of sums of bivariates is \text{NP}-hard. However, as previously noted, we have developed a method in \cref{thm:dual2primal} that can still be applied heuristically. Further, using linear programming, we can already obtain a dual solution---see \cref{subsec:duallinprog} for details. However, one may expect to use the closed-form solution on block coordinates of \cref{subsec:relaxation} to derive a different, possibly a computationally more efficient, solver. We attempted this already in \cref{subsec:coorddesc}\\
 
 We thus propose a block coordinate ascent for the dual of a generalization of the entropy-regularized star relaxation of \cref{thm:entropydual} to arbitrary graphs, termed \emph{entropy-regularized tree relaxation}. This generalization is similarly regularized by an entropy function and approximates tree-structured sums of bivariates.

 The main differences to the block coordinate ascent solver for the dual tree relaxation of \cref{subsec:coorddesc} lie in
 \begin{enumerate}
     \item the fact that the objective of the following solver only results in an approximation of the dual of the tree relaxation, and 
     \item that we are able to prove convergence of the solver to a solution with an objective value that is arbitrarily close to the optimal value of the dual tree relaxation. 
 \end{enumerate}
 
We present the algorithm in \hyperref[alg:bcadetr]{\texttt{BCADETR}}. Note that it is not obvious whether better solution candidates of the dual yield better solution candidates to the sum of bivariates, in general. Therefore, similarly as for \hyperref[alg:bcadtr]{\texttt{BCADTR}}, it is not clear in which iteration the method obtains the best solution to the original problem of minimizing a sum of bivariates.
 
\begin{remark}[Representation of Dual Entropy-Regularized Tree Relaxation]
\label{rem:reprdetr}%
    In the setting of \cref{def:biv}, let a sum of bivariates $F\equiv\sum_{(i,j)\in\mathcal{E}}f_{i,j}: \Omega \to \R$ and $\varepsilon>0$ be given. Then, by a similar argument as in \cref{thm:entropydual,lem:marginalize_entropy2}, we have
        \begingroup%
    \allowdisplaybreaks%
    \begin{align*}
    &\begin{cases}
    \displaystyle\min_{\substack{\mu_{\cdot, \cdot} \in \mathcal{M}^+(\Omega_\cdot \times \Omega_\cdot)\\ \mu_{\cdot} \in \mathcal{M}_1^+(\Omega_{\cdot})}} \bigg(\sum_{(i,j) \in \mathcal{E}} \big\langle f_{i,j} , \mu_{i,j}\big\rangle - \varepsilon H(\mu_{i,j})\bigg)+ \varepsilon \sum_{i \in \mathcal{V}}(\abs{\widetilde{\mathcal{N}}(i)}-1) \cdot H(\mu_{i}) 
        \\
        \,\,\text{s.t.}\quad
        \mu_i = \mu_{i,j} \circ \pi_i^{-1} \,,\quad \forall (i,j) \in \mathcal{E}
        \\
        \phantom{\,\,\text{s.t.}\quad}
        \mu_j = \mu_{i,j} \circ \pi_j^{-1} \,,\quad \forall (i,j) \in \mathcal{E} 
        \end{cases}%
\tag*{\small{(entropy-regularized tree relaxation)}}
    \\
    &\quad=
    \begin{cases}
    \displaystyle\max_{\rho_\cdot \in \R} \max_{\rho_{\cdot, \cdot} \in \R^{\Omega_\cdot}}  \bigg(\sum_{i \in \mathcal{V}} \rho_i \bigg)  -\varepsilon\!\!\!\displaystyle\sum_{(i,j) \in \mathcal{E} } \sum_{\substack{x_i \in \Omega_i\\x_j \in \Omega_j}}
        \exp \big(\big(\rho_{i,j}(x_{i})+\rho_{j,i}(x_{j})-f_{i,j}(x_{i},x_j)\big)/\varepsilon -1 \big)
    \\
    \text{s.t. }\big(\textstyle\sum_{j \in \widetilde{\mathcal{N}}(i)} \rho_{i, j}\big)(x_{i}) = \rho_i \,,\,\, \forall x_{i} \in \Omega_{i} \,,\,\, i \in \mathcal{V}\,.
    \end{cases}
\tag*{\small{(dual entropy-regularized tree relaxation)}}
\end{align*}%
\endgroup%
\end{remark}%
\begin{remark}%
\label{rem:entropyapprox}
    Due to boundedness of the entropy $H$, it is clear that for a given sum of bivariates, the \emph{entropy-regularized tree relaxation} uniformly approximates the \emph{tree relaxation} of \cref{thm:relaxation}. This holds, of course, in particular, if the graph associated with the sum of bivariates is a star graph, as in \cref{thm:entropydual}.
\end{remark}%

\begin{corollary}%
\label{cor:convergence}%
    The objective value of a block coordinate ascent w.r.t.~the blocks $\rho_{i,j}, \rho_i; j \in \widetilde{\mathcal{N}}(i)$ for all $i \in \mathcal{V}$ converges to the global maximal value of the \emph{dual entropy-regularized tree relaxation} of \cref{rem:reprdetr}.
\end{corollary}
\begin{proof}
    Clearly, due to the convexity of $\exp$, the \emph{dual entropy-regularized tree relaxation} is a sum of concave functions in $\rho_\cdot \in \R,\rho_{\cdot, \cdot} \in \R^{\Omega_\cdot}$. The objective function is smooth, and the domain of the function is a block-wise product space. Furthermore, by \cref{thm:entropydual}, block coordinate-wise minimization yields a unique minimum. Therefore, we can apply \cref{thm:bertsekas}.
\end{proof}

\begin{algorithm}[h]
\caption{Block Coordinate Ascent for the\\ Dual Entropy-Regularized Tree Relaxation (\texttt{BCADETR})%
}
\begin{algorithmic}[1]
    \INPUT vertex set $\mathcal{V} \coloneqq \N_{\leq n}$, $n \in \N$;
    \SUBINPUT edge set $\mathcal{E} \subseteq \{(i,j) \in \mathcal{V} \times \mathcal{V} \mid i < j\}$;
    \SUBINPUT candidate sets $\Omega_i \subset \R$, where $\abs{\Omega_i} \in \N$ and $i \in \mathcal{V}$;
    \SUBINPUT functions $f_{i,j}: \Omega_i \times \Omega_j \to \R$, where $(i,j) \in \mathcal{E}$;
    \SUBINPUT entropy regularization coefficient $\varepsilon \in (0, \infty)$;
    \SUBINPUT bijection $c: \mathcal{V} \to \N_{\leq n}$;
    \SUBINPUT budget $B \in \N$.
    \OUTPUT Element $x^* \in \Omega_1 \times \dots \times \Omega_n$ with \enquote{low} function value $\sum_{(i,j) \in \mathcal{E}} f_{i,j}(x_i, x_j)$.
    \INIT $t \coloneqq 1\in \R$;\quad $\widetilde{\mathcal{E}} \coloneqq \big\{ \{ i,j\} \mid (i,j) \in \mathcal{E}\}$;
    \SUBINPUT \quad\,\, $\widetilde{\mathcal{N}}(i) =\{j \mid \{i, j\} \in \widetilde{\mathcal{E}}\}$,
    \SUBINPUT \quad\,\, $\gamma_{i} \coloneqq \varepsilon\abs{\widetilde{\mathcal{N}}(i)} +\varepsilon\abs{\widetilde{\mathcal{N}}(i)}\log\big(\abs{\widetilde{\mathcal{N}}(i)}/\varepsilon\big)$, $i \in \mathcal{V}$;
    \SUBINPUT \quad\,\, $\rho_{i,j}, m_{i,j}, \underline{m}_{i,j}, \overline{m}_{i,j} \coloneqq 0 \in \R^{\Omega_i}, \rho_{j,i}, m_{j,i}, \underline{m}_{j,i}, \overline{m}_{j,i} \coloneqq 0 \in \R^{\Omega_j}, (i,j) \in \mathcal{E}$;
    \SUBINPUT \quad\,\, $F \equiv \sum_{(i,j) \in \mathcal{E}}f_{i,j}$;
    \SUBINPUT \quad\,\, $x^* \in \Omega_1\times \dots \times \Omega_n$ uniformly at random.
    \AlgBlankLine
    \FORALL{$i \in (c^{-1}(1), \dots, c^{-1}(n))$}           \algcomment{one loop in order defined by $c$}
        \FORALL{$\ell \in \widetilde{\mathcal{N}}(i)$}
        \algcomment{block coordinate maximization; see \cref{thm:entropydual}}
        \STATE $\overline{m}_{i,\ell}(x_i) \leftarrow \mathrm{lse}^{x_\ell}_\varepsilon\big(f_{\min\{i,\ell\}, \max\{i,\ell\}}(x_{\min\{i,\ell\}}, x_{\max\{i,\ell\}})-\rho_{\ell, i}(x_\ell)\big), \forall x_i \in \Omega_i$
        \STATE $\underline{m}_{i,\ell}(x_i) \leftarrow -\mathrm{lse}^{x_\ell}_\varepsilon\big(\rho_{\ell, i}(x_\ell)-f_{\min\{i,\ell\}, \max\{i,\ell\}}(x_{\min\{i,\ell\}}, x_{\max\{i,\ell\}})\big), \forall x_i \in \Omega_i$
        \ENDFOR
        \STATE $\rho \leftarrow  \gamma_i
        - \abs{\widetilde{\mathcal{N}}(i)}\mathrm{lse}_\varepsilon^{x_{i}, \ell}\bigg(\overline{m}_{i,\ell}(x_i) -\underline{m}_{i,\ell}(x_i)-\tfrac{1}{\abs{\widetilde{\mathcal{N}}(i)}}\!\!\!\!\!\!\!\displaystyle\sum_{(i,k)\in \widetilde{\mathcal{N}}(i)}\!\!\!\!\!\! \overline{m}_{i,k}(x_i)\bigg)$
        \STATE $\rho_{i,j}
        \leftarrow \overline{m}_{i,j} + \tfrac{1}{\abs{\widetilde{\mathcal{N}}(i)}} \bigg(\rho
        -\!\!\!\!\!\!\!\displaystyle\sum_{(i,k)\in \mathcal{N}(i)}\!\!\!\!\!\! \overline{m}_{i,k}\bigg), \forall j \in \widetilde{\mathcal{N}}(i)$
    \ENDFOR
    \STATE $m_{i,j} \leftarrow \min_{x_j \in \Omega_j} f_{i,j}(\cdot, x_j) - \rho_{j,i}(x_j), (i,j) \in \mathcal{E}$
         \algcomment{dual to primal; see \cref{thm:dual2primal}}
    \STATE $m_{j,i} \leftarrow \min_{x_i \in \Omega_i} f_{i,j}(x_i, \cdot) - \rho_{i,j}(x_i), (i,j) \in \mathcal{E}$
    \FOR{$i \in (c^{-1}(1), \dots, c^{-1}(n))$}
        \STATE $y^*_i \!\!\in \!\argmin_{x_i \in \Omega_i} \!\sum_{\!\!\!\substack{\{i,j\} \in \widetilde{\mathcal{E}}\\c(j)<c(i)}} \!\!\! f_{\min\{i,j\},\max\{i,j\}}(x_{\min\{i,j\}}, x_{\max\{i,j\}}) + \!\sum_{\substack{\{i,j\}\in\widetilde{\mathcal{E}}\\c(i)<c(j)}} \!\!m_{i,j}(x_i)$
    \ENDFOR
    \IF{$F(y^*)<F(x^*)$} \algcomment{keep best solution candidate}
        \STATE $x^* \leftarrow y^*$
    \ENDIF
    \IF{$B\leq t$} \algcomment{stopping criterion}
        \STATE \texttt{return} $x^*$
    \ENDIF
    \STATE $t \leftarrow t+1$
    \STATE \texttt{go to line 1}
    \algcomment{repeat the loop}
  \end{algorithmic}
\label{alg:bcadetr}
\end{algorithm}

\clearpage
\section{Experiments}
\label{sec:experiments}%

In this section, we compare implementations of the algorithms \hyperref[alg:cd]{\texttt{CD}}, \hyperref[alg:lpdlp]{\texttt{LPDLP}}, \hyperref[alg:bcadtr]{\texttt{BCADTR}}, \hyperref[alg:trw-s]{\texttt{TRW-S}} and \hyperref[alg:trw-s-leg]{\texttt{TRW-S-LEG}} for the optimization of several qualitatively distinct and practically relevant objective function classes.

The algorithms \hyperref[alg:trw-s]{\texttt{TRW-S}} and \hyperref[alg:trw-s-leg]{\texttt{TRW-S-LEG}} are considered state-of-the-art methods for the optimization of sums of bivariates. The experiments include two versions of \hyperref[alg:bcadtr]{\texttt{BCADTR}}. Namely one, where the weights are sampled uniformly on the simplex for each coordinate and each time step independently, termed \enquote{random}, and another version, where the weights are constant for each coordinate and each time step, termed \enquote{constant}. The numerically stable implementation of \hyperref[alg:bcadetr]{\texttt{BCADETR}} is considered to be beyond the scope of this work.\\

The analysis of the value of dual tree relaxation is restricted to \hyperref[alg:lpdlp]{\texttt{LPDLP}} and both versions of \hyperref[alg:bcadtr]{\texttt{BCADTR}} due to their consistent parameterization.

All algorithms are available in the Python programming language and are implemented in a vectorized, sparse as well as in-place manner, whenever beneficial.\footnote{An implementation can be found at \href{https://github.com/NiMlr/pybiv}{\texttt{https://github.com/NiMlr/pybiv}}.} The implementation of \hyperref[alg:lpdlp]{\texttt{LPDLP}} is based on the \emph{HiGHS solver} for linear programs \cite{huangfu2018}. For all iterative algorithms that operate on the objective, intermediate results are shown.

In addition, \emph{wall-clock time} is chosen for comparing the algorithms as it is a common index.

\subsection{Random Sums of Bivariates}
\label{sec:firstexp}
\begin{figure}[ht]
\begin{center}%
\hspace*{-.15\textwidth}\includegraphics[width=1.3\textwidth]{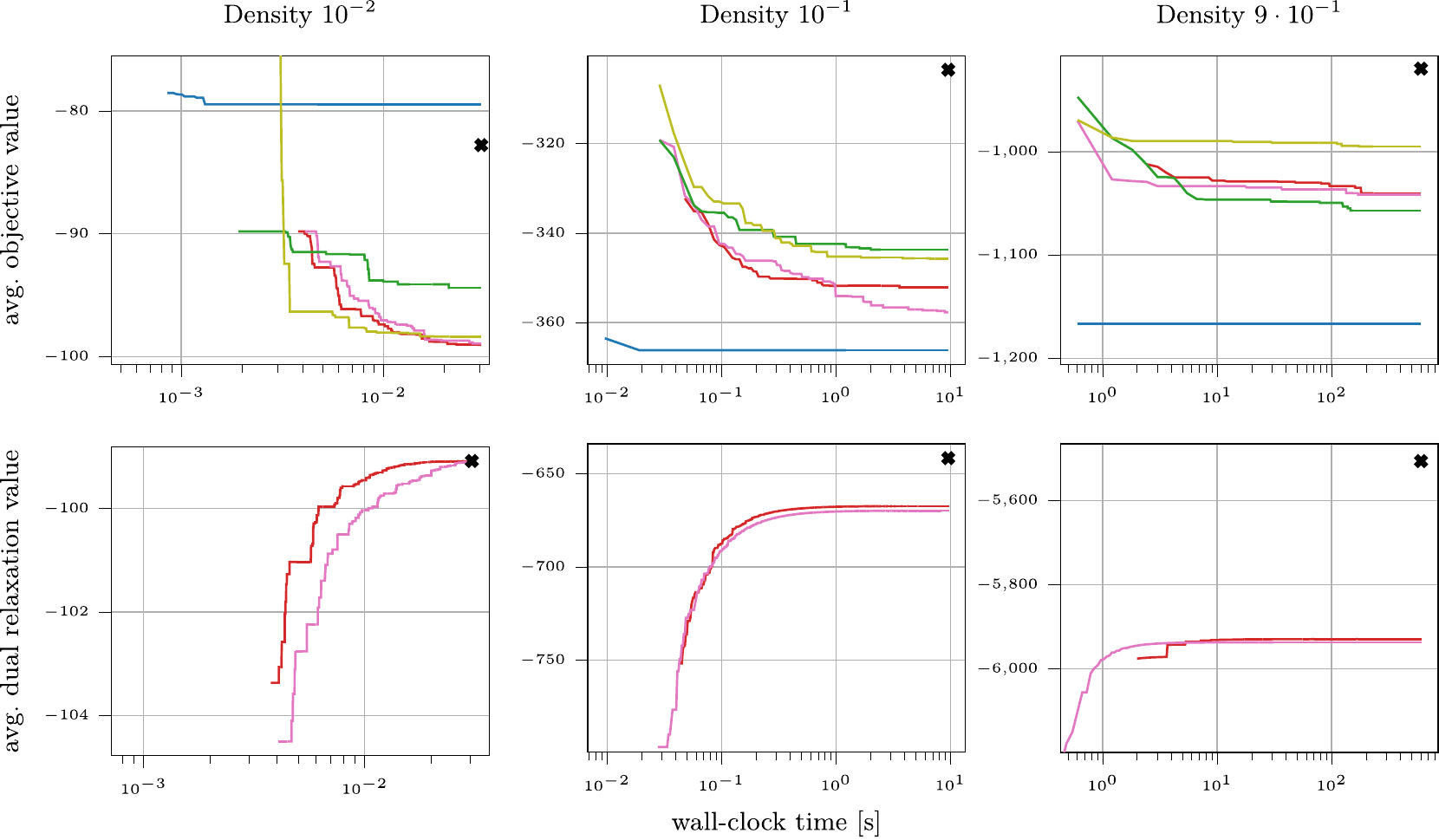}
\end{center}%
    \vspace*{-6mm}    \begin{center}\small{\raisebox{-0.165em}{\scalebox{2.0}{$\bullet$}} \hyperref[alg:lpdlp]{\texttt{LPDLP}}\quad \textcolor{steelblue31119180}{\raisebox{-0.175em}{\scalebox{2.0}{$\bullet$}}} \hyperref[alg:cd]{\texttt{CD}}\quad \textcolor{crimson2143940}{\raisebox{-0.175em}{\scalebox{2.0}{$\bullet$}}} \hyperref[alg:bcadtr]{\texttt{BCADTR}} (constant $w_\cdot$)\quad \textcolor{orchid227119194}{\raisebox{-0.175em}{\scalebox{2.0}{$\bullet$}}} \hyperref[alg:bcadtr]{\texttt{BCADTR}} (random $w_\cdot$)\quad \textcolor{forestgreen4416044}{\raisebox{-0.175em}{\scalebox{2.0}{$\bullet$}}} \hyperref[alg:trw-s]{\texttt{TRW-S}}\quad \textcolor{goldenrod18818934}{\raisebox{-0.175em}{\scalebox{2.0}{$\bullet$}}} \hyperref[alg:trw-s-leg]{\texttt{TRW-S-LEG}}}\end{center}%
    \caption{\textbf{Minimization of random sums of bivariates.} Benchmark of several algorithms for the optimization of sums of bivariates with $100$ arguments that can each take $5-15$ values for various densities of the Erd\H{o}s--R\'enyi random graph that indexes the sum of bivariates. The plots represent the average of $10$ independent runs and of the best value found after given wall-clock time. Primal and dual relaxed values are shown. See \cref{sec:firstexp} for details.}
    \label{fig:random-functions}
\end{figure}

First, we construct random sums of bivariates. As usual, we are in the setting of \cref{def:biv}, and we assume the function values of all bivariates to be independently normal distributed. In addition, the edges of the graph that indexes the sum of bivariates follow the distribution of the $G(n, m)$-model of Erd\H{o}s--R\'enyi random graphs. Precisely, we have
\begin{align*}
    &
    \Omega_i \coloneqq \{1, \dots, K_i \} \,,\,\, K_i \sim \mathrm{Unif}\big(\{5, \dots, 15\}\big) \quad \forall i \in \N_{\leq n} \,,
    \\
    &F \equiv \textstyle\sum_{(i,j) \in \mathcal{E}} f_{i,j} \,,
    \\
    &f_{i,j}(x_i, x_j) \sim \mathcal{N}(0,1) \text{ independently} \quad \forall x_i \in \Omega_i\,\, \forall x_j \in \Omega_j \,\, \forall (i,j) \in \mathcal{E} \,, \,\, \text{and}
    \\
    &\mathcal{E} \sim \mathrm{Unif}\big(\{ \mathcal{H} \subseteq \mathcal{J} \mid \abs{\mathcal{H}} = m\}\big) \,,\,\, \text{where} \,\, \mathcal{J} \coloneqq \{(i,j) \in \mathcal{V} \times \mathcal{V} \mid i < j\} \,, m \in \N \,.
\end{align*}
For the experiment, we select the number of arguments of the problem as $n=100$, the bijection $c: \mathcal{V} \to \N_{\leq n}$ uniformly at random and the number of edges $m$, such that the density $m/\abs{\mathcal{J}}$ of the graph $(\mathcal{V}, \mathcal{E})$ approximates the values $10^{-2}, 10^{-1}$ and $9\cdot10^{-1}$. We run the algorithms on $10$ independent samples of the objective function $F \equiv \textstyle\sum_{(i,j) \in \mathcal{E}} f_{i,j}$. The results are visualized in \cref{fig:random-functions}.\\

An interpretation of the results of the first experiment is that sufficiently dense random sums of bivariates contain structure that is hard to exploit for algorithms based on dual relaxations. In particular, even a simple method, such as (primal) coordinate descent (\hyperref[alg:cd]{\texttt{CD}}) outperforms all other algorithms on sufficiently dense instances. Additionally, random sums of bivariates seem to lie outside the design domain of of the legacy tree-reweighted message passing (\hyperref[alg:trw-s-leg]{\texttt{TRW-S-LEG}}), as it is outperformed by the block-coordinate ascent algorithm for the dual tree relaxation (\hyperref[alg:bcadtr]{\texttt{BCADTR}}).

Solving the dual tree relaxation more precisely does not in general yield better primal candidates, as can be seen by the exact solution to the dual linear program (\hyperref[alg:lpdlp]{\texttt{LPDLP}}). Interestingly, for sufficiently dense random instances, the dual block coordinate ascent algorithm (\hyperref[alg:bcadtr]{\texttt{BCADTR}}) does not converge to the maximal value of the dual linear program. It can also be seen that the quality of the primal candidate derived from the solution to the dual linear program deteriorates for dense instances at least as much as the quality of the best candidates of all other algorithms.

\subsection{Vertex Coloring}
\label{sec:secondexp}
For a given graph $\mathcal{G} = (\mathcal{E}, \mathcal{V}$), vertex coloring asks to determine whether we can assign one of $n_c$ colors for each vertex, such that neighboring vertices are assigned different colors. Mathematically, the goal is to determine the existence of a function $\mathrm{color}: \mathcal{V} \to \{1, \dots, n_c\}, n_c \in \N$, such that $\mathrm{color}(i) \neq \mathrm{color}(j)$ for all edges $(i,j) \in \mathcal{E}$.

The problem of determining whether a $n_c$-coloring with $n_c=4$ exists for the vertices of a graph is \text{NP}-complete. In the next experiment, we attempt to find a $4$-coloring for a given random graph with low density by minimizing the number of color defects between neighboring vertices.

\begin{figure}[ht]
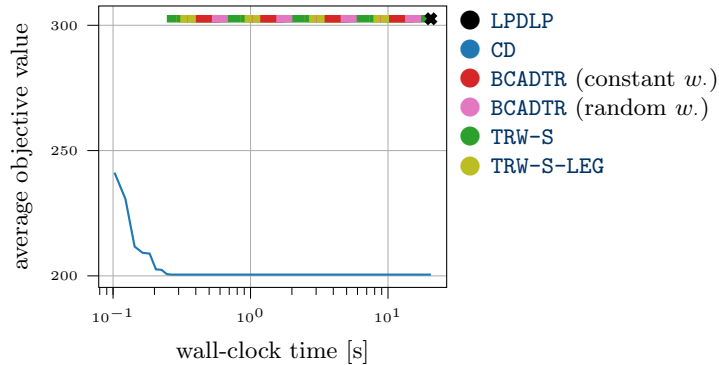

\begin{center}%
\begin{minipage}[t]{.4\textwidth}%
    \include{files/cf-den-001-prim}%
\end{minipage}%
\hspace*{4mm}%
\begin{minipage}[t]{.25\textwidth}%
\vspace*{1mm}
\small%
\raisebox{-0.175em}{\scalebox{2.0}{$\bullet$}} \hyperref[alg:lpdlp]{\texttt{LPDLP}}
\\
\textcolor{steelblue31119180}{\raisebox{-0.175em}{\scalebox{2.0}{$\bullet$}}} \hyperref[alg:cd]{\texttt{CD}}
\\
\textcolor{crimson2143940}{\raisebox{-0.175em}{\scalebox{2.0}{$\bullet$}}} \hyperref[alg:bcadtr]{\texttt{BCADTR}} (constant $w_\cdot$)
\\
\textcolor{orchid227119194}{\raisebox{-0.175em}{\scalebox{2.0}{$\bullet$}}} \hyperref[alg:bcadtr]{\texttt{BCADTR}} (random $w_\cdot$)
\\
\textcolor{forestgreen4416044}{\raisebox{-0.175em}{\scalebox{2.0}{$\bullet$}}} \hyperref[alg:trw-s]{\texttt{TRW-S}}
\\
\textcolor{goldenrod18818934}{\raisebox{-0.175em}{\scalebox{2.0}{$\bullet$}}} \hyperref[alg:trw-s-leg]{\texttt{TRW-S-LEG}}
\end{minipage}%
\end{center}%
\vspace*{-6.5mm}
\caption{\textbf{Minimization of the number of defects in vertex coloring.} Benckmark of several algorithms for the optimization of sums of bivariates with $10^{3}$ arguments. Each argument that represent one of $4$ colors in a vertex coloring of a given Erd\H{o}s--R\'enyi random graph with density $10^{-2}$ that indexes the sum of bivariates. The plots represent the average of $10$ independent runs and of the best value found after given wall-clock time. See \cref{sec:secondexp} for details.}%
\label{fig:vertex-coloring}%
\end{figure}%

Precisely, we have in the setting of \cref{def:biv} that
\begin{align*}
    &
    \Omega_i \coloneqq \{1, 2, 3, 4\} \,,
    \\
    &F \equiv \textstyle\sum_{(i,j) \in \mathcal{E}} f_{i,j} \,,
    \\
    &f_{i,j}(x_i, x_j) 
    =
    \begin{cases}
        1 & \text{if } x_i = x_j
        \\
        0 & \text{else}
    \end{cases}
    \,\,\, \forall x_i \in \Omega_i\,\, \forall x_j \in \Omega_j \,\, \forall (i,j) \in \mathcal{E} \,, \,\, \text{and}
    \\
    &\mathcal{E} \sim \mathrm{Unif}\big(\{ \mathcal{H} \subseteq \mathcal{J} \mid \abs{\mathcal{H}} = m\}\big) \,,\,\, \text{where} \,\, \mathcal{J} \coloneqq \{(i,j) \in \mathcal{V} \times \mathcal{V} \mid i < j\} \,, m \in \N \,.
\end{align*}

Clearly, if we determine $x \in \Omega$ such that $F(x) = 0$, then $x$ encodes a $4$-coloring of the graph $\mathcal{G}$.

For the experiment, we select the number of arguments of the problem as $n=1000$, the bijection $c: \mathcal{V} \to \N_{\leq n}$ uniformly at random, and the number of edges $m$, such that the density $m/\abs{\mathcal{J}}$ of the graph $(\mathcal{V}, \mathcal{E})$ approximates $10^{-2}$. We run the algorithms on $10$ independent samples of the objective function $F \equiv \textstyle\sum_{(i,j) \in \mathcal{E}} f_{i,j}$. The results are visualized in \cref{fig:vertex-coloring}.\\

An interpretation of the results of the second experiment is that all algorithms that are based on (tree) relaxations can only determine trivial solutions if the structure of the bivariates is sufficiently adverse. This holds even when the density of the graph is comparatively low. Only coordinate descent (\hyperref[alg:cd]{\texttt{CD}}) progresses by sequentially selecting each vertex color such that it has the least number of defects. However, we must conclude that all considered algorithms are not even heuristically competitive for the problem class of vertex coloring.

\subsection{Signal Reconstruction}
\label{sec:thirdexp}
The challenge of signal reconstruction describes a setting in which a mathematical function, termed \emph{signal}, is superimposed with (random) noise and where the goal is to determine the function without noise or a good approximation of it, termed \emph{reconstruction}. For our purposes, we consider a binary signal $\mathrm{sig}: \Z_{n} \to \{0,1\}$ defined on a discretized circle, the quotient group $\Z_n$. The noisy version of the signal is given by $\mathrm{sig}_{\mathrm{noisy}}: \Z_{n} \to \R$, where $\mathrm{sig}_{\mathrm{noisy}}(x) = \mathrm{sig}(x) + \nu_x$ and $\nu_x \sim \mathrm{Unif}([-1, 1])$ independently for all $x \in \Z_n$.

\begin{figure}[ht]
\begin{center}%
\begin{minipage}[t]{.4\textwidth}%
    \include{files/sf-noi-200-len-25000-prim}%
\end{minipage}%
\hspace*{4mm}%
\begin{minipage}[t]{.25\textwidth}%
\vspace*{4mm}
\small%
\raisebox{-0.175em}{\scalebox{2.0}{$\bullet$}} \hyperref[alg:lpdlp]{\texttt{LPDLP}}
\\
\textcolor{steelblue31119180}{\raisebox{-0.175em}{\scalebox{2.0}{$\bullet$}}} \hyperref[alg:cd]{\texttt{CD}}
\\
\textcolor{crimson2143940}{\raisebox{-0.175em}{\scalebox{2.0}{$\bullet$}}} \hyperref[alg:bcadtr]{\texttt{BCADTR}} (constant $w_\cdot$)
\\
\textcolor{orchid227119194}{\raisebox{-0.175em}{\scalebox{2.0}{$\bullet$}}} \hyperref[alg:bcadtr]{\texttt{BCADTR}} (random $w_\cdot$)
\\
\textcolor{forestgreen4416044}{\raisebox{-0.175em}{\scalebox{2.0}{$\bullet$}}} \hyperref[alg:trw-s]{\texttt{TRW-S}}
\\
\textcolor{goldenrod18818934}{\raisebox{-0.175em}{\scalebox{2.0}{$\bullet$}}} \hyperref[alg:trw-s-leg]{\texttt{TRW-S-LEG}}
\end{minipage}\\
\begin{minipage}[t]{.4\textwidth}%
    \vspace{-10mm}%
    \include{files/sf-noi-200-len-25000-dual}%
\end{minipage}%
\hspace*{4mm}%
\begin{minipage}[t]{.25\textwidth}%
\hfill%
\end{minipage}%
\end{center}%
\vspace*{-6.5mm}
\caption{\textbf{Minimization of the error in a (regularized) signal reconstruction problem.} Benckmark of several algorithms for the optimization of sums of bivariates with $25\cdot 10^3$ arguments. Each argument represents one of $2$ values the reconstruction of a binary signal can take. The graph that indexes the sum of bivariates representing the error has an approximate density of $8\cdot10^{-5}$. The plots represent the average of $10$ independent runs and of the best value found after given wall-clock time. Primal and dual relaxed values are shown. See \cref{sec:thirdexp} for details.}%
\label{fig:signal-reconstruction}%
\end{figure}%

Given the noisy version of the signal, our reconstruction of the signal is determined by minimizing a heuristic \emph{error function}, which commonly happens to have a representation as a sum of bivariates. The error function penalizes the reconstruction both when deviating from the noisy signal and when neighboring values are different. Precisely, we have in the setting of \cref{def:biv} that
\begin{align*}
    &
    \Omega_i \coloneqq \{0, 1\} \,,
    \\
    &F \equiv \textstyle\sum_{(i,j) \in \mathcal{E}} f_{i,j} \,,
    \\
    &f_{i,j}(x_i, x_j) 
    = \abs{\mathrm{sig_\mathrm{noisy}(x_i) - x_i}}
    \\
&\qquad\qquad\qquad\qquad+ (1/2)\cdot
    \begin{cases}
        0 & \text{if } x_i = x_j
        \\
        1 & \text{else}
    \end{cases}
    \quad \forall x_i \in \Omega_i\,\, \forall x_j \in \Omega_j \,\, \forall (i,j) \in \mathcal{E} \,, \,\, \text{and}
    \\
    &\mathcal{E} = \{ (x,x+1) \mid x \in \Z_n \}\,.
\end{align*}
For the experiment, we select the number of arguments of the problem as $n=25\cdot10^{3}$, the signal $\mathrm{sig}$ is sampled uniformly at random and the bijection $c: \mathcal{V} \to \N_{\leq n}$ as $c(i) = i$ for all $i \in \mathcal{V}$. The density is approximately $8\cdot10^{-5}$. We run the algorithms on $10$ independent samples of the objective function $F \equiv \textstyle\sum_{(i,j) \in \mathcal{E}} f_{i,j}$. The results are visualized in \cref{fig:signal-reconstruction}.\\

An interpretation of the results of the third experiment is that all algorithms that are based on (tree) relaxations perform well in reconstructing the signal. The fact that the relaxation gap is rather small and that all such methods outperform \hyperref[alg:cd]{\texttt{CD}} significantly suggests that the performance of the algorithms is non-trivial. Legacy tree-reweighted message passing (\hyperref[alg:trw-s-leg]{\texttt{TRW-S-LEG}}) works especially well---arguably because signal reconstruction methods are the design domain of this method. This class of problems seems to be the only class in which explicitly solving the dual linear program (\hyperref[alg:lpdlp]{\texttt{LPDLP}}) is competitive. In addition, the block-coordinate ascent methods converge to the maximum of the dual linear program.
\clearpage
\section{Conclusion and Future Work}
\label{sec:conc}%
We provided a comprehensive introduction to relaxation theory for the optimization of sums of bivariates that resulted in closed-form, linear programming, and block coordinate ascent solution approaches. The work was motivated by a fundamental perspective on both hard as well as easy instances of sums of bivariates, and a distributional perspective on the function class using a (no-)free-lunch theorem.

The theory featured insightful results on the limits of relaxation-based approaches, as well as fundamental restrictions on the reconstruction of global functions from their bivariate marginals, which are intimately related to the theory of sums of bivariates. A central technique when working with bivariates was the design of methods that not only provide the correct result for tree-structured sums of bivariates, but are applicable outside this narrow domain. We continued to apply this principle in entropy calculations and regularization.\\

Promising future work extends the theory of tractable subclasses of the sums of bivariates. This could be done, e.g., via a theory of separable lower bounds to bivariates as well as via approximation results that extend our analysis. Although suggested by the complexity results and the relaxation theory of this work, our experiments provide further evidence for this claim.

Generally, it is plausible that a theory of sums of bivariates shall precede a theory of the minimization thereof. It is also plausible that such generic structures as sums of bivariates are useful in other areas of applied mathematics.

Alternatives to relaxation-based optimization algorithms may be found via integer linear programming formulations of sums of bivariates. 
\clearpage
\printbibliography[heading=bibintoc]
\clearpage
\appendix
\section{Additional Proofs}

\subsection{{Lemmata for \cref{thm:entropydual}}}

\begin{lemma}
\label{lem:marginalize_entropy1}
    In the setting of \cref{def:biv} and \cref{thm:entropydual}, we have that
    \begin{align*}
        \min_{
        \mu_{i^*}\in \mathcal{M}_1^+(\Omega_{i^*})} 
     \bigg\langle -\varepsilon (\abs{\mathcal{N}(i^*)}-1) \log\mu_{i^*} + \!\!\!\!\!\!\!\sum_{(i^*,j) \in \mathcal{N}(i^*) } \!\!\!\!\rho_{i^*, j}, \mu_{i^*} \bigg\rangle = \min_{
        x_{i^*}\in  \Omega_{i^*}}\bigg(\sum_{(i^*,j) \in \mathcal{N}(i^*) } \!\!\!\!\rho_{i^*, j}\bigg)(x_{i^*}) \,,
    \end{align*}
    where $\varepsilon > 0$ and $\rho_{i^*,j} \in \R^{\Omega_{i^*}}, (i^*,j)\in\mathcal{N}(i^*)$.
\end{lemma}
\begin{proof}
    W.l.o.g., we assume that $\abs{\mathcal{N}(i^*)}-1 > 0$. Otherwise, the statement is clearly correct. Define
    \begin{align*}
        G(\mu_{i^*})&:=\bigg\langle -\varepsilon (\abs{\mathcal{N}(i^*)}-1) \log\mu_{i^*} + \sum_{(i^*,j) \in \mathcal{N}(i^*) } \rho_{i^*, j}, \mu_{i^*} \bigg\rangle
    \\
        &\phantom{:}
        =
        \sum_{x_{i^*} \in \Omega_{i^*}}\mu_{i^*}(x_{i^*})\bigg( -\varepsilon (\abs{\mathcal{N}(i^*)}-1) \log\mu_{i^*}(x_{i^*}) + \sum_{(i^*,j) \in \mathcal{N}(i^*) } \rho_{i^*, j}(x_{i^*})\bigg)
    \tag*{\small{(definition scalar product)}}
    \end{align*}
    for all $\mu_{i^*} \in \mathcal{M}_1^+(\Omega_{i^*})$\,. We have for all $x_{i^*},y_{i^*}\in \Omega_{i^*}$ and all strictly positive measures $\mu_{i^*} \in \mathcal{M}_1^+(\Omega_{i^*})$ that
    \begin{align*}
        &\frac{\partial}{\partial \mu_{i^*}(y_{i^*})}\frac{\partial}{\partial \mu_{i^*}(x_{i^*})} G(\mu_{i^*})
        \\
        &\qquad=
        \frac{\partial}{\partial \mu_{i^*}(y_{i^*})}\frac{\partial}{\partial \mu_{i^*}(x_{i^*})}\mu_{i^*}(x_{i^*})\bigg( -\varepsilon (\abs{\mathcal{N}(i^*)}-1) \log\mu_{i^*}(x_{i^*}) + \sum_{(i^*,j) \in \mathcal{N}(i^*) } \rho_{i^*, j}(x_{i^*})\bigg)
    \tag*{\small{(definition $G$)}}
        \\
        &\qquad=
        \frac{\partial}{\partial \mu_{i^*}(y_{i^*})}-\varepsilon (\abs{\mathcal{N}(i^*)}-1) \big(\log\mu_{i^*}(x_{i^*})+1\big) + \sum_{(i^*,j) \in \mathcal{N}(i^*) } \rho_{i^*, j}(x_{i^*})
    \tag*{\small{(product rule)}}
        \\
        &\qquad=
        \chi_{x_{i^*}=y_{i^*}}
        \frac{-\varepsilon(\abs{\mathcal{N}(i^*)}-1)}{\mu_{i^*}(x_{i^*})} \,.
    \tag*{\small{($\chi$ denotes indicator function)}}
    \end{align*}
    Therefore, the Hessian of $G$ is negative definite and $G$ is strictly concave. This implies that the minimum of $G$ lies on the boundary of $\mathcal{M}_1^+(\Omega_{i^*})$. Which in turn implies that at least one atom $x_{i^*} \in \Omega_{i^*}$ of the minimum has zero measure. Fixing $\mu_{i*}(x_{i^*})=0$ and minimizing $G$ over the remaining atoms yields a problem with the same properties, recursively. Therefore, the minimum of $G$ must be a Dirac measure.
    Thus, we have
    \begingroup
    \allowdisplaybreaks
    \begin{align*}
        &\min_{
        \mu_{i^*}\in \mathcal{M}_1^+(\Omega_{i^*})} 
     \bigg\langle -\varepsilon (\abs{\mathcal{N}(i^*)}-1) \log\mu_{i^*} + \sum_{(i^*,j) \in \mathcal{N}(i^*) } \rho_{i^*, j}, \mu_{i^*} \bigg\rangle
     \\
        &\qquad=
        \min_{
        \mu_{i^*}\in  \mathcal{D}(\Omega_{i^*})} 
     \bigg\langle -\varepsilon (\abs{\mathcal{N}(i^*)}-1) \log\mu_{i^*} + \sum_{(i^*,j) \in \mathcal{N}(i^*) } \rho_{i^*, j}, \mu_{i^*} \bigg\rangle
     \tag*{\small{($\mathcal{D}(\Omega_{i^*})$ denotes the Dirac measures on $\Omega_{i^*}$)}}
     \\
        &\qquad=
        \min_{
        \mu_{i^*}\in  \mathcal{D}(\Omega_{i^*})}\sum_{x_{i^*} \in \Omega_{i^*}}\mu_{i^*}(x_{i^*})\bigg( -\varepsilon (\abs{\mathcal{N}(i^*)}-1) \log\mu_{i^*}(x_{i^*}) + \sum_{(i^*,j) \in \mathcal{N}(i^*) } \rho_{i^*, j}(x_{i^*})\bigg)
    \tag*{\small{(definition scalar product)}}
     \\
        &\qquad=
        \min_{
        \mu_{i^*}\in  \mathcal{D}(\Omega_{i^*})}\sum_{x_{i^*} \in \Omega_{i^*}}\chi_{\mu_{i^*}(x_{i^*})=1}\sum_{(i^*,j) \in \mathcal{N}(i^*) } \rho_{i^*, j}(x_{i^*})
    \tag*{\small{($\chi$ denotes indicator function)}}
     \\
        &\qquad=
        \min_{
        x_{i^*}\in  \Omega_{i^*}}\bigg(\sum_{(i^*,j) \in \mathcal{N}(i^*) } \rho_{i^*, j}\bigg)(x_{i^*}) \,.
        \qedhere
    \end{align*}
    \endgroup
\end{proof}

\begin{lemma}
\label{lem:marginalize_entropy2}
In the setting of \cref{def:biv} and \cref{thm:entropydual}, we have for all $(i^*,j) \in \mathcal{N}(i^*)$ that
    \begin{align*}
        &\min_{\mu_{i^*, j} \in \mathcal{M}^+(\Omega_{i^*} \times \Omega_j)}\bigg\langle \varepsilon \log \mu_{i^*,j} + f_{i^*,j} -\rho_{i^*,j}, \mu_{i^*,j}\bigg\rangle
        \\
        &\qquad\qquad=
        -\varepsilon\sum_{x_{i^*} \in \Omega_{i^*}} \sum_{x_j \in \Omega_j}
        \exp \bigg(\big(\rho_{i^*,j}(x_{i^*}) -f_{i^*,j}(x_i,x_j)\big)/\varepsilon -1 \bigg)\,,
    \end{align*}
    where $\varepsilon > 0$ and $\rho_{i^*,j} \in \R^{\Omega_{i^*}}$.
\end{lemma}
\begin{proof}
We have
    \begin{align*}
        &\min_{\mu_{i^*, j} \in \mathcal{M}^+(\Omega_{i^*} \times \Omega_j)}\bigg\langle \varepsilon \log \mu_{i^*,j} + f_{i^*,j} -\rho_{i^*,j}, \mu_{i^*,j}\bigg\rangle
        \\
        &\,=\!\!\!\!\!
        \min_{\mu_{i^*, j} \in \mathcal{M}^+(\Omega_{i^*} \times \Omega_j)}
        \sum_{x_{i^*} \in \Omega_{i^*}} \!\!\sum_{x_j \in \Omega_j} \mu_{i^*,j}(x_{i^*},x_j)\bigg( \varepsilon \log \mu_{i^*,j}(x_{i^*},x_j) + f_{i^*,j}(x_{i^*},x_j) -\rho_{i^*,j}(x_{i^*}) \bigg)
        \\
        &\,=\!\!\!\!\!\!
        \sum_{x_{i^*} \in \Omega_{i^*}} \sum_{x_j \in \Omega_j}\!
        \min_{\mu_{i^*, j}(x_{i^*}, x_j) \in [0, \infty)}\!
        \mu_{i^*,j}(x_{i^*}, x_j)\bigg( \varepsilon \log \mu_{i^*,j}(x_{i^*}, x_j) + f_{i^*,j}(x_{i^*},x_j) -\rho_{i^*,j}(x_{i^*}) \bigg).
    \end{align*}
    It can be observed that for all $x_{i^*} \in \Omega_{i^*}$ and $x_j \in \Omega_j$, the minimizer of each subproblem will never be at $\infty$, since the objective convergences to $\infty$ for $\mu_{i^*,j}(x_{i^*}, x_j) \to \infty$. Further, at $\mu_{i^*,j}(x_{i^*}, x_j) \downarrow 0$ the objective has the value $0$.\\
    Therefore, we check a first order criterion to find critical points with $\mu_{i^*,j}(x_{i^*}, x_j) > 0$. We have for such critical points that
    \begin{align*}
        &0 \overset{!}{=}\frac{\mathrm{d}}{\mathrm{d}\mu_{i^*,j}(x_{i^*}, x_j)}\mu_{i^*,j}(x_{i^*}, x_j)\bigg( \varepsilon \log \mu_{i^*,j}(x_{i^*}, x_j) + f_{i^*,j}(x_{i^*},x_j) -\rho_{i^*,j}(x_{i^*}) \bigg)
        \\
        &\qquad\qquad=
        \varepsilon \big( \log \mu_{i^*,j}(x_{i^*}, x_j) + 1\big) +  f_{i^*,j}(x_{i^*},x_j) -\rho_{i^*,j}(x_{i^*})
    \tag*{\small{(product rule)}}
    \\
    &\iff
    \mu_{i^*,j}(x_{i^*}, x_j) = \exp \bigg(\big(\rho_{i^*,j}(x_{i^*}) -f_{i^*,j}(x_{i^*},x_j)\big)/\varepsilon -1 \bigg) \,.
    \end{align*}
    Plugging the critical point into the objective, we see that the associated objective value is
    \begin{align*}
        &\exp \bigg(\big(\rho_{i^*,j}(x_{i^*}) -f_{i^*,j}(x_{i^*},x_j)\big)/\varepsilon -1 \bigg)\bigg( \varepsilon \log \exp \bigg(\big(\rho_{i^*,j}(x_{i^*}) -f_{i^*,j}(x_{i^*},x_j)\big)/\varepsilon -1 \bigg)
        \\
        &\qquad\qquad
        + f_{i^*,j}(x_{i^*},x_j) -\rho_{i^*,j}(x_{i^*}) \bigg)
        \\
        &\qquad=
        -\exp \bigg(\big(\rho_{i^*,j}(x_{i^*}) -f_{i^*,j}(x_{i^*},x_j)\big)/\varepsilon -1 \bigg)\varepsilon \,.
    \end{align*}
    The critical point, which has a negative function value, must be a minimizer, as at both boundaries the function has non-negative values and there are no further critical points.\\
    Plugging this representation of the minimal value into the initial equation, we get
    \begin{align*}
        &\min_{\mu_{i^*, j} \in \mathcal{M}^+(\Omega_{i^*} \times \Omega_j)}\bigg\langle \varepsilon \log \mu_{i^*,j} + f_{i^*,j} -\rho_{i^*,j}, \mu_{i^*,j}\bigg\rangle
        \\
        &\qquad\qquad=
        -\varepsilon\sum_{x_{i^*} \in \Omega_{i^*}} \sum_{x_j \in \Omega_j}
        \exp \bigg(\big(\rho_{i^*,j}(x_{i^*}) -f_{i^*,j}(x_{i^*},x_j)\big)/\varepsilon -1 \bigg)\,.
        \qedhere
    \end{align*}%
\end{proof}%

\begin{lemma}
\label{lem:entropy_closed_form}
In the setting of \cref{def:biv} and \cref{thm:entropydual}, we have
\begin{align*}
     &\displaystyle\max_{\rho_{i^*, \cdot} \in \R^{\Omega_{i^*}}}\bigg(\bigg(\min_{
        x_{i^*}\in  \Omega_{i^*}}\big(\textstyle\sum_{(i^*,j) \in \mathcal{N}(i^*) } \rho_{i^*, j}\big)(x_{i^*}) \bigg)
    \\
    &\qquad\qquad\qquad\quad
    -\varepsilon\!\!\!\displaystyle\sum_{x_{i^*} \in \Omega_{i^*}} \sum_{(i^*,j) \in \mathcal{N}(i^*) }\exp\big(\rho_{i^*,j}(x_{i^*})/\varepsilon\big)\sum_{x_j \in \Omega_j}
        \exp \big(-f_{i^*,j}(x_{i^*},x_j)/\varepsilon -1 \big)\bigg) \,.
    \\
    &=
    \begin{cases}
     \rho - \abs{\mathcal{N}(i^*)}
    \\
    \text{where }
    \\
    \rho_{i^*,j}(x_{i^*})
     = \mathrm{lse}_\varepsilon\big(f_{i^*,j}(x_{i^*}, \cdot)\big) + \tfrac{1}{\abs{\mathcal{N}(i^*)}} \bigg(\rho
     -\!\!\!\!\!\!\!\!\!\!\displaystyle\sum_{(i^*,k)\in \mathcal{N}(i^*)} \!\!\!\!\!\!\!\!\!\mathrm{lse}_\varepsilon\big(f_{i^*,k}(x_{i^*}, \cdot)\big)\bigg)\,,
    \\
     \qquad\qquad\forall x_{i^*} \in \Omega_{i^*} \,,\,\, (i^*, j) \in \mathcal{N}(i^*)
    \\
    \rho = \varepsilon\abs{\mathcal{N}(i^*)} +\varepsilon\abs{\mathcal{N}(i^*)}\log\big(\abs{\mathcal{N}(i^*)}/\varepsilon\big)
    \\
    \qquad
    - \abs{\mathcal{N}(i^*)}\mathrm{lse}_\varepsilon^{x_{i^*}, j}\bigg(\mathrm{lse}^{x_j}_\varepsilon\big(f_{i^*,j}(x_{i^*}, x_j)\big) +\mathrm{lse}^{x_j}_\varepsilon\big(-f_{i^*, j}(x_{i^*}, x_j)\big)\\
    \qquad\qquad\qquad\qquad\qquad
    -\tfrac{1}{\abs{\mathcal{N}(i^*)}}\!\!\!\!\!\!\!\!\!\!\displaystyle\sum_{(i^*,k)\in \mathcal{N}(i^*)} \!\!\!\!\!\!\!\!\!\mathrm{lse}^{x_k}_\varepsilon\big(f_{i^*,k}(x_{i^*}, x_k)\big)\bigg) \,.
    \end{cases}
\end{align*}
\end{lemma}
\begin{proof}
We have
\begingroup
\allowdisplaybreaks
\begin{align*}
     &\displaystyle\max_{\rho_{i^*, \cdot} \in \R^{\Omega_{i^*}}}\bigg(\bigg(\min_{
        x_{i^*}\in  \Omega_{i^*}}\big(\textstyle\sum_{(i^*,j) \in \mathcal{N}(i^*) } \rho_{i^*, j}\big)(x_{i^*}) \bigg)
    \\
    &\qquad\qquad\qquad\qquad
    -\varepsilon\!\!\!\displaystyle\sum_{x_{i^*} \in \Omega_{i^*}} \sum_{(i^*,j) \in \mathcal{N}(i^*) }\exp\big(\rho_{i^*,j}(x_{i^*})/\varepsilon\big)\sum_{x_j \in \Omega_j}
        \exp \big(-f_{i^*,j}(x_{i^*},x_j)/\varepsilon -1 \big)\bigg)
    \\
    &=
    \begin{cases}
    \displaystyle\max_{\rho \in \R} \displaystyle\max_{\rho_{i^*, \cdot} \in \R^{\Omega_{i^*}}} \! \rho  -\varepsilon\!\!\!\!\!\displaystyle\sum_{x_{i^*} \in \Omega_{i^*}} \sum_{(i^*,j) \in \mathcal{N}(i^*) } \!\!\!\!\!\! \exp\big(\rho_{i^*,j}(x_{i^*})/\varepsilon\big)\!\!\!\sum_{x_j \in \Omega_j}\!\!
        \exp \big(\!-f_{i^*,j}(x_{i^*},x_j)/\varepsilon -1 \big)
    \\
    \text{s.t. }\big(\textstyle\sum_{(i^*,j) \in \mathcal{N}(i^*) } \rho_{i^*, j}\big)(x_{i^*}) = \rho \,,\,\, \forall x_{i^*} \in \Omega_{i^*}
    \end{cases}
\tag*{\small{(separability \& monotonicity)}}
    \\
    &=
    \max_{\rho \in \R} \rho - \varepsilon\cdot\!
    \begin{cases}
    \min_{\rho_{i^*, \cdot}\!\! \in \R^{\Omega_{i^*}}}  \!\!\!\!\!\displaystyle\sum_{x_{i^*} \in \Omega_{i^*}} \sum_{(i^*,j) \in \mathcal{N}(i^*) } \!\!\!\!\!\!\exp\big(\rho_{i^*,j}(x_{i^*})/\varepsilon\big)\!\!\!\sum_{x_j \in \Omega_j}\!\!
        \exp \big(\!-f_{i^*,j}(x_{i^*},x_j)/\varepsilon \! -\!1 \big)
    \\
    \text{s.t. }\big(\textstyle\sum_{(i^*,j) \in \mathcal{N}(i^*) } \rho_{i^*, j}\big)(x_{i^*}) = \rho \,,\,\, \forall x_{i^*} \in \Omega_{i^*}
    \end{cases}
\tag*{\small{(separability)}}
    \\
    &=
    \begin{cases}
    \max_{\rho \in \R} \rho - \varepsilon\!\!\!\displaystyle\sum_{x_{i^*} \in \Omega_{i^*}} \sum_{(i^*,j) \in \mathcal{N}(i^*) }\exp\big(\rho_{i^*,j}(x_{i^*})/\varepsilon\big)\sum_{x_j \in \Omega_j}
        \exp \big(-f_{i^*,j}(x_{i^*},x_j)/\varepsilon -1 \big)
    \\
    \text{where }\rho_{i^*,j}(x_{i^*})
     = \mathrm{lse}_\varepsilon\big(f_{i^*,j}(x_{i^*}, \cdot)\big) + \tfrac{1}{\abs{\mathcal{N}(i^*)}} \bigg(\rho
     -\!\!\!\!\!\!\!\!\!\!\displaystyle\sum_{(i^*,k)\in \mathcal{N}(i^*)} \!\!\!\!\!\!\!\!\!\mathrm{lse}_\varepsilon\big(f_{i^*,k}(x_{i^*}, \cdot)\big)\bigg)\,,
     \\
     \qquad\qquad\forall x_{i^*} \in \Omega_{i^*} \,,\,\, (i^*, j) \in \mathcal{N}(i^*)
    \end{cases}
\tag*{\small{(see \enquote{inner problem} below)}}
    \\
    &=
    \begin{cases}
    \rho - \varepsilon\!\!\!\displaystyle\sum_{x_{i^*} \in \Omega_{i^*}} \sum_{(i^*,j) \in \mathcal{N}(i^*) }\exp\big(\rho_{i^*,j}(x_{i^*})/\varepsilon\big)\sum_{x_j \in \Omega_j}
        \exp \big(-f_{i^*,j}(x_{i^*},x_j)/\varepsilon -1 \big)
    \\
    \text{where }
    \\
    \rho_{i^*,j}(x_{i^*})
     = \mathrm{lse}_\varepsilon\big(f_{i^*,j}(x_{i^*}, \cdot)\big) + \tfrac{1}{\abs{\mathcal{N}(i^*)}} \bigg(\rho
     -\!\!\!\!\!\!\!\!\!\!\displaystyle\sum_{(i^*,k)\in \mathcal{N}(i^*)} \!\!\!\!\!\!\!\!\!\mathrm{lse}_\varepsilon\big(f_{i^*,k}(x_{i^*}, \cdot)\big)\bigg)\,,
     \\
     \qquad\qquad\forall x_{i^*} \in \Omega_{i^*} \,,\,\, (i^*, j) \in \mathcal{N}(i^*)
    \\
    \rho = \varepsilon\abs{\mathcal{N}(i^*)} +\varepsilon\abs{\mathcal{N}(i^*)}\log\big(\abs{\mathcal{N}(i^*)}/\varepsilon\big)
    \\
    \qquad
    - \abs{\mathcal{N}(i^*)}\mathrm{lse}_\varepsilon^{x_{i^*}, j}\bigg(\mathrm{lse}^{x_j}_\varepsilon\big(f_{i^*,j}(x_{i^*}, x_j)\big) +\mathrm{lse}^{x_j}_\varepsilon\big(-f_{i^*, j}(x_{i^*}, x_j)\big)\\
    \qquad\qquad\qquad\qquad\qquad
     -\tfrac{1}{\abs{\mathcal{N}(i^*)}}\!\!\!\!\!\!\!\!\!\!\displaystyle\sum_{(i^*,k)\in \mathcal{N}(i^*)} \!\!\!\!\!\!\!\!\!\mathrm{lse}^{x_k}_\varepsilon\big(f_{i^*,k}(x_{i^*}, x_k)\big)\bigg)
    \end{cases}
\tag*{\small{(see \enquote{outer problem} below)}}
    \\
    &=
    \begin{cases}
     \rho - \abs{\mathcal{N}(i^*)}
    \\
    \text{where }
    \\
    \rho_{i^*,j}(x_{i^*})
     = \mathrm{lse}_\varepsilon\big(f_{i^*,j}(x_{i^*}, \cdot)\big) + \tfrac{1}{\abs{\mathcal{N}(i^*)}} \bigg(\rho
     -\!\!\!\!\!\!\!\!\!\!\displaystyle\sum_{(i^*,k)\in \mathcal{N}(i^*)} \!\!\!\!\!\!\!\!\!\mathrm{lse}_\varepsilon\big(f_{i^*,k}(x_{i^*}, \cdot)\big)\bigg)\,,
     \\
     \qquad\qquad\forall x_{i^*} \in \Omega_{i^*} \,,\,\, (i^*, j) \in \mathcal{N}(i^*)
    \\
    \rho = \varepsilon\abs{\mathcal{N}(i^*)} +\varepsilon\abs{\mathcal{N}(i^*)}\log\big(\abs{\mathcal{N}(i^*)}/\varepsilon\big)
    \\
    \qquad
    - \abs{\mathcal{N}(i^*)}\mathrm{lse}_\varepsilon^{x_{i^*}, j}\bigg(\mathrm{lse}^{x_j}_\varepsilon\big(f_{i^*,j}(x_{i^*}, x_j)\big) +\mathrm{lse}^{x_j}_\varepsilon\big(-f_{i^*, j}(x_{i^*}, x_j)\big)\\
    \qquad\qquad\qquad\qquad\qquad
    -\tfrac{1}{\abs{\mathcal{N}(i^*)}}\!\!\!\!\!\!\!\!\!\!\displaystyle\sum_{(i^*,k)\in \mathcal{N}(i^*)} \!\!\!\!\!\!\!\!\!\mathrm{lse}^{x_k}_\varepsilon\big(f_{i^*,k}(x_{i^*}, x_k)\big)\bigg) \,.
    \end{cases}
\tag*{\small{(see \enquote{simplify maximal value} below)}}
\end{align*}%
\endgroup%

\underline{Simplify maximal value:} In the above setting, we have
\begingroup
\allowdisplaybreaks
\begin{align*}%
    &
    \rho - \varepsilon\!\!\!\displaystyle\sum_{x_{i^*} \in \Omega_{i^*}} \sum_{(i^*,j) \in \mathcal{N}(i^*) }\exp\big(\rho_{i^*,j}(x_{i^*})/\varepsilon\big)\sum_{x_j \in \Omega_j}
        \exp \big(-f_{i^*,j}(x_{i^*},x_j)/\varepsilon -1 \big)
    \\
    &=
    \rho - \varepsilon\exp\tfrac{1}{\varepsilon}\varepsilon\log\displaystyle\sum_{x_{i^*} \in \Omega_{i^*}} \sum_{(i^*,j) \in \mathcal{N}(i^*) }\exp\bigg(\big(\rho_{i^*,j}(x_{i^*})+ \mathrm{lse}^{x_j}_\varepsilon\big(-f_{i^*,j}(x_{i^*},x_j)- \varepsilon\big)\big)/\varepsilon\bigg)
    \tag*{\small{(definition $\mathrm{lse}_\varepsilon$; properties of $\exp$; $\exp \log \equiv \mathrm{id}$)}}
    \\
    &=
    \rho - \varepsilon\exp\bigg(\mathrm{lse}_\varepsilon^{x_{i^*},j}\bigg(\rho_{i^*,j}(x_{i^*})+ \mathrm{lse}^{x_j}_\varepsilon\big(-f_{i^*,j}(x_{i^*},x_j)\big)\bigg)/\varepsilon-1\bigg) \,.
    \tag*{\small{(definition $\mathrm{lse}_\varepsilon$; properties of $\mathrm{lse}_\varepsilon$)}}
    \\
    &=
    \rho - \varepsilon\exp\bigg(\mathrm{lse}_\varepsilon^{x_{i^*},j}\bigg(\mathrm{lse}^{x_j}_\varepsilon\big(f_{i^*,j}(x_{i^*}, x_j)\big) + \tfrac{1}{\abs{\mathcal{N}(i^*)}} \bigg(\rho
     -\!\!\!\!\!\!\!\!\!\!\displaystyle\sum_{(i^*,k)\in \mathcal{N}(i^*)} \!\!\!\!\!\!\!\!\!\mathrm{lse}^{x_k}_\varepsilon\big(f_{i^*,k}(x_{i^*}, x_k)\big)\bigg)\\
     &\qquad\qquad\qquad\qquad\qquad
     +\mathrm{lse}^{x_j}_\varepsilon\big(-f_{i^*,j}(x_{i^*},x_j)\big)\bigg)/\varepsilon-1\bigg) \,.
    \tag*{\small{(representation of $\rho_{i^*,j}(x_{i^*})$; renaming variables)}}
    \\
    &=
    \rho - \varepsilon\exp\bigg(\mathrm{lse}_\varepsilon^{x_{i^*},j}\bigg(\mathrm{lse}^{x_j}_\varepsilon\big(f_{i^*,j}(x_{i^*}, x_j)\big) - \tfrac{1}{\abs{\mathcal{N}(i^*)}}\!\!\!\!\!\!\!\!\!\!\displaystyle\sum_{(i^*,k)\in \mathcal{N}(i^*)} \!\!\!\!\!\!\!\!\!\mathrm{lse}^{x_k}_\varepsilon\big(f_{i^*,k}(x_{i^*}, x_k)\big)\\
     &\qquad\qquad\qquad\qquad\qquad
     +\mathrm{lse}^{x_j}_\varepsilon\big(-f_{i^*,j}(x_{i^*},x_j)\big)\bigg)/\varepsilon +  \tfrac{\rho}{\varepsilon\abs{\mathcal{N}(i^*)}} -1\bigg)
    \tag*{\small{(constant addition to $\mathrm{lse}_\varepsilon$)}}
    \\
    &=
    \rho - \varepsilon\exp\bigg(1 + \log\big(\abs{\mathcal{N}(i^*)}/\varepsilon\big)  -1\bigg)
    \tag*{\small{(representation of $\rho$)}}
    \\
    &=
    \rho - \abs{\mathcal{N}(i^*)} \,.
    \tag*{\small{(simplifying)}}
\end{align*}%
\endgroup%

\underline{Outer problem:} Since the outer problem is unbounded from below at infinity, a necessary condition for the maximizers of the outer problem can be found using a first order criterion, i.e.~in critical points. Define
\begin{align*}
    R(\rho) := \rho - \varepsilon\cdot\displaystyle\sum_{x_{i^*} \in \Omega_{i^*}} \sum_{(i^*,j) \in \mathcal{N}(i^*) }\exp\big(\rho_{i^*,j}(x_{i^*})/\varepsilon\big)\sum_{x_j \in \Omega_j}
        \exp \big(-f_{i^*,j}(x_{i^*},x_j)/\varepsilon -1 \big)
    \\
    \text{where }
    \rho_{i^*,j}(x_{i^*})
     = \underbrace{\mathrm{lse}_\varepsilon\big(f_{i^*,j}(x_{i^*}, \cdot)\big) -\tfrac{1}{\abs{\mathcal{N}(i^*)}}\!\!\!\!\!\!\!\!\!\!\displaystyle\sum_{(i^*,k)\in \mathcal{N}(i^*)} \!\!\!\!\!\!\!\!\!\mathrm{lse}_\varepsilon\big(f_{i^*,k}(x_{i^*}, \cdot)\big)}_{=: \gamma(x_{i^*},j)}+ \tfrac{\rho}{\abs{\mathcal{N}(i^*)}}
      \,.
\end{align*}
We have
\begingroup
\allowdisplaybreaks
\begin{align*}
    0&=\frac{\mathrm{d}}{\mathrm{d}\rho} R(\rho)
    \\
    &= 1\!-\! \tfrac{\varepsilon}{\abs{\mathcal{N}(i^*)}}\cdot\!\!\!\displaystyle\sum_{x_{i^*} \in \Omega_{i^*}} \sum_{(i^*,j) \in \mathcal{N}(i^*) }\!\!\!\!\!\!\!\exp\big(\rho_{i^*,j}(x_{i^*})/\varepsilon\big)\!\!\sum_{x_j \in \Omega_j}\!\!\!
        \exp \big(\!-\!f_{i^*,j}(x_{i^*},x_j)/\varepsilon \!-\!1 \big)
    \\
    &= 1\!-\! \tfrac{\varepsilon}{\abs{\mathcal{N}(i^*)}}\cdot\!\!\!\displaystyle\sum_{x_{i^*} \in \Omega_{i^*}} \sum_{(i^*,j) \in \mathcal{N}(i^*) }\!\!\!\!\!\!\!\exp\big(\gamma(x_{i^*},j)/\varepsilon + \tfrac{\rho}{\varepsilon\abs{\mathcal{N}(i^*)}}\big)\!\!\sum_{x_j \in \Omega_j}\!
        \!\!\exp \big(\!-\!f_{i^*,j}(x_{i^*},x_j)/\varepsilon \!-\!1 \big)
\tag*{\small{(definition of $\rho_{i^*,\cdot}$)}}
    \\
    \iff&
    \\
    1= &\exp\big(\tfrac{\rho}{\varepsilon\abs{\mathcal{N}(i^*)}}\big)\tfrac{\varepsilon}{\abs{\mathcal{N}(i^*)}}\cdot\!\!\!\!\displaystyle\sum_{x_{i^*} \in \Omega_{i^*}} \sum_{(i^*,j) \in \mathcal{N}(i^*) }\!\!\!\!\exp\big(\gamma(x_{i^*},j)/\varepsilon\big)\!\!\sum_{x_j \in \Omega_j}\!\!\!
        \exp \big(\!-\!f_{i^*,j}(x_{i^*},x_j)/\varepsilon -1 \big)
\tag*{\small{(rearranging)}}
    \\
    \iff&
    \rho
    =
    \varepsilon\abs{\mathcal{N}(i^*)}\log\bigg(\frac{\abs{\mathcal{N}(i^*)}}{\varepsilon\!\!\!\!\!\!\displaystyle\sum_{x_{i^*} \in \Omega_{i^*}} \sum_{(i^*,j) \in \mathcal{N}(i^*) }\!\!\!\!\!\!\exp\big(\gamma(x_{i^*},j)/\varepsilon\big)\!\!\!\sum_{x_j \in \Omega_j}\!\!\!
        \exp \big(-f_{i^*,j}(x_{i^*},x_j)/\varepsilon -1 \big)}\bigg)
\tag*{\small{(rearranging)}}
    \\
    \phantom{\iff}&
    \phantom{\rho}
    =
    \varepsilon\abs{\mathcal{N}(i^*)}\log\bigg(\frac{\abs{\mathcal{N}(i^*)}}{\varepsilon\!\!\!\!\!\!\displaystyle\sum_{x_{i^*} \in \Omega_{i^*}} \sum_{(i^*,j) \in \mathcal{N}(i^*) }\!\!\!\!\!\!\exp\big(\gamma(x_{i^*},j)/\varepsilon\big)\exp \big(\mathrm{lse}_\varepsilon\big(-f_{i^*, j}(x_{i^*}, \cdot)-\varepsilon\big)/\varepsilon\big)}\bigg)
\tag*{\small{(definition $\mathrm{lse}_\varepsilon$)}}
    \\
    \phantom{\iff}&
    \phantom{\rho}
    =
    \varepsilon\abs{\mathcal{N}(i^*)}\log\bigg(\frac{\abs{\mathcal{N}(i^*)}}{\varepsilon\!\!\!\!\!\!\displaystyle\sum_{x_{i^*} \in \Omega_{i^*}} \sum_{(i^*,j) \in \mathcal{N}(i^*) }\!\!\!\!\!\!\exp\bigg(\tfrac{1}{\varepsilon}\big(\gamma(x_{i^*},j)+\mathrm{lse}_\varepsilon\big(-f_{i^*, j}(x_{i^*}, \cdot)-\varepsilon\big)\big)\bigg)}\bigg)
\tag*{\small{(property of $\exp$)}}
    \\
    \phantom{\iff}&
    \phantom{\rho}
    =
    \varepsilon\abs{\mathcal{N}(i^*)}\log\big(\abs{\mathcal{N}(i^*)}/\varepsilon\big)
    \\
    &\quad\qquad- \varepsilon\abs{\mathcal{N}(i^*)}\log\bigg(\displaystyle\sum_{x_{i^*} \in \Omega_{i^*}} \sum_{(i^*,j) \in \mathcal{N}(i^*) }\!\!\!\!\!\!\exp\bigg(\tfrac{1}{\varepsilon}\big(\gamma(x_{i^*},j)+\mathrm{lse}_\varepsilon\big(-f_{i^*, j}(x_{i^*}, \cdot)-\varepsilon\big)\big)\bigg)\bigg)
\tag*{\small{(property of $\log$)}}
    \\
    \phantom{\iff}&
    \phantom{\rho}
    =
    \varepsilon\abs{\mathcal{N}(i^*)}\log\big(\abs{\mathcal{N}(i^*)}/\varepsilon\big)
    - \abs{\mathcal{N}(i^*)}\mathrm{lse}_\varepsilon^{x_{i^*}, j}\bigg(\gamma(x_{i^*},j)+\mathrm{lse}^{x_j}_\varepsilon\big(-f_{i^*, j}(x_{i^*}, x_j)-\varepsilon\big)\bigg)
\tag*{\small{(definition of $\mathrm{lse}_\varepsilon$; superscript of $\mathrm{lse}_\varepsilon$ denotes the variable being maximized)}}
    \\
    \phantom{\iff}&
    \phantom{\rho}
    =
    \varepsilon\abs{\mathcal{N}(i^*)} +\varepsilon\abs{\mathcal{N}(i^*)}\log\big(\abs{\mathcal{N}(i^*)}/\varepsilon\big)
    \\
    &\qquad\qquad\qquad\qquad
    - \abs{\mathcal{N}(i^*)}\mathrm{lse}_\varepsilon^{x_{i^*}, j}\bigg(\gamma(x_{i^*},j)+\mathrm{lse}^{x_j}_\varepsilon\big(-f_{i^*, j}(x_{i^*}, x_j)\big)\bigg)
\tag*{\small{(constant addition to $\mathrm{lse}_\varepsilon$)}}
    \\
    \phantom{\iff}&
    \phantom{\rho}
    =
    \varepsilon\abs{\mathcal{N}(i^*)} +\varepsilon\abs{\mathcal{N}(i^*)}\log\big(\abs{\mathcal{N}(i^*)}/\varepsilon\big)
    \\
    &\qquad
    - \abs{\mathcal{N}(i^*)}\mathrm{lse}_\varepsilon^{x_{i^*}, j}\bigg(\mathrm{lse}^{x_j}_\varepsilon\big(f_{i^*,j}(x_{i^*}, x_j)\big) 
    \\
    &\qquad\qquad\qquad\qquad\qquad\quad
    -\tfrac{1}{\abs{\mathcal{N}(i^*)}}\!\!\!\!\!\!\!\!\!\!\displaystyle\sum_{(i^*,k)\in \mathcal{N}(i^*)} \!\!\!\!\!\!\!\!\!\mathrm{lse}^{x_k}_\varepsilon\big(f_{i^*,k}(x_{i^*}, x_k)\big)+\mathrm{lse}^{x_j}_\varepsilon\big(-f_{i^*, j}(x_{i^*}, x_j)\big)\bigg)
\tag*{\small{(definition $\gamma(x_{i^*}, j)$; superscript of $\mathrm{lse}_\varepsilon$ denotes the variable being maximized)}}
\end{align*}%
\endgroup%

\underline{Inner problem:} Since the inner problem is unbounded from above at infinity, a necessary condition for the minimizers of the inner problem can be found using Lagrange multipliers, i.e.~in critical points of the function
\begin{align*}
    &G\big(\lambda_\cdot, \rho_{i^*, \cdot}(\cdot) \big) :=\displaystyle\sum_{x_{i^*} \in \Omega_{i^*}} \sum_{(i^*,j) \in \mathcal{N}(i^*) }\exp\big(\rho_{i^*,j}(x_{i^*})/\varepsilon\big)\sum_{x_j \in \Omega_j}
        \exp \big(-f_{i^*,j}(x_{i^*},x_j)/\varepsilon -1 \big)
    \\
    &\qquad\qquad\qquad\qquad- \sum_{x_{i^*} \in \Omega_{i^*}} \lambda_{x_{i^*}} \bigg( \big(\textstyle\sum_{(i^*,j) \in \mathcal{N}(i^*) } \rho_{i^*, j}\big)(x_{i^*}) - \rho  \bigg) \,,
\end{align*}
where $\lambda_\cdot \in \R^{\Omega_{i^*}}$.
We have for a critical point that
\begin{align*}
    0
    &= \frac{\partial}{\partial \rho_{i^*, j}(x_{i^*})} G\big(\lambda_\cdot, \rho_{i^*, \cdot}(\cdot) \big)
    \\
    &= (1/\varepsilon)\exp\big(\rho_{i^*,j}(x_{i^*})/\varepsilon\big)\bigg(\sum_{x_j \in \Omega_j}
        \exp \big(-f_{i^*,j}(x_{i^*},x_j)/\varepsilon -1 \big)\bigg) - \lambda_{x_{i^*}}
    \\
    \iff&
    \varepsilon\log\varepsilon\lambda_{x_{i^*}}
    \\
    &\,\,=\!
    \rho_{i^*,j}(x_{i^*})\!+\!\varepsilon\log\bigg(\!\sum_{x_j \in \Omega_j}\!\!
        \exp \big(\!-\!f_{i^*,j}(x_{i^*},x_j)/\varepsilon \!-\!1 \big)\!\bigg)\,,
    \,\, \forall (i^*, j) \in \mathcal{N}(i^*)\,,\,\, \forall x_{i^*} \in \Omega_{i^*} \,.
\tag*{\small{(applying the function $\varepsilon\log (\varepsilon(\cdot+\lambda_{x_{i^*}}))$)}}
\end{align*}
From the constraint of the inner problem, we can find the minimizer, as
\begingroup
\allowdisplaybreaks
\begin{align*}
    &\qquad\big(\textstyle\sum_{(i^*,k) \in \mathcal{N}(i^*) } \rho_{i^*, k}\big)(x_{i^*}) = \rho\,,\,\, \forall x_{i^*} \in \Omega_{i^*} 
    \\
    &\iff
    \sum_{(i^*,k)\in \mathcal{N}(i^*)}\!\!\!\!\!\!\!\!\varepsilon\log\varepsilon\lambda_{x_{i^*}}
    -\varepsilon\log\bigg(\sum_{x_k \in \Omega_k}
        \exp \big(-f_{i^*,k}(x_{i^*},x_k)/\varepsilon -1 \big)\bigg) = \rho\,,\,\, \forall x_{i^*} \in \Omega_{i^*} 
\tag*{\small{(necessary condition on critical points)}}
    \\
    &\iff
    \abs{\mathcal{N}(i^*)}\varepsilon\log\varepsilon\lambda_{x_{i^*}}
    -\!\!\!\!\!\!\!\sum_{(i^*,k)\in \mathcal{N}(i^*)}\!\!\!\!\!\!\!\!\!\varepsilon\log\bigg(\sum_{x_k \in \Omega_k}\!\!
        \exp \big(\!-\!f_{i^*,k}(x_{i^*},x_k)/\varepsilon \!-\!1 \big)\bigg) = \rho\,,\,\, \forall x_{i^*} \in \Omega_{i^*} 
\tag*{\small{(constant summation)}}
    \\
    &\iff
    \abs{\mathcal{N}(i^*)}\bigg( \rho_{i^*,j}(x_{i^*})+\varepsilon\log\bigg(\sum_{x_j \in \Omega_j}
        \exp \big(-f_{i^*,j}(x_{i^*},x_j)/\varepsilon -1 \big)\bigg)\bigg)
    \\
    &\qquad\qquad
    -\!\!\!\!\!\!\!\!\!\!\sum_{(i^*,k)\in \mathcal{N}(i^*)}\!\!\!\!\!\!\!\!\!\varepsilon\log\bigg(\sum_{x_k \in \Omega_k}\!\!
        \exp \big(\!-\!f_{i^*,k}(x_{i^*},x_k)/\varepsilon\! -\!1 \big)\bigg) = \rho\,,\,\, \forall x_{i^*} \in \Omega_{i^*} \,,\,\, (i^*, j) \in \mathcal{N}(i^*)
\tag*{\small{(necessary condition on critical points)}}
    \\
    &\iff
    \rho_{i^*,j}(x_{i^*})
     = \tfrac{1}{\abs{\mathcal{N}(i^*)}} \bigg(\rho
     +\!\!\!\!\!\!\sum_{(i^*,k)\in \mathcal{N}(i^*)}\!\!\!\!\!\!\!\!\varepsilon\log\bigg(\sum_{x_k \in \Omega_k}
        \exp \big(-f_{i^*,k}(x_{i^*},x_k)/\varepsilon -1 \big)\bigg)\bigg)
    \\
    &\qquad\qquad\qquad\quad
    -
    \varepsilon\log\bigg(\sum_{x_j \in \Omega_j}
        \exp \big(-f_{i^*,j}(x_{i^*},x_j)/\varepsilon -1 \big)\bigg)
     \,,\,\, \forall x_{i^*} \in \Omega_{i^*} \,,\,\, (i^*, j) \in \mathcal{N}(i^*)
\tag*{\small{(rearranging)}}
    \\
    &\iff
    \rho_{i^*,j}(x_{i^*})
    \\
    &\quad\,\,\,
     = \tfrac{1}{\abs{\mathcal{N}(i^*)}} \bigg(\!\rho
     \!+\!\varepsilon\!\!\!\!\!\!\!\!\!\!\!\sum_{(i^*,k)\in \mathcal{N}(i^*)}\!\!\!\!\!\!\!\!\!\!\log\!\bigg(\!\frac{\sum_{x_j \in \Omega_j}\!\!
        \exp \big(f_{i^*,j}(x_{i^*},x_j)/\varepsilon \cancel{+1} \big)}{\sum_{x_k \in \Omega_k}\!\!
        \exp \big(f_{i^*,k}(x_{i^*},x_k)/\varepsilon \cancel{+1} \big)}\bigg)\!\bigg),\, \forall x_{i^*} \in \Omega_{i^*} ,\, (i^*, j) \in \mathcal{N}(i^*) .   
\tag*{\small{(simplifying)}}
    \\
    &\iff
    \rho_{i^*,j}(x_{i^*})
    \\
    &\qquad
     = \mathrm{lse}_\varepsilon\big(f_{i^*,j}(x_{i^*}, \cdot)\big) + \tfrac{1}{\abs{\mathcal{N}(i^*)}} \bigg(\rho
     -\!\!\!\!\!\!\!\!\!\sum_{(i^*,k)\in \mathcal{N}(i^*)} \!\!\!\!\mathrm{lse}_\varepsilon\big(f_{i^*,k}(x_{i^*}, \cdot)\big)\bigg),\, \forall x_{i^*} \in \Omega_{i^*} ,\, (i^*, j) \in \mathcal{N}(i^*) .   
\tag*{\small{(simplifying)}}
\end{align*}%
\endgroup%
\end{proof}%

\clearpage
\section{Additional Code}

\subsection{{Verifying that $F\circ\Tilde{p}$ of \cref{ex:freelunch} is not a sum of bivariates}}
\label{code:freelunch}
\input{mintedcache/codesnip.pygtex}

\clearpage
\subsection{Pseudocode for TRW-S}

\begin{algorithm}[ht]
\caption{Legacy Sequential Tree-Reweighted Message Passing (\texttt{TRW-S-LEG}) by \cite{kolmogorov2005}%
}
\begin{algorithmic}[1]
    \INPUT vertex set $\mathcal{V} := \N_{\leq n}$, $n \in \N$;
    \SUBINPUT edge set $\mathcal{E} \subseteq \{(i,j) \in \mathcal{V} \times \mathcal{V} \mid i < j\}$;
    \SUBINPUT candidate sets $\Omega_i \subset \R$, where $\abs{\Omega_i} \in \N$ and $i \in \mathcal{V}$;
    \SUBINPUT functions $f_{i,j}: \Omega_i \times \Omega_j \to \R$, where $(i,j) \in \mathcal{E}$;
    \SUBINPUT bijection $c: \mathcal{V} \to \N_{\leq n}$;
    \SUBINPUT budget $B \in \N$.
    \OUTPUT Element $x^* \in \Omega_1 \times \dots \times \Omega_n$ with \enquote{low} function value $\sum_{(i,j) \in \mathcal{E}} f_{i,j}(x_i, x_j)$.
    \INIT $t,\delta := 1,0 \in \R$;\quad $\Tilde{\mathcal{E}} := \big\{ \{ i,j\} \mid (i,j) \in \mathcal{E}\}$;
    \SUBINPUT \quad\,\, $m_{i,j} := 0 \in \R^{\Omega_i}, m_{j,i} := 0 \in \R^{\Omega_j}, (i,j) \in \mathcal{E}$;\quad $m_i := 0 \in \R^{\Omega_i}, i \in \mathcal{V}$;
    \SUBINPUT \quad\,\, $p_i := \max\big\{ \abs{\{\{i,j\} \in \Tilde{\mathcal{E}} : c(i) < c(j) \}}, \abs{\{\{i,j\} \in \Tilde{\mathcal{E}} : c(j) < c(i)\}}\big\}, i \in \mathcal{V}$;
    \SUBINPUT \quad\,\, $\gamma_{i,j} := 1/p_i, \gamma_{j,i} := 1/p_j, (i,j) \in \mathcal{E}$.
    \AlgBlankLine
    \FORALL{$i \in (c^{-1}(1), \dots, c^{-1}(n))$}           \algcomment{loop vertices in order defined by $c$}
        \STATE $m_i \leftarrow \sum_{\{i,j\} \in \Tilde{\mathcal{E}}} m_{i,j}$
        \STATE $\delta \leftarrow \min_{x_i \in \Omega_i} m_i(x_i)$
        \STATE $m_i(x_i) \leftarrow m_i(x_i) - \delta \,, \,\, \forall x_i\in \Omega_i$
        \FORALL{$\{i,j\} \in \Tilde{\mathcal{E}}$ with $c(i) < c(j)$}     \algcomment{restrict to \enquote{increasing} edges}
            \STATE \!\!\!$m_{j,i}(x_j) \!\!\leftarrow \!\!\min_{x_i \in \Omega_i} \!\gamma_{i,j} m_i(x_i)\! - \!m_{i,j}(x_i) \!+\! \!f_{\min\{i,j\},\max\{i,j\}}\!(x_{\min\{i,j\}},\!x_{\max\{i,j\}}\!),\!\forall x_j\! \in\! \Omega_j$
            \STATE \!\!\!$\delta \leftarrow \min_{x_j \in \Omega_j} m_{j,i}(x_j)$
            \STATE \!\!\!$m_{j,i}(x_j) \leftarrow m_{j,i}(x_j) - \delta, \forall x_j \in \Omega_j$
        \ENDFOR
    \ENDFOR
    \FOR{$i \in (c^{-1}(1), \dots, c^{-1}(n))$}
         \algcomment{determine the solution candidate}
        \STATE $y^*_i \!\!\in \!\argmin_{x_i \in \Omega_i} \!\sum_{\!\!\!\substack{\{i,j\} \in \widetilde{\mathcal{E}}\\c(j)<c(i)}} \!\!\! f_{\min\{i,j\},\max\{i,j\}}(x_{\min\{i,j\}}, x_{\max\{i,j\}}) + \!\sum_{\substack{\{i,j\}\in\widetilde{\mathcal{E}}\\c(i)<c(j)}} \!\!m_{i,j}(x_i)$
    \ENDFOR
    \IF{$F(y^*)<F(x^*)$} \algcomment{keep best solution candidate}
        \STATE $x^* \leftarrow y^*$
    \ENDIF
    \IF{$B\leq t$} \algcomment{stopping criterion}
        \STATE \texttt{return} $x^*$
    \ENDIF
    \STATE $t \leftarrow t+1$         \algcomment{run again in reversed order}
    \STATE $c(i) \leftarrow \abs{\mathcal{V}}+1-i \,,\,\, i\in\mathcal{V}$
    \STATE \texttt{go to line 1}
  \end{algorithmic}
\label{alg:trw-s-leg}
\end{algorithm}

\begin{algorithm}[ht]
\caption{Sequential Tree-Reweighted Message Passing (\texttt{TRW-S}) by \cite{tourani2020}%
}
\begin{algorithmic}[1]
    \INPUT vertex set $\mathcal{V} := \N_{\leq n}$, $n \in \N$;
    \SUBINPUT edge set $\mathcal{E} \subseteq \{(i,j) \in \mathcal{V} \times \mathcal{V} \mid i < j\}$;
    \SUBINPUT candidate sets $\Omega_i \subset \R$, where $\abs{\Omega_i} \in \N$ and $i \in \mathcal{V}$;
    \SUBINPUT functions $f_{i,j}: \Omega_i \times \Omega_j \to \R$, where $(i,j) \in \mathcal{E}$;
    \SUBINPUT bijection $c: \mathcal{V} \to \N_{\leq n}$;
    \SUBINPUT budget $B \in \N$.
    \OUTPUT Element $x^* \in \Omega_1 \times \dots \times \Omega_n$ with \enquote{low} function value $\sum_{(i,j) \in \mathcal{E}} f_{i,j}(x_i, x_j)$.
    \INIT $t := 1 \in \N$;\quad $\widetilde{\mathcal{E}} := \big\{ \{ i,j\} \mid (i,j) \in \mathcal{E}\}$;
    \SUBINPUT \quad\,\, $m_{i,j},\rho_{i,j} := 0 \in \R^{\Omega_i}, m_{j,i},\rho_{j,i} := 0 \in \R^{\Omega_j}, (i,j) \in \mathcal{E}$;\quad $m_i := 0 \in \R^{\Omega_i}, i \in \mathcal{V}$;
    \SUBINPUT \quad\,\, $p_i := \max\big\{ \abs{\{\{i,j\} \in \widetilde{\mathcal{E}} : c(i) < c(j) \}}, \abs{\{\{i,j\} \in \widetilde{\mathcal{E}} : c(j) < c(i)\}}\big\}, i \in \mathcal{V}$;
    \SUBINPUT \quad\,\, $\gamma_{i,j} := 1/p_i, \gamma_{j,i} := 1/p_j, (i,j) \in \mathcal{E}$.
    \AlgBlankLine
    \FORALL{$i \in (c^{-1}(1), \dots, c^{-1}(n))$}           \algcomment{loop vertices in order defined by $c$}
        \STATE $m_i \leftarrow \sum_{\{i,j\} \in \widetilde{\mathcal{E}}} m_{i,j}$
        \FORALL{$\{i,j\} \in \widetilde{\mathcal{E}}$ with $c(i) < c(j)$}     \algcomment{restrict to \enquote{increasing} edges}
            \STATE $m_{i,j}(x_i) \leftarrow \min_{x_j \in \Omega_j} f_{\min\{i,j\},\max\{i,j\}}(x_{\min\{i,j\}}, x_{\max\{i,j\}}) - \rho_{j,i}(x_j), \forall x_i \in \Omega_i$
            \STATE $\rho_{i,j} \leftarrow m_{i,j} - \gamma_{i,j} m_i$
        \ENDFOR
    \ENDFOR
    \FOR{$i \in (c^{-1}(1), \dots, c^{-1}(n))$}
         \algcomment{determine the solution candidate}
        \STATE $y^*_i \!\!\in \!\argmin_{x_i \in \Omega_i} \!\sum_{\!\!\!\substack{\{i,j\} \in \widetilde{\mathcal{E}}\\c(j)<c(i)}} \!\!\! f_{\min\{i,j\},\max\{i,j\}}(x_{\min\{i,j\}}, x_{\max\{i,j\}}) + \!\sum_{\substack{\{i,j\}\in\widetilde{\mathcal{E}}\\c(i)<c(j)}} \!\!m_{i,j}(x_i)$
    \ENDFOR
    \IF{$F(y^*)<F(x^*)$} \algcomment{keep best solution candidate}
        \STATE $x^* \leftarrow y^*$
    \ENDIF
    \IF{$B\leq t$} \algcomment{stopping criterion}
        \STATE \texttt{return} $x^*$
    \ENDIF
    \STATE $t \leftarrow t+1$         \algcomment{run again in reversed order}
    \STATE $c(i) \leftarrow \abs{\mathcal{V}}+1-i \,,\,\, i\in\mathcal{V}$
    \STATE \texttt{go to line 1}
  \end{algorithmic}
\label{alg:trw-s}
\end{algorithm}

\clearpage
\section{Referenced Results}

\begin{theorem}[{Convergence~of~Block~Coordinate~Descent \cite[p.~324~Prop.~3.7.1]{bertsekas2018nonlinear}}]%
\label{thm:bertsekas}%
    Let
    \begin{itemize}
        \item a set $X = X_1 \times X_2 \times \dots \times X_m$ be given, where $m, n_i \in \N$ and $X_i \subseteq \R^{n_i}$ is closed and convex for all $i \in \N_{\leq m}$,
        \item the function $f: X \to \R$ be continuously differentiable, and
        \item the function $\xi \in X_i \mapsto f(x_1, \dots, x_{i-1}, \xi, x_{i+1}, \dots, x_m)$ attain a unique minimum $\overline{\xi}$ for all $x \in X$ and all $i \in \N_{\leq m}$, and be monotonically non-increasing in the interval from $x_i$ to $\overline{\xi}$.
    \end{itemize}
    Then, every limit point of the sequence $(x_k)_{k \in \N}$ in $X$, where $x_0 \in X$, and where for all $k \in \N$, we have
    \[
        x_i^{k+1} \in \argmin_{\xi \in X_i} f(x_1^{k+1}, \dots, x_{i-1}^{k+1}, \xi, x_{i+1}^k, \dots, x_m^k) \,, \quad i=1, \dots, m \,,
    \]
    is a stationary point of $f$.
\end{theorem}

\begin{theorem}[{Strong~Lagrangian~Duality~\cite{slater}~\&~\cite[p.~322~Thm.~6.13]{geiger2002}}]%
\label{thm:slater}%
    Consider 
    \begin{itemize}
        \item the convex functions $f: \R^n \to \R$ and $g_i: \R^n \to \R, i \in \N_{\leq m}$, as well as
        \item the affine functions $h_i:\R^n \to \R, i \in \N_{\leq p}$, where $n,m,p \in \N$, and
        \item the non-empty convex set $X \subseteq \R^n$.
    \end{itemize}
    Then, if there exists $\hat{x}$ in the relative interior of $X$ and, such that,
    \[
     g_i(\hat{x}) < 0 \quad \text{for}\,\, i = 1, \dots, m \quad \text{and} \quad h(\hat{x}) = 0\,,
    \]
    we have the equality
    \[
        \begin{cases}
            \inf_{x \in X} f(x)
            \\
            \,\,\text{s.t.}\quad
            g(x) \leq 0
            \\
            \phantom{\,\,\text{s.t.}\quad} h(x) = 0
        \end{cases}
        =
        \begin{cases}
            \sup_{\substack{\lambda \in \R^m\\ \mu \in \R^p}} \inf_{x \in X} f(x) + \sum_{i=1}^m \lambda_i g_i(x) + \sum_{j=1}^p \mu_j h_j(x)
            \\
            \,\,\text{s.t.}\quad
            \lambda \geq 0 \,,
        \end{cases}
    \]
    if the left-hand side of it is finite.
\end{theorem}

\begin{theorem}[{Consistency~of~Global~Relaxations~\cite[Corollary~3.1.1]{mueller2021}}]%
\label{thm:globalrelaxation}
    Given
    \begin{itemize}
        \item a family of probability measures $\{ \Pr_\theta : \theta \in \Theta\}$ on a measure space $(\Omega, \mathcal{A})$,
        \item an optimization problem $\min_{x \in \Omega} f(x)$, where $f: \Omega \to \R$ is $\mathcal{A}$-continuous at its unique global minimum $x^* \in \Omega$,
        \item $f\in \mathcal{L}^1(\Omega, \mathcal{A}, \Pr_\theta)$ for all $\theta \in \Theta$, and
        \item that for all $\varepsilon, \gamma >0$, we have
        \[
            \int_{\Omega-U_\gamma(x^*)} \max\{\abs{f}, 1\} \, \mathrm{d} \Pr_\theta < \varepsilon \,,
        \]
        where $U_\gamma$ denotes the $\gamma$-ball centered at $x$.
        \end{itemize}
        Then, if a minimum $\theta^*$ of
        \[
            \theta \in \Theta \mapsto \int_\Omega f \, \mathrm{d}\Pr_\theta
        \]
        exists, it is unique and $\Pr_{\theta^*} = \delta_{x^*}$, i.e.~the Dirac measure at $x^*$.
\end{theorem}
\end{document}